\documentclass[11pt]{article}
\usepackage[dvipdfmx]{color}
\usepackage[dvipdfmx]{graphicx}
\usepackage[top=30truemm,bottom=30truemm,left=20truemm,right=20truemm]{geometry}
\usepackage{amsmath}
\usepackage{amssymb}
\usepackage{amsthm}
\usepackage{ascmac}
\usepackage{bm}
\usepackage{braket}
\usepackage{color}
\usepackage{comment}
\usepackage{here}

\theoremstyle{definition}
\newtheorem{theorem}{Theorem}[section]
\newtheorem*{theorem*}{Theorem}
\newtheorem{definition}[theorem]{Definition}
\newtheorem*{definition*}{Definition}
\newtheorem{prop}[theorem]{Proposition}
\newtheorem*{prop*}{Proposition}
\newtheorem{lemma}[theorem]{Lemma}
\newtheorem*{lemma*}{Lemma}

\newtheorem*{cor*}{Corollary}
\newtheorem{remark}[theorem]{Remark}
\newtheorem*{remark*}{Remark}

\newtheorem*{example*}{Example}

\newtheorem*{fact*}{Fact}

\newtheorem*{prf*}{\underline{Proof}}

\newtheoremstyle{mystyle}
    {}
    {}
    {\normalfont}
    {}
    {\bfseries}
    {}
    { }
    {}
\theoremstyle{mystyle}

\begin{document}

\title{
  Explicit formula of deformation quantization\\
  with separation of variables for complex\\
  two-dimensional locally symmetric K\"{a}hler manifold
}

\author{
  Taika Okuda
  \thanks{
    Graduate School of Science, Department of Mathematics and Science Education, Tokyo University of Science, 1-3 Kagurazaka, Shinjuku-ku, Tokyo 162-8601, Japan
  }
  \and
  Akifumi Sako
  \thanks{
    Faculty of Science Division II, Department of Mathematics, Tokyo University of Science, 1-3 Kagurazaka, Shinjuku-ku, Tokyo 162-8601, Japan
  }
}

\date{\empty}

\maketitle
\allowdisplaybreaks

\begin{abstract}
We give a complex two-dimensional noncommutative locally symmetric K\"{a}hler manifold via a deformation quantization with separation of variables. We present an explicit formula of its star product by solving the system of recurrence relations given by Hara-Sako. In the two-dimensional case, this system of recurrence relations gives two types of equations corresponding to the two coordinates. From the two types of recurrence relations, symmetrized and antisymmetrized recurrence relations are obtained. The symmetrized one gives the solution of the recurrence relation. From the antisymmetrized one, the identities satisfied by the solution are obtained. The star products for \(\mathbb{C}^{2}\) and \(\mathbb{C}P^{2}\) are constructed by the method obtained in this study, and we verify that these star products satisfy the identities.
\end{abstract}

\section{Introduction}
\label{INTRO}
\ Deformation quantization is one of the quantization method based on a deformation for a Poisson algebra and is known as a method of constructing noncommutative differentiable manifolds. There are two types of this, ``formal deformation quantization" proposed by Bayen et al. \cite{BFFLS}, and ``strict deformation quantization", based on \(C^{\star}\)-algebra proposed by Rieffel \cite{Rief1}. In this paper, we study formal deformation quantization, and in the following, ``deformation quantization" is used in the sense of ``formal deformation quantization".

\begin{definition}
  Let 
  \(M\)
  be a Poisson manifold, 
  \(
    C^{\infty}\left(M\right)
  \)
  be a set of \(C^{\infty}\) functions on \(M\),
  \(
  \left\{
    \cdot,
    \cdot
  \right\}
  :
  C^{\infty}\left(M\right)
  \times
  C^{\infty}\left(M\right)
  \to
  C^{\infty}\left(M\right)
  \)
  be a Poisson bracket, and 
  \(
    C^{\infty}
    \left(M\right)
    \left[\!
        \left[
          \hbar
        \right]\!
      \right]
    :=
    \left\{
      f
      \middle|
      f
      =
      \sum_{k}
      f_{k}
      \hbar^{k},
      \
      f_{k}
      \in
      C^{\infty}\left(M\right)
    \right\}
  \)
  be the ring of formal power series over
  \(
    C^{\infty}\left(M\right),
  \)
   where \(\hbar\) is a formal parameter. Let a product \(\ast\) on
  \(
    C^{\infty}
    \left(M\right)
    \left[\!
      \left[
        \hbar
      \right]\!
    \right],
  \)
   called the star product, be a product denoted by
  \begin{align*}
    f
    \ast
    g
    =
    \sum_{k=0}^{\infty}
    C_{k}\left(f,g\right)
    \hbar^{k}
  \end{align*}
  satisfying the following conditions:
  \begin{enumerate}
    \item
    \(\ast\)
    is associative, i.e. for any 
    \(
      f,g,h
      \in
      C^{\infty}
      \left(M\right)
      \left[\!
        \left[
          \hbar
        \right]\!
      \right]
    \)
    , 
    \(
      f
      \ast
      \left(
        g
        \ast
        h
      \right)
      =
      \left(
        f
        \ast
        g
      \right)
      \ast
      h.
    \)
    \item
    Each 
    \(
      C_{k}
      \left(
        \cdot,
        \cdot
      \right)
      :
      C^{\infty}
      \left(M\right)
      \left[\!
        \left[
          \hbar
        \right]\!
      \right]
      \times
      C^{\infty}
      \left(M\right)
      \left[\!
        \left[
          \hbar
        \right]\!
      \right]
      \to
      C^{\infty}
      \left(M\right)
      \left[\!
        \left[
          \hbar
        \right]\!
      \right]
    \)
    is a bi-differential operator, i.e. for any 
    \(
      f,g
      \in
      C^{\infty}
      \left(M\right)
      \left[\!
        \left[
          \hbar
        \right]\!
      \right],
    \)
    \(
      C_{k}
      \left(
        f,
        g
      \right)
    \) 
    can be written as 
    \begin{align*}
      C_{k}
      \left(
        f,
        g
      \right)
      =
      \sum_{I,J}
      a_{I,J}
      \partial^{I}
      f
      \partial^{J}
      g,
    \end{align*}
    where \(I,J\) are multi-indices.
    \item
    For any
    \(
      f,g
      \in
      C^{\infty}
      \left(M\right)
    \)
    ,
    \begin{align*}
      C_{0}
      \left(
        f,
        g
      \right)
      &=
      fg,
      \\
      C_{1}
      \left(
        f,
        g
      \right)
      -
      C_{1}
      \left(
        g,
        f
      \right)
      &=
      \left\{f,g\right\}.
    \end{align*}
    \item
    For any
    \(
      f
      \in
      C^{\infty}
      \left(M\right)
      \left[\!
        \left[
          \hbar
        \right]\!
      \right]
    \)
    ,
    \(
      f
      \ast
      1
      =
      1
      \ast
      f
      =
      f.
    \)
  \end{enumerate}
  The pair 
  \(
    \left(
      C^{\infty}
      \left(M\right)
      \left[\!
        \left[
          \hbar
        \right]\!
      \right],
      \ast
    \right)
  \) 
  is called a deformation quantization for \(M.\)
  \label{dq_def}
\end{definition}
\noindent
For a more detailed review of deformation quantization, see e.g. \cite{GUT}. The construction method of deformation quantization for symplectic manifolds has been known by de Wilde-Lecomte \cite{DL}, Omori-Maeda-Yoshioka \cite{OMY1} and Fedosov \cite{Fedo}. After these works, a method for Poisson manifolds was proposed by Kontsevich \cite{Kont1}. For any K\"{a}hler manifold, Karabegov studied a construction method of deformation quantization with separation of variables \cite{Kara1,Kara2}.

\begin{definition}
  Let \(M\) be a K\"{a}hler manifold. A star product \(\ast\) on \(M\) is the separation of variables if the following two conditions are satisfied for any open set \(U\) of \(M\) and 
  \(
    f
    \in
    C^{\infty}
    \left(
      U
    \right)
  \) 
  :
  \begin{enumerate}
    \item
    For a holomorphic function
    \(
      a
    \) 
    on \(U,\) 
    \(
      a
      \ast
      f
      =
      af.
    \)
    \item
    For an anti-holomorphic function 
    \(
      b
    \) 
    on \(U,\) 
    \(
      f
      \ast
      b
      =
      fb.
    \)
  \end{enumerate}
  \label{sep_of_bar_starprd}
\end{definition}

\noindent
\ Furthermore, inspired by Karabegov's idea, a construction method for a locally symmetric K\"{a}hler manifold, i.e. a K\"{a}hler manifold such that 
\(
  \nabla_{\partial_{E}}
  R_{ABC}^{\quad \ \ \, D}
  =
  0
\) 
for 
\(
  A,B,C,D,E
  \in
  \left\{
    1,
    \cdots,
    N,
    \overline{1},
    \cdots,
    \overline{N}
  \right\},
\) 
was later proposed by Sako-Suzuki-Umetsu \cite{SSU1,SSU2} and Hara-Sako \cite{HS1,HS2}. Some notations in this paper are explained in more detail in Appendix \ref{app_kaehler}.
\bigskip
\\
\ In this paper, we propose an explicit formula that gives a deformation quantization with separation of variables for a complex two-dimensional locally symmetric K\"{a}hler manifold. This main result that is given in Theorem \ref{main_thm_N2} in Section \ref{NEW_RESULT} gives the explicit star product which is expanded in differential operators whose coefficients consist of covariantly constant. Each coefficient is explicitly determined by some matrix multiplications, and it contains the Riemann curvature tensor. This theorem is shown by solving the recurrence relations given by Hara-Sako \cite{HS1,HS2}. To explain our main result, we must introduce several definitions. So we shall not state our main theorem concretely, here.
\bigskip
\\
\
This paper is organized mainly into four Sections and three Appendices. In Section \ref{DQ_Review}, we review the previous works by Karabegov \cite{Kara1,Kara2} and Hara-Sako \cite{HS1,HS2}, as the background concerning a deformation quantization with separation of variables for K\"{a}hler manifolds.
In Section \ref{NEW_RESULT}, we show our main results which are the explicit formula to give the star product and the identities. In Section \ref{conc_ex}, we construct concrete examples for \(\mathbb{C}^{2}\) and \(\mathbb{C}P^{2}\) and they reproduce the previous results. In Section \ref{DQ_Discussion}, we state future works related to our main results from both mathematical and physical perspectives. In each Appendices \ref{app_kaehler}--\ref{prf_Tn_dagger}, we describe the properties and detailed calculations used in this paper. In Appendix \ref{app_kaehler}, we summarize some properties of K\"{a}hler manifolds used in this paper. In Appendix \ref{app_minus_n2}, we calculate in detail the identity \eqref{minus_thm_eq_ens} in Subsection \ref{NR_minus} for the 2nd order. In Appendix \ref{prf_Tn_dagger}, the Hermiteness of the coefficients of a star product is shown.

\section{Review of Noncommutative K\"{a}hler manifolds}
\label{DQ_Review}

\ The Quantization of K\"{a}hler manifolds was studied by Berezin \cite{Bere1,Bere2}, Moreno \cite{Mor1,Mor2}, Cahen-Gutt-Rawnsley \cite{CGR1,CGR2,CGR3,CGR4}, Karabegov \cite{Kara1,Kara2}, Omori-Maeda-Miyazaki-Yoshioka \cite{OMMY}, Schlichenmaier \cite{Schli1,Schli2,Schli3,Schli4,Schli5,Schli6,Schli7}, Karabegov-Schlichenmaier \cite{KS1}, Sako-Suzuki-Umetsu \cite{SSU1,SSU2} and Hara-Sako \cite{HS1,HS2}. In particular, Karabegov's method was proposed as a way to give noncommutative K\"{a}hler manifolds via a deformation quantization with separation of variables. After that, Sako-Suzuki-Umetsu and Hara-Sako methods were proposed for a locally symmetric case, inspired by this method. In addition, Sako-Suzuki-Umetsu was mentioned the Fock representations of noncommutative \(\mathbb{C}P^{N}\) and \(\mathbb{C}H^{N}\). Moreover, these previous results were generalized for any noncommutative K\"{a}hler manifolds by Sako-Umetsu \cite{SU1,SU2,SU3}. In Section \ref{NCK_History}, we review the methods by Karabegov as the background of this paper, and in Section \ref{DQ_LSK}, we review the method by Hara-Sako since our result is obtained from the recurrence relations in this method.

\subsection{
Noncommutative K\"{a}hler manifolds}
\label{NCK_History}

\ Berezin proposed a general definition of quantization and constructed the quantization of K\"{a}hler manifolds in the case of phase space via symbol algebras \cite{Bere1,Bere2}.  
The coherent states of K\"{a}hler manifolds arising from the geometric quantization of Kostant \cite{Kos} and Souriau \cite{Sou} have also been studied by Rawnsley \cite{Rawn}. It is known that this coherent state is related to Berezin quantization. See \cite{Pere} for more detail. 
After that, the deformation quantization of K\"{a}hler manifolds have been provided by Moreno \cite{Mor1,Mor2} and Omori-Maeda-Miyazaki-Yoshioka \cite{OMMY}. The relations between deformation quantization and Berezin quantization have been studied by Cahen-Gutt-Rawnsley \cite{CGR1,CGR2,CGR3,CGR4}. 
It has also studied the quantization of K\"{a}hler manifolds via Toeplitz quantization by Bordemann et al. \cite{BMS}. 
Furthermore, Karabegov and Schlichenmaier have provided Berezin-Toeplitz quantization in the case of compact K\"{a}hler ones \cite{Schli1,Schli2,Schli3,Schli4,Schli5,Schli7}. These previous works related to Berezin-Toeplitz quantizations for K\"{a}hler manifolds were reviewed by Schlichenmaier \cite{Schli6}. 
From the other angle of the quantization, the construction method of noncommutative K\"{a}hler manifolds was studied via the deformation quantization with separation of variables by Karabegov \cite{Kara1,Kara2}. 
Moreover, for any noncommutative K\"{a}hler manifolds obtained by Karabegov's construction, Sako-Umetsu constructed Fock representations of them \cite{SU1,SU2}. In this subsection, we review 
Karabegov's method and Fock representations of noncommutative K\"{a}hler manifolds by Sako-Umetsu.
\bigskip
\\
\ Let \(M\) be an \(N\)-dimensional K\"ahler manifold and 
\(
  U
  \subset
  M
\) be a holomorphic coordinate neighborhood of \(M\). We choose a local holomorphic coordinates by \(\left(z^{1},\cdots,z^{N}\right)\). For a K\"{a}hler manifolds, a K\"{a}hler 2-form \(\omega\) and a K\"{a}hler metric \(g\) can be locally expressed by using a K\"{a}hler potential \(\Phi\) as follows:
\begin{align*}
  \omega
  =
  i
  g_{k\bar{l}}
  dz^{k}
  \wedge
  d\bar{z}^{l},
  \qquad
  g_{k\bar{l}}
  =
  \frac{
    \partial^2
    \Phi
  }{
    \partial
    z^{k}
    \partial
    \bar{z}^{l}
  }.
\end{align*}
\noindent
Note that we use the Einstein summation convention on the above. We also denote the inverse matrix of \(\left(g_{\overline{k}l}\right)\) by \(\left(g^{\overline{k}l}\right)\). Here we introduce the differential operators
\(
  D^{k},
  D^{\overline{k}}
\) 
defined by
\begin{align*}
  &
  D^{k}
  :=
  g^{k\overline{l}}
  \partial_{\overline{l}},
  \quad
  D^{\overline{k}}
  :=
  g^{\overline{k}l}
  \partial_{l}.
\end{align*}
\noindent
We define the set of differential operators 
\begin{align*}
  {\mathcal S}
  :=
  \left\{ 
    A
    \middle|
    A
    =
    \sum_{\alpha}
    a_{\alpha}
    D^{\alpha},
    \,
    a_{\alpha}
    \in
    C^{\infty}
    \left(U\right)
  \right\},
\end{align*}
\noindent
where 
\(
  \alpha
  =
  \left(
    \alpha_{1},
    \dots,
    \alpha_{N}
  \right)
\) is a multi-index
, i.e. 
\(
  D^{\alpha}
  :=
  \left(
    D^{\bar{1}}
  \right)^{\alpha_{1}}
  \cdots
  \left(
    D^{\bar{N}}
  \right)^{\alpha_{N}}
\). 
We can construct the left star-multiplication operator \(L_{f}\) for 
\(
  f
  \in
  C^{\infty}
  \left(U\right)
\) 
such that 
\(
  L_{f}
  g
  :=
  f
  \ast
  g
\)
:
\begin{theorem}[Karabegov\cite{Kara1}]
Let \(M\) be an \(N\)-dimensional K\"{a}hler manifold, \(U\) be a 
holomorphic coordinate neighborhood on \(M\), 
and \(\omega\) be a K\"{a}hler form on \(M\). 
Then, there is 
the left star-multiplication operator
\begin{align*}
  L_{f}
  =
  \sum_{n=0}^{\infty}
  \hbar^{n}
  A_{n},
  \qquad
  f
  \in
  C^{\infty}
  \left(U\right),
\end{align*}
where 
\(
  A_{n}
  :=
  a_{n,\alpha}\left(f\right)
  D^{\alpha}
  \in
  \mathcal{S}
\) 
are differential operators whose coefficients
\(
  a_{n,\alpha}\left(f\right)
  \in
  C^{\infty}
  \left(U\right)
\) 
 depend on \(f\).
\(L_{f}\) is determined by the following conditions:
\begin{enumerate}
\item
\(
  \left[
    L_{f},
    R_{\partial_{\overline{l}}\Phi}
  \right]
  =
  0
\), 
where 
\(
  R_{\partial_{\overline{l}}\Phi}
  =
  \partial_{\overline{l}}\Phi
  +
  \hbar
  \partial_{\overline{l}}
\),
\item
\(
  L_{f}1
  =
  f
  \ast
  1
  =
  f
\),
\item
For any 
\(
  g,h
  \in
  C^{\infty}
  \left(U\right)
\), 
the left star-multiplication operator is associative, i.e. 
\begin{align*}
  L_{f}
  \left(
    L_{g}
    h
  \right)
  =
  f
  \ast
  \left(
    g
    \ast
    h
  \right)
  =
  \left(
    f
    \ast
    g
  \right)
  \ast
  h
  =
  L_{
    L_{f}g
  }
  h
  .
\end{align*}
\end{enumerate}
\label{thm_Kara_star_opr}
\end{theorem}
\noindent
By using the definition of the separation of variables and the commutation relations of the star-multiplication operators, we obtain the following commutation relations 
concerning \(z^{i}\), \(\overline{z}^{i}\), \(\partial_{i}\Phi\) and \(\partial_{\overline{i}}\Phi\).
\begin{align}
  \left[
    \frac{1}{\hbar} 
    \partial_{i}
    \Phi,
    z^{j}
  \right]_{\ast}
  &=
  \delta_{ij}, 
  \qquad
  \left[
    z^{i},
    z^{j}
  \right]_{\ast}
  =
  0,
  \qquad
  \left[
    \partial_{i}
    \Phi,
    \partial_{j}
  \Phi
  \right]_{\ast}  
  =
  0
  \label{star_comm_rel_1}
  \\
  \left[
    \overline{z}^{i},
    \frac{1}{\hbar}
    \partial_{\overline{j}}
    \Phi
  \right]_{\ast}
  &=
  \delta_{ij},
  \qquad
  \left[
    \overline{z}^{i},
    \overline{z}^{j}
  \right]_{\ast}
  =
  0,
  \qquad
  \left[
    \partial_{\overline{i}}
    \Phi,
    \partial_{\overline{j}}
    \Phi
  \right]_{\ast}
  =
  0.
  \label{star_comm_rel_2}
\end{align}
\noindent
Note that the commutator \(\left[\cdot,\cdot\right]_{\ast}\) is defined by 
\(
  \left[
    A,B
  \right]_{\ast}
  :=
  A
  \ast
  B
  -
  B
  \ast
  A
\). 
Here, we introduce the creation and annihilation operators 
as follows:
\begin{align}
  a^{\dagger}_{i}
  =
  z^{i},
  \qquad
  \underline{a}_{i}
  =
  \frac{1}{\hbar}
  \partial_{i}
  \Phi,
  \qquad
  a_{i}
  =
  \bar{z}^{i},
  \qquad
  \underline{a}_{i}^{\dagger}
  =
  \frac{1}{\hbar}
  \partial_{\bar{i}}
  \Phi
  \qquad
  \left(
    i=1,\cdots,N
  \right).
  \label{cre_ann_Kahler}
\end{align}
\noindent
Then, the commutation relations \eqref{star_comm_rel_1} and \eqref{star_comm_rel_2} can be rewritten as
\begin{align}
  \left[
    \underline{a}_{i},
    a^{\dagger}_{j}
  \right]_{\ast}
  &=
  \delta_{ij}, 
  \qquad
  \left[
    a^{\dagger}_{i},
    a^{\dagger}_{j}
  \right]_{\ast}
  =
  0,
  \qquad
  \left[
    \underline{a}_{i},
    \underline{a}_{j}
  \right]_{\ast}  
  =
  0
  \label{star_comm_rel_CA_1}
  \\
  \left[
    a_{i},
    \underline{a}_{j}^{\dagger}
  \right]_{\ast}
  &=
  \delta_{ij},
  \qquad
  \left[
    a_{i},
    a_{j}
  \right]_{\ast}
  =
  0,
  \qquad
  \left[
    \underline{a}_{i}^{\dagger},
    \underline{a}_{j}^{\dagger}
  \right]_{\ast}
  =
  0.
  \label{star_comm_rel_CA_2}
\end{align}
\noindent
Since 
\(
  \left[
    a_{i},
    a_{i}^{\dagger}
  \right]_{\ast}
  \neq
  0
\), 
\(
  \left[
    \underline{a}_{i},
    \underline{a}_{j}^{\dagger}
  \right]_{\ast}
  \neq
  0
\) 
in general, these relations are slightly different from the ordinary canonical commutation relation. From the 
above operators, the (twisted) Fock space is defined by a vector space spanned by the basis
\begin{align}
  |{\vec{n}}\rangle
  =
  |n_{1},\cdots,n_{N}\rangle
  =
  \frac{1}{
    \sqrt{
      \vec{n}!
    }
  }
 {\left(a^{\dagger}_{1}\right)}_{\ast }^{n_{1}}
 \ast
 \cdots
 \ast
 {\left(a^{\dagger}_{N}\right)}_{\ast }^{n_{N}}
 \ast  
 |\vec{0}\rangle,
 \label{twFock_base_n}
\end{align}
\noindent
where a vacuum 
\(
  |{\vec{0}}\rangle
  =
  |0,\cdots,0\rangle 
\) 
is the vector such that
\begin{align}
  \underline{a}_{i}
  \ast
  |{\vec{0}}\rangle
  =
  0,
  \qquad
  (i=1,\cdots,N),
  \label{twFock_vac_n}
\end{align}
\noindent
and  
\(
  \vec{n}!
  =
  n_{1}!
  \cdots
  n_{N}!
\). 
Note that 
\(
  \left(A\right)_{\ast}^{n}
\) 
is the product of multiplying \(n\) by the star product \(\ast\), i.e. 
\(
  \left(A\right)_{\ast}^{n}
  :=
  \underbrace{
    A
    \ast
    \cdots
    \ast
    A
  }_{
    n
    \
    \rm{
      times
      \
      by
    }
    \
    \ast
  }
\). Similarly, the dual basis for \(|{\vec{n}}\rangle\) is defined by
\begin{align}
  \underline{\langle{\vec{m}}|}
  =
  \underline{
    \langle
    m_{1},
    \cdots,
    m_{N}|
  }
  =
  \langle
  \vec{0}|
  \ast
  \left(\underline{a}_{1}\right)_{\ast }^{m_{1}}
  \ast
  \cdots
  \ast  
  \left(\underline{a}_{N}\right)_{\ast }^{m_{N}} 
  \frac{1}{
    \sqrt{
      \vec{m}!
    }
  },
 \label{twFock_dbase_m}
\end{align}
\noindent
where \(\langle\vec{0}|\) is the (dual) vector for a vacuum \(|{\vec{0}}\rangle\) such that
\begin{align}
  \langle
  \vec{0}|
  \ast 
  a_{i}^{\dagger}
  =
  0,
  \qquad
  (i=1,\cdots,N),
\end{align}
\noindent
and 
\(
  \vec{m}!
  =
  m_{1}!
  \cdots
  m_{N}!
\).
Note that \(\underline{\langle{\vec{m}}|}\) does not imply Hermitian conjugate of \(|{\vec{m}}\rangle\), i.e. 
\(
  \underline{\langle{\vec{m}}|}
  \neq
  |{\vec{m}}\rangle^{\dagger}
\).
\begin{definition}
Let \(M\) be a K\"{a}hler manifold and \(U\) be a holomorphic coordinate neighborhood on \(M\). Then, the (local) twisted Fock algebra (representation) \(F_{U}\) is defined by
\begin{align}
  F_{U}
  :=
  \left\{
    \sum_{\vec{n},\vec{m}}
    A_{\vec{n}\vec{m}}
    |\vec{n}\rangle
    \underline{\langle\vec{m}|}
    \middle|
    A_{
      \vec{n}
      \vec{m}
    }
    \in
    \mathbb{C}
  \right\}.
\end{align}
\label{twFock_def}
\end{definition}
\noindent
\(F_{U}\) is defined as the algebra which is given by the creation and annihilation operators in \eqref{cre_ann_Kahler} and star-multiplication between each element of \(F_{U}\). Moreover, we can concretely express the coefficient functions \(A_{\vec{n}\vec{m}}\) which are the elements of \(F_{U}\). 
We expand a function 
$\exp \Phi(z, \bar{z})/\hbar$ as a power series,
\begin{align}
 e^{\Phi(z, \bar{z})/\hbar}
 = \sum_{\vec{m},\vec{n}} H_{\vec{m}, \vec{n}}
 (z)^{\vec{m}} (\bar{z})^{\vec{n}}
 \label{H}
\end{align}
where $(z)^{\vec{n}} = (z^1)^{n_1} \cdots (z^N)^{n_N}$ and $(\bar{z})^{\vec{n}} = (\bar{z}^1)^{n_1} \cdots (\bar{z}^N)^{n_N}$.
The creation and annihilation operators $a_i^\dagger, \underline{a}_i$
act on the bases as follows,
\begin{align}
 a_i^\dagger \ast  |\vec{m} \rangle \underline{\langle \vec{n}|}
 = \sqrt{m_i+1} 
 |\vec{m} + \vec{e}_i \rangle \underline{\langle \vec{n}|},
 \qquad
 &
 \underline{a}_i \ast  |\vec{m} \rangle \underline{\langle \vec{n}|}
 = \sqrt{m_i} 
 |\vec{m} - \vec{e}_i \rangle \underline{\langle \vec{n}|}, \\
 |\vec{m} \rangle \underline{\langle \vec{n}|} \ast  a_i^\dagger 
 = \sqrt{n_i} 
 |\vec{m} \rangle \underline{\langle \vec{n} - \vec{e}_i|},
 \qquad
 &
 |\vec{m} \rangle \underline{\langle \vec{n}|} \ast  \underline{a}_i
 = \sqrt{n_i+1} 
 |\vec{m} \rangle \underline{\langle \vec{n} + \vec{e}_i|}, 
\end{align}
where $\vec{e}_i$ is a unit vector, $(\vec{e}_i)_j = \delta_{ij}$.
The action of $a_i$ and $\underline{a}_i^\dagger$ is derived by the
Hermitian conjugation of the above equations.
%
%
The results of the twisted Fock representation of the
noncommutative K\"{a}hler manifolds are summarized as the 
dictionary in the Table \ref{table:result}:
\begin{table}[H]
 \caption{Functions - Fock operators Dictionary }
 \label{table:result}
 \centering
 \begin{tabular}{c|c}
  \hline
  Functions & Fock operators  \\
  \hline \hline
  $\displaystyle e^{-\Phi/\hbar}$  &  $|\vec{0}\rangle \langle \vec{0}|$
      \\
  \hline
  $z_i $    &   $a_i^{\dagger}$ 
      \\ \hline
  $\displaystyle \frac{1}{\hbar}\partial_i \Phi$   &  $\underline{a}_i$
      \\ \hline
  $\bar{z}^i$   
  & $\displaystyle a_i = 
      \sum \sqrt{\frac{\vec{m}!}{\vec{n}!}} H_{\vec{m},\vec{k}}
      H^{-1}_{\vec{k}+\vec{e}_i , \vec{n}} 
      |\vec{m}\rangle \underline{\langle \vec{n} |} $ 
      \\ \hline
  $\displaystyle \frac{1}{\hbar} \partial_{\bar i} \Phi$ 
  & 
      $\displaystyle \underline{a}_i^{\dagger} = \sum 
      \sqrt{\frac{\vec{m}!}{\vec{n}!}}
      (k_i+1) H_{\vec{m},\vec{k}+\vec{e}_i}
      H^{-1}_{\vec{k} , \vec{n}} |\vec{m}\rangle 
      \underline{ \langle \vec{n} |} $
      \\
  \hline
 \end{tabular}
\end{table}
For physics, it is difficult to interpret formal power series. So this Fock representation is useful to construct physics theories like field theories on noncommutative K\"{a}hler manifolds. We will discuss this point in Section \ref{DQ_Discussion}.

\subsection{Deformation quantization with separation of variables for locally symmetric K\"{a}hler manifold}
\label{DQ_LSK}

\ Let \(M\) be a complex \(N\)-dimensional locally symmetric K\"{a}hler manifold, \(U\) be a holomorphic coordinate neighborhood of \(M\). For any 
\(
  f,g
  \in
  C^{\infty}
  \left(
    U
  \right),
\) 
we assume the following form for a star product with separation of variables on \(M\) by Sako-Suzuki-Umetsu \cite{SSU1,SSU2} and Hara-Sako \cite{HS1,HS2} : 
\begin{align}
  f\ast g
  =
  L_{f}g
  :=
  \sum_{n=0}^{\infty}
  \sum_{
    \overrightarrow{\alpha_{n}},
    \overrightarrow{\beta_{n}^{\ast}}
  }
  T_{\overrightarrow{\alpha_{n}},\overrightarrow{\beta_{n}^{\ast}}}^{n}
  \left(
    D^{\overrightarrow{\alpha_{n}}}f
  \right)
  \left(
    D^{\overrightarrow{\beta_{n}^{\ast}}}g
  \right).
  \label{star_expand}
\end{align}
\noindent
Here \(L_{f}\) is a left star-multiplication operator with respect to \(f\), and
\(
  D^{\overrightarrow{\alpha_{n}}},
  D^{\overrightarrow{\beta_{n}^{\ast}}}
\)
are differential operators defined by
\begin{align*}
  &
  D^{k}
  :=
  g^{k\overline{l}}
  \partial_{\overline{l}},
  \
  D^{\overline{k}}
  :=
  g^{\overline{k}l}
  \partial_{l},
  \\
  &
  D^{\overrightarrow{\alpha_{n}}}
  :=
  D^{\alpha_{1}^{n}}\cdots D^{\alpha_{N}^{n}},
  \quad
  D^{\alpha_{k}^{n}}
  :=
  \left(
    D^{k}
  \right)^{\alpha_{k}^{n}},
  \\
  &
  D^{\overrightarrow{\beta_{n}^{\ast}}}
  :=
  \overline{
    D^{\overrightarrow{\beta_{n}}}
  }
  =
  \overline{D^{\beta_{1}^{n}}}
  \cdots
  \overline{D^{\beta_{N}^{n}}},
  \quad
  \overline{D^{\beta_{k}^{n}}}
  :=
  \left(
    D^{\overline{k}}
  \right)^{\beta_{k}^{n}},
  \\
  &
  \overrightarrow{\alpha_{n}},
  \overrightarrow{\beta_{n}}
  \in
  \left\{
    \left(\gamma_{1}^{n},\cdots,\gamma_{N}^{n}\right)
    \in
    \mathbb{Z}^{N}
    \middle|
    \sum_{k=1}^{N}
    \gamma_{k}^{n}
    =
    n
  \right\},
\end{align*}
\noindent
respectively, where
\(
  \left\{
    \left(\gamma_{1}^{n},\cdots,\gamma_{N}^{n}\right)
    \in
    \mathbb{Z}^{N}
    \middle|
    \sum_{k=1}^{N}
    \gamma_{k}^{n}
    =
    n
  \right\}
\)
is a \(N\)-dimensional module such that a sum of all components is a non-negative integer \(n\). If 
there exists at least one negative 
\(
  \alpha_{k}^{n}
  \notin\mathbb{Z}_{\geq0}
\) 
for 
\(
  k
  \in
  \left\{
    1,\cdots,N
  \right\}
\), then we define 
\(
  D^{\overrightarrow{\alpha_{n}}}:=0
\). 
We define 
\(
  D^{\overrightarrow{\beta_{n}^{\ast}}}:=0
\) 
in the same way when \(\overrightarrow{\beta_{n}^{\ast}}\) has negative components. Note that \eqref{star_expand} is not a power series in the formal parameter \(\hbar\), 
but a power series 
in the differential operators 
\(
  D^{k}
\) 
and 
\(
  D^{\bar{k}}.
\)
 Since \(M\) is locally symmetric, the coefficients of the star product 
\(
  T_{\overrightarrow{\alpha_{n}},\overrightarrow{\beta_{n}^{\ast}}}^{n}
\) can be assumed to be covariantly constants.
%
If 
\(
  \overrightarrow{\alpha_{n}}
  \notin
  \mathbb{Z}_{\geq0}^{n}
\) 
or
\(
  \overrightarrow{\beta_{n}}
  \notin
  \mathbb{Z}_{\geq0}^{n},
\)
then we define
\(
  T_{\overrightarrow{\alpha_{n}},\overrightarrow{\beta_{n}^{\ast}}}^{n}
  :=
  0
\) 
as well as 
\(
  D^{\overrightarrow{\alpha_{n}}}
\) 
and 
\(
  D^{\overrightarrow{\beta_{n}}}
\). 
In the following discussion, we shall omit \(\ast\) when 
\(
  \overrightarrow{\beta_{n}^{\ast}}
\) 
is expressed in explicit form. 
For example, if 
\(
  \overrightarrow{\beta_{n}^{\ast}}
  =
  \left(1,2,3\right)^{\ast},
\) 
then we denote 
\(
  \left(1,2,3\right)^{\ast}
  =
  \left(1,2,3\right)
\). 
It is known that 
the coefficient 
\(
  T_{\overrightarrow{\alpha_{n}},\overrightarrow{\beta_{n}^{\ast}}}^{n}
\) 
for the zeroth and first orders in 
\(
  f
  \ast
  g
\) 
are given below.

\begin{prop}[Hara-Sako\cite{HS1,HS2}]
  For a star product with separation of variables \(\ast\) on \(U,\)
  \begin{align*}
    T_{\overrightarrow{\alpha_{0}},\overrightarrow{\beta_{0}^{\ast}}}^{0}
    =
    1,
    \
    T_{\overrightarrow{e_{i}},\overrightarrow{e_{j}}}^{1}
    =
    \hbar g_{i\overline{j}},
  \end{align*}
  \noindent
  where 
  \(
    \overrightarrow{e_{i}}
    =
    \left(\delta_{1i},\cdots,\delta_{Ni}\right).
  \)
  \label{prop:covar_const_0and1}
\end{prop}

\noindent
In other words, Proposition \ref{prop:covar_const_0and1} states that for any complex \(N\)-dimensional locally symmetric K\"{a}hler manifold, the zeroth and first orders are completely determined.

\begin{theorem}[Hara-Sako\cite{HS1,HS2}]
  For 
  \(
    f,g
    \in
    C^{\infty}\left(M\right),
  \) 
  there exists 
  a star product with separation of variables \(\ast\) such that
  \begin{align*}
    L_{f}
    g
    :=
    f
    \ast
    g
    =
    \sum_{n=0}^{\infty}
    \sum_{\overrightarrow{\alpha_{n}},\overrightarrow{\beta_{n}^{\ast}}}
    T_{\overrightarrow{\alpha_{n}},\overrightarrow{\beta_{n}^{\ast}}}^{n}
    \left(
      D^{\overrightarrow{\alpha_{n}}}f
    \right)
    \left(
      D^{\overrightarrow{\beta_{n}^{\ast}}}g
    \right),
    \
    \left(
      f,g
      \in
      C^{\infty}
      \left(
        M
      \right)
    \right),
  \end{align*}
  \noindent
  where the coefficient 
  \(
    T_{
      \overrightarrow{\alpha_{n}},
      \overrightarrow{\beta_{n}^{\ast}}
    }^{n}
  \) 
  satisfies the following recurrence relation :
  \begin{align*}
    &
    \sum_{d=1}^{N}
    \hbar
    g_{\overline{i}d}
    T_{\overrightarrow{\alpha_{n}}-\overrightarrow{e_{d}},\overrightarrow{\beta_{n}^{\ast}}-\overrightarrow{e_{i}}}^{n-1}
    \\
    &=
    \beta_{i}^{n}
    T_{\overrightarrow{\alpha_{n}},\overrightarrow{\beta_{n}^{\ast}}}^{n}
    +
    \sum_{k=1}^{N}
    \sum_{\rho=1}^{N}
    \frac{
      \hbar
      \left(\beta_{k}^{n}-\delta_{k\rho}-\delta_{ik}+1\right)
      \left(\beta_{k}^{n}-\delta_{k\rho}-\delta_{ik}+2\right)
    }{2}
    R_{\overline{\rho} \ \ \ \,\overline{i}}^{\ \,\overline{k}\,\overline{k}}
    T_{\overrightarrow{\alpha_{n}},\overrightarrow{\beta_{n}^{\ast}}-\overrightarrow{e_{\rho}}+2\overrightarrow{e_{k}}-\overrightarrow{e_{i}}}^{n}
    \\
    &
    \mbox{}
    \quad
    +
    \sum_{k=1}^{N-1}
    \sum_{l=1}^{N-k}
    \sum_{\rho=1}^{N}
    \hbar
    \left(
      \beta_{k}^{n}
      -
      \delta_{k\rho}
      -
      \delta_{ik}
      +
      1
    \right)
    \left(
      \beta_{k+l}^{n}
      -
      \delta_{k+l,\rho}
      -
      \delta_{i,k+l}
      +
      1
    \right)
    R_{\overline{\rho} \quad \ \ \ \,\overline{i}}^{\ \,\overline{k+l}\,\overline{k}}
    T_{
      \overrightarrow{\alpha_{n}},
      \overrightarrow{\beta_{n}^{\ast}}
      -
      \overrightarrow{e_{\rho}}
      +
      \overrightarrow{e_{k}}
      +
      \overrightarrow{e_{k+l}}
      -
      \overrightarrow{e_{i}}
    }^{n}.
  \end{align*}
  \label{theorem:star_prd}
\end{theorem}

\noindent
These recurrence relations in Theorem \ref{theorem:star_prd} are equivalent to the equations (6.9) in \cite{SSU1}. Hence, from Theorem \ref{theorem:star_prd}, a star product with separation of variables on any complex \(N\)-dimensional locally symmetric K\"{a}hler manifold is obtained. However, finding a general term 
\(
  T_{\overrightarrow{\alpha_{n}},\overrightarrow{\beta_{n}^{\ast}}}^{n}
\) 
that satisfies this system of recurrence relations is not easy, except in the one-dimensional case.

\begin{prop}[Hara-Sako\cite{HS1,HS2}]
  Let \(M\) be a one-dimensional locally symmetric K\"{a}hler manifold, and \(U\) be an open set of \(M\). For 
  \(
    f,g
    \in
    C^{\infty}\left(U\right)
  \), 
  the star product 
  \(
    f
    \ast
    g
  \) 
  is given by
  \begin{align*}
    f
    \ast
    g
    =
    \sum_{n=0}^{\infty}
    \left[
      \left(
        g_{1\overline{1}}
      \right)^{n}
      \left\{
        \prod_{k=1}^{n}
        \frac{4\hbar}{4k+\hbar k\left(k-1\right)R}
      \right\}
      \left\{
        \left(
        g^{1\overline{1}}
        \frac{\partial}{\partial \overline{z}}
      \right)^{n}
      f
      \right\}
      \left\{
        \left(
        g^{1\overline{1}}
        \frac{\partial}{\partial z}
      \right)^{n}
	g
      \right\}
    \right],
  \end{align*}
  where 
  \(
    R
    =
    2
    R_{\overline{1} \ \ \, \overline{1}}^{\ \, \overline{1}\overline{1}}
  \)
  is the scalar curvature on \(M\).
  \label{prop:cplx_1}
\end{prop}

\noindent
Note that Proposition \ref{prop:cplx_1} is corrected 
some errata for the one-dimensional formula in \cite{HS1,HS2}. This proposition can be shown by direct calculations 
since the recurrence relation is a simple expression in the one-dimensional case.
\bigskip
\\
\ On the other hand,  
\(
  T_{\overrightarrow{\alpha_{n}},\overrightarrow{\beta_{n}^{\ast}}}^{n}
\) 
were 
only obtained for 
\(
  n
  =
  0,
  \
  1
\) 
and 
\(
  2
\) 
for the two-dimensional case. The following proposition is shown by directly solving the recurrence relation in Theorem \ref{theorem:star_prd} for \(n=2\).

\begin{prop}[Hara-Sako\cite{HS1,HS2}]
  Let \(M\) be a complex two-dimensional locally symmetric K\"{a}hler manifold. Then the coefficients 
  \(
    T_{\overrightarrow{\alpha_{2}},\overrightarrow{\beta_{2}^{\ast}}}^{2}
  \)  
  are given by
  \begin{align*}
    \left(
      T^{2}_{ij}
    \right)
    =
    \hbar^{2}
    A_{2}
    X_{2}^{-1},
  \end{align*}
  
  \noindent
  where
  \begin{align*}
    &
    T^{2}_{ij}
    :=
    T_{
      \left(
        3-i,
        i-1
      \right),
      \left(
        3-j,
        j-1
      \right)
    }^{2},
    \\
    &
    A_{2}
    :=
    \left(
      \begin{array}{ccc}
        \left(
          g_{\bar{1}1}
        \right)^{2}
        &
        g_{\bar{1}1}
        g_{\bar{2}1}
        &
        \left(
          g_{\bar{2}1}
        \right)^{2}
        \\
        2
        g_{\bar{1}1}
        g_{\bar{1}2}
        &
        g_{\bar{1}2}
        g_{\bar{2}1}
        +
        g_{\bar{1}1}
        g_{2\bar{2}}
        &
        2
        g_{\bar{2}1}
        g_{\bar{2}2}
        \\
        \left(
          g_{\bar{1}2}
        \right)^{2}
        &
        g_{\bar{2}1}
        g_{\bar{2}2}
        &
        \left(
          g_{\bar{2}2}
        \right)^{2}
      \end{array}
    \right),
    \\
    &
    X_{2}
    :=
    \left(
      \begin{array}{ccc}
        2
        +
        \hbar
        R_{\bar{1} \ \ \, \bar{1}}^{\ \, \bar{1}\bar{1}}
        &
        \hbar
        R_{\bar{2} \ \ \, \bar{1}}^{\ \, \bar{1}\bar{1}}
        &
        \hbar
        R_{\bar{2} \ \ \, \bar{2}}^{\ \, \bar{1}\bar{1}}
        \\
        \hbar
        R_{\bar{1} \ \ \, \bar{1}}^{\ \, \bar{2}\bar{1}}
        &
        1
        +
        \hbar
        R_{\bar{2} \ \ \, \bar{1}}^{\ \, \bar{2}\bar{1}}
        &
        \hbar
        R_{\bar{2} \ \ \, \bar{2}}^{\ \, \bar{2}\bar{1}}
        \\
        \hbar
        R_{\bar{1} \ \ \, \bar{1}}^{\ \, \bar{2}\bar{2}}
        &
        \hbar
        R_{\bar{2} \ \ \, \bar{1}}^{\ \, \bar{2}\bar{2}}
        &
        2
        +
        \hbar
        R_{\bar{2} \ \ \, \bar{2}}^{\ \, \bar{2}\bar{2}}
      \end{array}
    \right).
  \end{align*}
  \label{prop:n2}
\end{prop}

\noindent
We have reviewed the previous works. In general, solving recurrence relations is not easy. In particular, when we attempt to obtain the general term using a matrix representation, we need a square matrix of order \(n+1\). In addition, the matrix size increases with increasing order. Therefore, it had been considered difficult to obtain the general term. In this paper, we find that the expression of the general term does not become unlimitedly complex, and succeed in getting the general term by using this fact. We shall see that in the following sections. In Section \ref{NEW_RESULT}, we are going to describe the coefficients 
\(
  T_{
    \overrightarrow{\alpha_{n}},
    \,
    \overrightarrow{\beta_{n}^{\ast}}
  }^{n}
\) 
for any \(n\in\mathbb{Z}_{\geq0}\).

\section{Star product with separation of variables for a complex\\
two-dimensional locally symmetric K\"{a}hler manifold}
\label{NEW_RESULT}
\ In this section, we construct 
the formula that explicitly determines a star product with separation of variables 
for a complex two-dimensional locally symmetric K\"{a}hler manifold. 
In other words, we construct the solution of the recurrence relations in Theorem \ref{theorem:star_prd} for the two-dimensional case in this section.

\subsection{Complex two-dimensional formula}
\label{NR_plus}
\ The explicit formula for the coefficient 
\(
  T_{
    \overrightarrow{\alpha_{n}},
    \,
    \overrightarrow{\beta_{n}^{\ast}}
  }^{n}
\) 
for any order 
\(
  n
  \in
  \mathbb{Z}_{\geq0}
\) 
is obtained by not dealing with the recurrence relations independently but attributing them to one recurrence relation.
\begin{theorem}
  Let \(M\) be a complex two-dimensional locally symmetric K\"{a}hler manifold, \(U\) be an open set of \(M\), and 
  \(
    \ast
  \) 
  be a star product with separation of variables on \(U\) such that
  \begin{align*}
    f
    \ast
    g
    =&
    \sum_{n=0}^{\infty}
    \sum_{\overrightarrow{\alpha_{n}},\,\overrightarrow{\beta_{n}^{\ast}}}
    T_{\overrightarrow{\alpha_{n}},\,\overrightarrow{\beta_{n}^{\ast}}}^{n}
    \left(
      D^{\overrightarrow{\alpha_{n}}}f
    \right)
    \left(
      D^{\overrightarrow{\beta_{n}^{\ast}}}g
    \right)
    \\
    =&
    \sum_{n=0}^{\infty}
    \sum_{i,j=1}^{n+1}
    T_{ij}^{n}
    \left\{
      \left(
        D^{1}
      \right)^{n-i+1}
      \left(
        D^{2}
      \right)^{i-1}
      f
    \right\}
    \left\{
      \left(
        D^{\overline{1}}
      \right)^{n-j+1}
      \left(
        D^{\overline{2}}
      \right)^{j-1}
      g
    \right\}
  \end{align*}
  for 
  \(
    f,g
    \in
    C^{\infty}
    \left(U\right).
  \) 
  Then, 
  \begin{align}
    T_{n}
    X_{n}
    =
    \hbar
    A'_{n},
    \label{main_them_N2_T}
  \end{align}
  where
  \begin{align*}
    T_{n}
    &
    :=
    \left(
      T^{n}_{ij}
    \right)
    \in
    M_{n+1}
    \left(
    \mathbb{C}
    \left[\!
      \left[
        \hbar
      \right]\!
    \right]
  \right),
    \\
    T^{n}_{ij}
    &
    :=
    T_{
      \left(
        n-i+1,
        i-1
      \right),
      \left(
        n-j+1,
        j-1
      \right)
    }^{n},
    \\
    A'_{n}
    &
    =
    \left(A_{ij}^{'n}\right)
    =
    \left(
      \sum_{k,l=1}^{2}
      g_{\overline{k}l}
      T^{n-1}_{i-\delta_{2l},j-\delta_{2k}}
    \right)
    \in
    M_{n+1}
    \left(
    \mathbb{C}
    \left[\!
      \left[
        \hbar
      \right]\!
    \right]
  \right),
\end{align*}
and 
\(
  X_{n}
  \in
  M_{n+1}
  \left(
    \mathbb{C}
    \left[\!
      \left[
        \hbar
      \right]\!
    \right]
  \right)
\)
is a pentadiagonal matrix such that these components are given as follows :
\begin{align*}
    X_{j-2,j}^{n}
    &=
    \tbinom{n-j+3}{2}
    \hbar
    R_{\bar{2} \ \ \, \bar{2}}^{\ \, \bar{1}\bar{1}},
    \\
    X_{j-1,j}^{n}
    &=
    2
    \tbinom{n-j+2}{2}
    \hbar
    R_{\bar{2} \ \ \, \bar{1}}^{\ \, \bar{1}\bar{1}}
    +
    \left(
      n-j+2
    \right)
    \left(
      j-2
    \right)
    \hbar
    R_{\bar{2} \ \ \, \bar{2}}^{\ \, \bar{2}\bar{1}},
    \\
    X_{j,j}^{n}
    &=
    n
    +
    \tbinom{n-j+1}{2}
    \hbar
    R_{\bar{1} \ \ \, \bar{1}}^{\ \, \bar{1}\bar{1}}
    +
    \tbinom{j-1}{2}
    \hbar
    R_{\bar{2} \ \ \, \bar{2}}^{\ \, \bar{2}\bar{2}}
    +
    2
    \left(
      n-j+1
    \right)
    \left(
      j-1
    \right)
    \hbar
    R_{\bar{2} \ \ \, \bar{1}}^{\ \, \bar{2}\bar{1}},
    \\
    X_{j+1,j}^{n}
    &=
    2
    \tbinom{j}{2}
    \hbar
    R_{\bar{2} \ \ \, \bar{1}}^{\ \, \bar{2}\bar{2}}
    +
    j
    \left(
      n-j
    \right)
    \hbar
    R_{\bar{1} \ \ \, \bar{1}}^{\ \, \bar{2}\bar{1}},
    \\
    X_{j+2,j}^{n}
    &=
    \tbinom{j+1}{2}
    \hbar
    R_{\bar{1} \ \ \, \bar{1}}^{\ \, \bar{2}\bar{2}},
    \\
    X_{j,k}^{n}
    &=
    0
    \
    \left(
      \left|
        j-k
      \right|
      >
      2
    \right),
  \end{align*}
  \noindent
  where 
  \(
    \tbinom{m}{n}
  \) 
  is a binomial coefficient.
  \label{prop:general}
\end{theorem}

\begin{prf*}
  Note that, 
  \(
    \overrightarrow{\alpha_{n}}
    =
    \left(
      \alpha_{1}^{n},
      \alpha_{2}^{n}
    \right)
  \) 
  and 
  \(
    \overrightarrow{\beta_{n}}
    =
    \left(
      \beta_{1}^{n},
      \beta_{2}^{n}
    \right)
  \) 
  satisfy 
  \(
    \alpha_{1}^{n}
    +
    \alpha_{2}^{n}
    =
    \beta_{1}^{n}
    +
    \beta_{2}^{n}
    =
    n
  \). 
  All possible combinations of 
  \(
  T_{\overrightarrow{\alpha_{n}},\overrightarrow{\beta_{n}^{\ast}}}^{n}
  \)
  are
  \(
    \left(
      n+1
    \right)^{2}
  \) 
  ways :
  \begin{align*}
    \begin{array}{ccc}
      T_{
        \left(n,0\right),
        \left(n,0\right)
      }^{n}
      &
      \cdots
      &
      T_{
        \left(n,0\right),
        \left(0,n\right)
      }^{n},
      \\
      \vdots
      &
      \ddots
      &
      \vdots
      \\
      T_{
        \left(0,n\right),
        \left(n,0\right)
      }^{n}
      &
      \cdots
      &
      T_{
        \left(0,n\right),
        \left(0,n\right)
      }^{n}.
    \end{array}
  \end{align*}
  
  \noindent
  The recurrence relations for these 
  \(
  T_{\overrightarrow{\alpha_{n}},\overrightarrow{\beta_{n}^{\ast}}}^{n}
  \) 
  satisfies
  \begin{flalign}
    &
    \sum_{d=1}^{2}
    \hbar
    g_{\overline{i}d}
    T_{\overrightarrow{\alpha_{n}}-\overrightarrow{e_{d}},\overrightarrow{\beta_{n}^{\ast}}-\overrightarrow{e_{i}}}^{n-1}
    &
    \nonumber
    \\
    &=
    \beta_{i}^{n}
    T_{\overrightarrow{\alpha_{n}},\overrightarrow{\beta_{n}^{\ast}}}^{n}
    +
    \sum_{k=1}^{2}
    \sum_{c=1}^{2}
    \frac{
      \hbar
      \left(\beta_{k}^{n}-\delta_{kc}-\delta_{ik}+1\right)
      \left(\beta_{k}^{n}-\delta_{kc}-\delta_{ik}+2\right)
    }{2}
    R_{\overline{c} \ \ \ \,\overline{i}}^{\ \,\overline{k}\,\overline{k}}
    T_{\overrightarrow{\alpha_{n}},\overrightarrow{\beta_{n}^{\ast}}-\overrightarrow{e_{c}}+2\overrightarrow{e_{k}}-\overrightarrow{e_{i}}}^{n}
    &
    \nonumber
    \\
    &
    \mbox{}
    \quad
    +
    \sum_{c=1}^{2}
    \hbar
    \left(
      \beta_{1}^{n}
      -
      \delta_{1c}
      -
      \delta_{i1}
      +
      1
    \right)
    \left(
      \beta_{2}^{n}
      -
      \delta_{2c}
      -
      \delta_{i2}
      +
      1
    \right)
    R_{\overline{c} \ \ \ \,\overline{i}}^{\ \,\overline{2}\,\overline{1}}
    T_{\overrightarrow{\alpha_{n}},\overrightarrow{\beta_{n}^{\ast}}-\overrightarrow{e_{c}}+\overrightarrow{e_{1}}+\overrightarrow{e_{2}}-\overrightarrow{e_{i}}}^{n}
    &
    \label{rec_rel}
\end{flalign}

  \noindent
  from Theorem \ref{theorem:star_prd}, and two types of them for \(i=1,2\) exist. Both sides of these recurrence relations are linear combinations of some 
  \(
    T_{\overrightarrow{\alpha_{n}},\overrightarrow{\beta_{n}^{\ast}}}^{n}.
  \) 
  So we now make one new recurrence relation below by summing over the index \(i\) on both sides of them.
  \begin{flalign}
    &
    \sum_{i,d=1}^{2}
    \hbar
    g_{\overline{i}d}
    T_{\overrightarrow{\alpha_{n}}-\overrightarrow{e_{d}},\overrightarrow{\beta_{n}^{\ast}}-\overrightarrow{e_{i}}}^{n-1}
    &
    \nonumber
    \\
    &=
    \sum_{i=1}^{2}
    \beta_{i}^{n}
    T_{\overrightarrow{\alpha_{n}},\overrightarrow{\beta_{n}^{\ast}}}^{n}
    +
    \sum_{i,k,c=1}^{2}
    \frac{
      \hbar
      \left(\beta_{k}^{n}-\delta_{kc}-\delta_{ik}+1\right)
      \left(\beta_{k}^{n}-\delta_{kc}-\delta_{ik}+2\right)
    }{2}
    R_{\overline{c} \ \ \ \,\overline{i}}^{\ \,\overline{k}\,\overline{k}}
    T_{\overrightarrow{\alpha_{n}},\overrightarrow{\beta_{n}^{\ast}}-\overrightarrow{e_{c}}+2\overrightarrow{e_{k}}-\overrightarrow{e_{i}}}^{n}
    &
    \nonumber
    \\
    &
    \mbox{}
    \quad
    +
    \sum_{i,c=1}^{2}
    \hbar
    \left(
      \beta_{1}^{n}
      -
      \delta_{1c}
      -
      \delta_{i1}
      +
      1
    \right)
    \left(
      \beta_{2}^{n}
      -
      \delta_{2c}
      -
      \delta_{i2}
      +
      1
    \right)
    R_{\overline{c} \ \ \ \,\overline{i}}^{\ \,\overline{2}\,\overline{1}}
    T_{\overrightarrow{\alpha_{n}},\overrightarrow{\beta_{n}^{\ast}}-\overrightarrow{e_{c}}+\overrightarrow{e_{1}}+\overrightarrow{e_{2}}-\overrightarrow{e_{i}}}^{n}.
    &
    \label{rec_rel_sum}
\end{flalign}
  
  \noindent
  Using the fact 
  \(
    \alpha_{1}^{n}
    +
    \alpha_{2}^{n}
    =
    n
  \) 
  and 
  \(
    \beta_{1}^{n}
    +
    \beta_{2}^{n}
    =
    n
  \), 
we can 
redefine the coefficients 
\(
  T_{\overrightarrow{\alpha_{n}},\overrightarrow{\beta_{n}^{\ast}}}
\) 
as
\begin{align*}
  T^{n}_{i,j}
  :=
  T_{\overrightarrow{\alpha_{n}},\overrightarrow{\beta_{n}^{\ast}}}^{n},
  \
  i
  :=
  \alpha_{2}^{n}+1,
  \
  j
  :=
  \beta_{2}^{n}+1.
\end{align*}
  \noindent
  Then the 
  recurrence relation \eqref{rec_rel_sum} is rewritten as
  \begin{align}
    \hbar
    \sum_{k,l=1}^{2}
    g_{\overline{k}l}
    T^{n-1}_{i-\delta_{2l},j-\delta_{2k}}
    =
    X_{j-2}^{n}
    T_{i,j-2}^{n}
    +
    X_{j-1}^{n}
    T_{i,j-1}^{n}
    +
    X_{j}^{n}
    T_{i,j}^{n}
    +
    X_{j+1}^{n}
    T_{i,j+1}^{n}
    +
    X_{j+1}^{n}
    T_{i,j+2}^{n},
    \label{rec_rel_lhs_symb}
\end{align}

\noindent
where each 
\(
  X_{j}^{n}
\) 
is given as
\begin{align*}
    X_{j-2}^{n}
    &=
    \tbinom{n-j+3}{2}
    \hbar
    R_{\bar{2} \ \ \, \bar{2}}^{\ \, \bar{1}\bar{1}},
    \\
    X_{j-1}^{n}
    &=
    2
    \tbinom{n-j+2}{2}
    \hbar
    R_{\bar{2} \ \ \, \bar{1}}^{\ \, \bar{1}\bar{1}}
    +
    \left(
      n
      -
      j+2
    \right)
    \left(
      j
      -
      2
    \right)
    \hbar
    R_{\bar{2} \ \ \, \bar{2}}^{\ \, \bar{2}\bar{1}},
    \\
    X_{j}^{n}
    &=
    n
    +
    \tbinom{n-j+1}{2}
    \hbar
    R_{\bar{1} \ \ \, \bar{1}}^{\ \, \bar{1}\bar{1}}
    +
    \tbinom{j-1}{2}
    \hbar
    R_{\bar{2} \ \ \, \bar{2}}^{\ \, \bar{2}\bar{2}}
    +
    2
    \left(
      n
      -
      j
      +
      1
    \right)
    \left(
      j-1
    \right)
    \hbar
    R_{\bar{2} \ \ \, \bar{1}}^{\ \, \bar{2}\bar{1}},
    \\
    X_{j+1}^{n}
    &=
    2
    \tbinom{j}{2}
    \hbar
    R_{\bar{2} \ \ \, \bar{1}}^{\ \, \bar{2}\bar{2}}
    +
    j
    \left(
      n
      -
      j
    \right)
    \hbar
    R_{\bar{1} \ \ \, \bar{1}}^{\ \, \bar{2}\bar{1}},
    \\
    X_{j+2}^{n}
    &=
    \tbinom{j+1}{2}
    \hbar
    R_{\bar{1} \ \ \, \bar{1}}^{\ \, \bar{2}\bar{2}},
  \end{align*}

\noindent
respectively. We put 
that 
\(
  X_{j}^{n}
  :=
  0
\) 
when \(j<0\) or \(j>n+1\).
Here, introducing 
each 
\(
  X_{k,j}^{n}
\) 
for 
\(
  k
  =
  1,
  \cdots,
  j,
  \cdots,
  n+1
\) 
by
\begin{flalign*}
  &
  \left(
      \begin{array}{ccccccccccc}
        0
        &
        \cdots
        &
        0
        &
        X_{j-2,j}^{n}
        &
        X_{j-1,j}^{n}
        &
        X_{j,j}^{n}
        &
        X_{j+1,j}^{n}
        &
        X_{j+2,j}^{n}
        &
        0
        &
        \cdots
        &
        0
      \end{array}
  \right)^{T}
  &
  \\
  &
  :=
  \left(
      \begin{array}{ccccccccccc}
        0
        &
        \cdots
        &
        0
        &
        X_{j-2}^{n}
        &
        X_{j-1}^{n}
        &
        X_{j}^{n}
        &
        X_{j+1}^{n}
        &
        X_{j+2}^{n}
        &
        0
        &
        \cdots
        &
        0
      \end{array}
  \right)^{T}
  ,
  &
\end{flalign*}

\noindent
then 
the right-hand side of the recurrence relation \eqref{rec_rel_lhs_symb} can be written as
  \begin{align*}
    \left(r.h.s.\right)
    &=
    \left(
      \begin{array}{ccc}
        T^{n}_{i,1}
        &
        \cdots
        &
        T^{n}_{i,n+1}
      \end{array}
    \right)
    \left(
      \begin{array}{ccccccccccc}
        0
        &
        \cdots
        &
        0
        &
        X_{j-2,j}^{n}
        &
        X_{j-1,j}^{n}
        &
        X_{j,j}^{n}
        &
        X_{j+1,j}^{n}
        &
        X_{j+2,j}^{n}
        &
        0
        &
        \cdots
        &
        0
      \end{array}
    \right)^{T}.
  \end{align*}
  \noindent
  Furthermore, to summarize the recurrence relation \eqref{rec_rel_lhs_symb} by a matrix representation, we introduce \(A_{ij}^{'n}\) for the left-hand side and 
  \(
    \mathbf{X}_{j}
  \)
  for the right-hand side as 
  \begin{align*}
    A_{ij}^{'n}
    :=&
    \sum_{k,l=1}^{2}
    g_{\overline{k}l}
    T^{n-1}_{i-\delta_{2l},j-\delta_{2k}}
    \\
    \mathbf{X}_{j}
    :=&
    \left(
      \begin{array}{ccccccccccc}
        0
        &
        \cdots
        &
        0
        &
        X_{j-2,j}^{n}
        &
        X_{j-1,j}^{n}
        &
        X_{j,j}^{n}
        &
        X_{j+1,j}^{n}
        &
        X_{j+2,j}^{n}
        &
        0
        &
        \cdots
        &
        0
      \end{array}
    \right)^{T},
  \end{align*}
  \noindent
  respectively. Thus, by summarizing the recurrence relation \eqref{rec_rel_lhs_symb} for each \(i\) and using a matrix representation, we have
  \begin{align*}
    \hbar
    A'_{n}
    =
    T_{n}
    X_{n}
  \end{align*}
  \noindent
  for the coefficient 
  \(
    T^{n}_{i,j}
    =
    T_{\overrightarrow{\alpha_{n}},\overrightarrow{\beta_{n}^{\ast}}}^{n},
  \) 
  where 
  \(
    T_{n}
    :=
    \left(
      T^{n}_{ij}
    \right),
    \
    A'_{n}
    =
    \left(A_{ij}^{'n}\right)
    \in
    M_{n+1}
    \left(
    \mathbb{C}
    \left[\!
      \left[
        \hbar
      \right]\!
    \right]
    \right)
  \), 
  and
  \begin{align*}
    X_{n}
    &=
    \left(
      \mathbf{X}_{1},
      \cdots,
      \mathbf{X}_{n+1}
    \right)
    \\
    &=
    \left(
    \left.
      \begin{array}{ccccccccccc}
        X_{11}^{n} & X_{12}^{n} & X_{13}^{n} & 0 & 0 & & & & & & \text{\huge{0}} \\
        X_{21}^{n} & X_{22}^{n} & X_{23}^{n} & X_{24}^{n} & 0 & \ddots & & & & & \\
        X_{31}^{n} & X_{32}^{n} & X_{33}^{n} & X_{34}^{n} & X_{35}^{n} & \ddots & \ddots & & & & \\
        0          & X_{22}^{n} & X_{43}^{n} & X_{44}^{n} & X_{45}^{n} & \ddots & \ddots & \ddots & & & \\
        0          & 0          & X_{53}^{n} & X_{54}^{n} & X_{55}^{n} & \ddots & \ddots & \ddots & \ddots & & \\
                              & \ddots     & \ddots     & \ddots     & \ddots     & \ddots & \ddots & \ddots & \ddots & \ddots & \\
                              & & \ddots     & \ddots & \ddots & \ddots & X_{n-3,n-3}^{n} & X_{n-3,n-2}^{n} & X_{n-3,n-1}^{n} & 0 & 0 \\
                   &  &  & \ddots & \ddots & \ddots & X_{n-2,n-3}^{n} & X_{n-2,n-2}^{n} & X_{n-2,n-1}^{n} & X_{n-2,n}^{n} & 0 \\
                   &            &            &        & \ddots & \ddots & X_{n-1,n-3}^{n} & X_{n-1,n-2}^{n} & X_{n-1,n-1}^{n} & X_{n-1,n}^{n} & X_{n-1,n+1}^{n} \\
                   &            &            & & & \ddots & 0 & X_{n,n-2}^{n} & X_{n,n-1}^{n} & X_{n,n}^{n} & X_{n,n+1}^{n} \\
        \text{\huge{0}} &            & & & & & 0 & 0 & X_{n+1,n-1}^{n} & X_{n+1,n}^{n} & X_{n+1,n+1}^{n}
      \end{array}
    \right.
    \right).
  \end{align*}
  
  \noindent
  This proof was completed.\(_{\blacksquare}\)
\end{prf*}

\noindent
\ Here, we note 
the fact with respect to \(X_{k}\) 
that each \(X_{k}^{-1}\), the inverse matrix of each \(X_{k}\), is determined by a formal power series with respect to a matrix 
\(
  H_{k}
  \in
  M_{k+1}
  \left(
    \mathbb{C}
  \right).
\) 
Here \(H_{k}\) is a pentadiagonal matrix such that each component depends on Riemann curvature tensors on \(M\).
We 
decompose \(X_{k}\) into the two matrices:
\begin{align*}
  X_{k}
  =
  k
  \mathrm{Id}_{k+1}
  +
  \hbar
  H_{k},
\end{align*}

\noindent
where 
\(
  H_{k}
\)
 is given as follows :
\begin{align*}
    H_{j-2,j}^{k}
    &=
    \tbinom{k-j+3}{2}
    R_{\bar{2} \ \ \, \bar{2}}^{\ \, \bar{1}\bar{1}},
    \\
    H_{j-1,j}^{k}
    &=
    2
    \tbinom{k-j+2}{2}
    R_{\bar{2} \ \ \, \bar{1}}^{\ \, \bar{1}\bar{1}}
    +
    \left(
      k-j+2
    \right)
    \left(
      j-2
    \right)
    R_{\bar{2} \ \ \, \bar{2}}^{\ \, \bar{2}\bar{1}},
    \\
    H_{j,j}^{k}
    &=
    \tbinom{k-j+1}{2}
    R_{\bar{1} \ \ \, \bar{1}}^{\ \, \bar{1}\bar{1}}
    +
    \tbinom{j-1}{2}
    R_{\bar{2} \ \ \, \bar{2}}^{\ \, \bar{2}\bar{2}}
    +
    2
    \left(
      k-j+1
    \right)
    \left(
      j-1
    \right)
    R_{\bar{2} \ \ \, \bar{1}}^{\ \, \bar{2}\bar{1}},
    \\
    H_{j+1,j}^{k}
    &=
    2
    \tbinom{j}{2}
    R_{\bar{2} \ \ \, \bar{1}}^{\ \, \bar{2}\bar{2}}
    +
    j
    \left(
      k-j
    \right)
    R_{\bar{1} \ \ \, \bar{1}}^{\ \, \bar{2}\bar{1}},
    \\
    H_{j+2,j}^{k}
    &=
    \tbinom{j+1}{2}
    R_{\bar{1} \ \ \, \bar{1}}^{\ \, \bar{2}\bar{2}},
    \\
    H_{j,l}^{k}
    &=
    0
    \
    \left(
      \left|
        j-l
      \right|
      >
      2
    \right).
  \end{align*}
%
%
%
%
\noindent
Then, \(X_{k}^{-1}\) is given by
\begin{align}
  X_{k}^{-1}
  =&
  \sum_{p=0}^{\infty}
  \frac{
    \left(
      -\hbar
    \right)^{p}
  }{
    k^{p+1}
  }
  \left(H_{k}\right)^{p}.
  \label{Xk_inv_fps}
\end{align}

\noindent
Note that, since a power of a pentadiagonal matrix is not a pentadiagonal matrix in general, \(X_{k}\) is pentadiagonal but \(X_{k}^{-1}\) is not always pentadiagonal.
\bigskip
\\
\ For each 
\(
  n
  \in
  \mathbb{Z}_{\geq0}
\), 
\(
  T_{n}
  \in
  M_{n+1}
  \left(
  \mathbb{C}
  \left[\!
    \left[
      \hbar
    \right]\!
  \right]
  \right),
\) 
\eqref{main_them_N2_T} is also a square matrix of order \(n+1\). 
These matrices are not a unified expression, as the size of matrices depends on \(n\), then it is inconvenient to solve the general term. This problem can be solved by embedding ``finite-dimensional" matrices into ``infinite-dimensional" matrices, and such a procedure provides a unified expression for \(n\geq2\). 
For matrices 
\(
  A
  \in
  M_{n}
  \left(
    \mathbb{C}
    \left[\!
      \left[
        \hbar
      \right]\!
    \right]
  \right)
\) 
the embedding of \(A\) is carried out so that the component whose row or column is greater than \(n\) is 0. That is, the embedding is
\begin{align}
  A
  =
  \left(
    \begin{array}{ccc}
      a_{11}  &  \cdots  &  a_{1n}  \\
      \vdots  &  \ddots  &  \vdots  \\
      a_{n1}  &  \cdots  &  a_{nn}
    \end{array}
  \right)
  \mapsto
  A\left(\infty\right)
  :=
  \left(
    \begin{array}{ccc|ccc}
      a_{11}  &  \cdots  &  a_{1n} &   &  \\
      \vdots  &  \ddots  &  \vdots &  &  \text{\huge{0}}  & \\
      a_{n1}  &  \cdots  &  a_{nn}  &  &  &  \\ \hline
        &  &  &  &  & \\
        &  \text{\huge{0}}  &  &  &  \text{\huge{0}}  &  \\
        &  &  &  &  &  \\
    \end{array}
  \right)
  =
  \left(
    \begin{array}{c|c}
      A & 0 \\  \hline
      0 & 0
    \end{array}
  \right).
  \label{emb_gene}
\end{align}

\noindent
%
In the following calculations, we use \(\mathrm{Id}_{n}\), \(\mathrm{Id}\left(\infty\right)\), \(F_{d,n}\), \(F_{d}\left(\infty\right)\) \(F_{r,n}\) and \(F_{r}\left(\infty\right)\) given by 
\begin{align}
  &
  \mathrm{Id}_{n}
  =
  \left(
    \begin{array}{ccc}
      1 & & 0 \\
       & \ddots &  \\
      0 & & 1
    \end{array}
  \right),
  &
  &
  \mathrm{Id}\left(\infty\right)
  =
  \left(
    \begin{array}{ccc}
       1 & 0 & \\
       0 & 1 & \\
         &  & \ddots
    \end{array}
  \right),
  &
  \label{emb_id}
  \\
  &
  F_{d,n}
  =
  \left(
    \begin{array}{cccc}
      0 &  & & \\
      1 & 0 &  & \\
       & \ddots & \ddots \\
       & &  1  &  0
    \end{array}
  \right),
  &
  &
  F_{d}\left(\infty\right)
  =
  \left(
    \begin{array}{cccc}
      0 &  & & \\
      1 & 0 &  & \\
       & 1 & 0 \\
       & &  \ddots  &  \ddots
    \end{array}
  \right),
  &
  \label{emb_Fd}
  \\
  &
  F_{r,n}
  =
  F_{d,n}^{T},
  &
  &
  F_{r}\left(\infty\right)
  =
  F_{d}\left(\infty\right)^{T}.
  &
  \label{emb_Fr}
\end{align}
%
%
%
%
%
%
%

\begin{remark}
  Note that multiplying \(F_{d,n}\) from the left corresponds to 
  ``shifting the components of the matrix downward by one position, with zeros appearing in the top row". And multiplying \(F_{r,n}\) from the right corresponds to ``shifting the components of the matrix rightward by one position, with zeros appearing in the first column". 
  \(F_{d}\left(\infty\right)\) and \(F_{r}\left(\infty\right)\) correspond 
  to their infinite-dimensional versions.
  \label{rem_Fd_Fr}
\end{remark}

\ 
From Proposition \ref{prop:covar_const_0and1}, the embeddings of 
\(
  T_{0}\left(\infty\right)
\) 
and 
\(
  T_{1}\left(\infty\right)
\) 
are 
\begin{align*}
  T_{0}
  =
  \left(1\right)
  &
  \mapsto
  T_{0}\left(\infty\right)
  :=
  \left(
  \begin{array}{c|c}
    T_{0} & 0 \\  \hline
    0 & 0
  \end{array}
  \right)
  =
  \left(
  \begin{array}{ccc}
    1 & 0 & \cdots \\
    0 & 0 & \\
    \vdots & & \ddots
  \end{array}
  \right),
  \\
  T_{1}
  =
  \hbar
  \left(
    \begin{array}{cc}
      g_{1\bar{1}}
      &
      g_{1\bar{2}}
      \\
      g_{2\bar{1}}
      &
      g_{2\bar{2}}
    \end{array}
  \right)
  &
  \mapsto
  T_{1}\left(\infty\right)
  :=
  \left(
  \begin{array}{c|c}
    T_{1} & 0 \\  \hline
    0 & 0
  \end{array}
  \right)
  =
  \hbar
  \left(
    \begin{array}{cccc}
      g_{1\bar{1}}
      &
      g_{1\bar{2}}
      &
      0
      &
      \cdots
      \\
      g_{2\bar{1}}
      &
      g_{2\bar{2}}
      &
      0
      &
      \\
      0
      &
      0
      &
      0
      &
      \\
      \vdots
      &
      &
      &
      \ddots
    \end{array}
  \right),
\end{align*}
\noindent
respectively. Similarly, \eqref{main_them_N2_T} in Theorem \ref{prop:general} is embedded as
\begin{align}
  T_{n}\left(\infty\right)
  =
  \hbar
  A'_{n}\left(\infty\right)
  X^{-1}_{n}\left(\infty\right),
  \label{T_emb_inf}
\end{align}

\noindent
where \(A'_{n}\left(\infty\right)\) and 
\(X^{-1}_{n}\left(\infty\right)\) are
\begin{align*}
  A'_{n}\left(\infty\right)
  :=
  \left(
  \begin{array}{c|c}
    A'_{n} & 0 \\  \hline
    0 & 0
  \end{array}
  \right),
  \qquad
  X^{-1}_{n}\left(\infty\right)
  :=
  \left(
  \begin{array}{c|c}
    X^{-1}_{n} & 0 \\  \hline
    0 & 0
  \end{array}
  \right).
\end{align*}

\noindent
Note that \(X_{n}^{-1}\left(\infty\right)\) is not the 
inverse matrix of \(X_{n}\left(\infty\right)\). As we saw in Theorem \ref{prop:general}, 
\(
  A'_{n}
  \in
  M_{n+1}
  \left(
    \mathbb{C}
    \left[\!
      \left[
        \hbar
      \right]\!
    \right]
  \right)
\) 
is 
expressed as
\begin{align*}
  A'_{n}
  =
  \left(
    \sum_{k,l=1}^{2}
    g_{\overline{k}l}
    T^{n-1}_{i-\delta_{2l},j-\delta_{2k}}
  \right)
  =
  \left(
    g_{\overline{1}1}
    T^{n-1}_{i,j}
    +
    g_{\overline{1}2}
    T^{n-1}_{i-1,j}
    +
    g_{\overline{2}1}
    T^{n-1}_{i,j-1}
    +
    g_{\overline{2}2}
    T^{n-1}_{i-1,j-1}
  \right).
\end{align*}

\noindent
Recall 
that \(T^{n-1}_{i,j}=0\)  if 
\(
  i<0,
  \
  i>n-1,
  \
  j<0
\) 
or 
\(
  j>n-1
\) 
by definition. 
By a matrix representation
\begin{align*}
  &
  \left(
    g_{\overline{1}1}
    T^{n-1}_{i,j}
    +
    g_{\overline{1}2}
    T^{n-1}_{i-1,j}
    +
    g_{\overline{2}1}
    T^{n-1}_{i,j-1}
    +
    g_{\overline{2}2}
    T^{n-1}_{i-1,j-1}
  \right)
  \\
  &=
  g_{\overline{1}1}
  \left(
  \begin{array}{c|c}
    T_{n-1} & \mathbf{0}_{n-1} \\  \hline
    \mathbf{0}_{n-1}^{T} & 0
  \end{array}
  \right)
  +
  g_{\overline{1}2}
  \left(
  \begin{array}{c|c}
    \mathbf{0}_{n-1}^{T} & 0 \\  \hline
    T_{n-1} & \mathbf{0}_{n-1}
  \end{array}
  \right)
  +
  g_{\overline{2}1}
  \left(
  \begin{array}{c|c}
    \mathbf{0}_{n-1} & T_{n-1} \\  \hline
    0 & \mathbf{0}_{n-1}^{T}
  \end{array}
  \right)
  +
  g_{\overline{2}2}
  \left(
  \begin{array}{c|c}
    0 & \mathbf{0}_{n-1}^{T} \\  \hline
    \mathbf{0}_{n-1} & T_{n-1}
  \end{array}
  \right)
\end{align*}

\noindent
with \(T_{n-1}\) which 
is the square matrix of order \(n\), 
where \(\mathbf{0}_{n-1}\) is an \(n-1\) dimensional zero vector. Since these matrices can be expressed as

\begin{align*}
  \left(
  \begin{array}{c|c}
    T_{n-1} & \mathbf{0}_{n-1} \\  \hline
    \mathbf{0}_{n-1}^{T} & 0
  \end{array}
  \right)
  =&
  \mathrm{Id_{n+1}}
  \left(
  \begin{array}{c|c}
    T_{n-1} & \mathbf{0}_{n-1} \\  \hline
    \mathbf{0}_{n-1}^{T} & 0
  \end{array}
  \right)
  \mathrm{Id_{n+1}},
  \\
  \left(
  \begin{array}{c|c}
    \mathbf{0}_{n-1}^{T} & 0 \\  \hline
    T_{n-1} & \mathbf{0}_{n-1}
  \end{array}
  \right)
  =&
  F_{d,n+1}
  \left(
  \begin{array}{c|c}
    T_{n-1} & \mathbf{0}_{n-1} \\  \hline
    \mathbf{0}_{n-1}^{T} & 0
  \end{array}
  \right)
  \mathrm{Id_{n+1}},
  \\
  \left(
  \begin{array}{c|c}
    \mathbf{0}_{n-1} & T_{n-1} \\  \hline
    0 & \mathbf{0}_{n-1}^{T}
  \end{array}
  \right)
  =&
  \mathrm{Id_{n+1}}
  \left(
  \begin{array}{c|c}
    T_{n-1} & \mathbf{0}_{n-1} \\  \hline
    \mathbf{0}_{n-1}^{T} & 0
  \end{array}
  \right)
  F_{r,n+1},
  \\
  \left(
  \begin{array}{c|c}
    0 & \mathbf{0}_{n-1}^{T} \\  \hline
    \mathbf{0}_{n-1} & T_{n-1}
  \end{array}
  \right)
  =&
  F_{d,n+1}
  \left(
  \begin{array}{c|c}
    T_{n-1} & \mathbf{0}_{n-1} \\  \hline
    \mathbf{0}_{n-1}^{T} & 0
  \end{array}
  \right)
  F_{r,n+1},
\end{align*}

\noindent
using \(\mathrm{Id}_{n+1},\) \(F_{d,n+1}\) and \(F_{r,n+1},\) respectively, then we have

\begin{align*}
  &
  g_{\overline{1}1}
  \left(
  \begin{array}{c|c}
    T_{n-1} & \mathbf{0}_{n-1} \\  \hline
    \mathbf{0}_{n-1}^{T} & 0
  \end{array}
  \right)
  +
  g_{\overline{1}2}
  \left(
  \begin{array}{c|c}
    \mathbf{0}_{n-1}^{T} & 0 \\  \hline
    T_{n-1} & \mathbf{0}_{n-1}
  \end{array}
  \right)
  +
  g_{\overline{2}1}
  \left(
  \begin{array}{c|c}
    \mathbf{0}_{n-1} & T_{n-1} \\  \hline
    0 & \mathbf{0}_{n-1}^{T}
  \end{array}
  \right)
  +
  g_{\overline{2}2}
  \left(
  \begin{array}{c|c}
    0 & \mathbf{0}_{n-1}^{T} \\  \hline
    \mathbf{0}_{n-1} & T_{n-1}
  \end{array}
  \right)
  \\
  &=
  g_{\overline{1}1}
  \mathrm{Id_{n+1}}
  \left(
  \begin{array}{c|c}
    T_{n-1} & \mathbf{0}_{n-1} \\  \hline
    \mathbf{0}_{n-1}^{T} & 0
  \end{array}
  \right)
  \mathrm{Id_{n+1}}
  +
  g_{\overline{1}2}
  F_{d,n+1}
  \left(
  \begin{array}{c|c}
    T_{n-1} & \mathbf{0}_{n-1} \\  \hline
    \mathbf{0}_{n-1}^{T} & 0
  \end{array}
  \right)
  \mathrm{Id_{n+1}}
  \\
  &
  \mbox{}
  \quad
  +
  g_{\overline{2}1}
  \mathrm{Id_{n+1}}
  \left(
  \begin{array}{c|c}
    T_{n-1} & \mathbf{0}_{n-1} \\  \hline
    \mathbf{0}_{n-1}^{T} & 0
  \end{array}
  \right)
  F_{r,n+1}
  +
  g_{\overline{2}2}
  F_{d,n+1}
  \left(
  \begin{array}{c|c}
    T_{n-1} & \mathbf{0}_{n-1} \\  \hline
    \mathbf{0}_{n-1}^{T} & 0
  \end{array}
  \right)
  F_{r,n+1}.
\end{align*}

\noindent
By introducing some functions 
\(
  \Theta_{\bar{\mu}}^{p},
  \Theta_{\nu}^{p}
  \in
  C^{\infty}
  \left(
    U
  \right)
  \
  \left(
    p
    \in
    \left\{
      1,
      \cdots,
      4
    \right\}
  \right)
\) 
such that 
\(
  g_{\overline{\mu}\nu}
  =
  \Theta_{\bar{\mu}}^{p}
  \Theta_{\nu}^{p},
\) 
\(A'_{n}\) can be expressed 
as
\begin{align*}
  A'_{n}
  =
  \sum_{p=1}^{4}
  \left(
    \Theta_{1}^{p}
    \mathrm{Id}_{n+1}
    +
    \Theta_{2}^{p}
    F_{d,n+1}
  \right)
  \left(
  \begin{array}{c|c}
    T_{n-1} & \mathbf{0}_{n-1} \\  \hline
    \mathbf{0}_{n-1}^{T} & 0
  \end{array}
  \right)
  \left(
    \Theta_{\bar{1}}^{p}
    \mathrm{Id}_{n+1}
    +
    \Theta_{\bar{2}}^{p}
    F_{r,n+1}
  \right).
\end{align*}
\noindent
These 
\(
  \Theta_{\bar{\mu}}^{p}
\) 
and 
\(
  \Theta_{\nu}^{p}
\) 
play a similar role of 
``a vierbein", but not a vierbein because they are not necessarily to be orthonormal. Furthermore, by introducing the new matrices 
\(
  F^{p}_{n},
  \
  B^{p}_{n}
  \in
  M_{n+1}
  \left(
  \mathbb{C}
  \left[\!
    \left[
      \hbar
    \right]\!
  \right]
  \right)
\) 
as
\begin{align*}
  F^{p}_{n}
  :=
  \Theta_{1}^{p}
  \mathrm{Id}_{n+1}
  +
  \Theta_{2}^{p}
  F_{d,n+1},
  \qquad
  B^{p}_{n}
  :=
  \Theta_{\bar{1}}^{p}
  \mathrm{Id}_{n+1}
  +
  \Theta_{\bar{2}}^{p}
  F_{r,n+1}
\end{align*}
\noindent
respectively, \(A'_{n}\) is expressed as
\begin{align*}
  A'_{n}
  =
  \sum_{p=1}^{4}
  F^{p}_{n}
  \left(
  \begin{array}{c|c}
    T_{n-1} & \mathbf{0}_{n-1} \\  \hline
    \mathbf{0}_{n-1}^{T} & 0
  \end{array}
  \right)
  B^{p}_{n},
\end{align*}

\noindent
so \(A'_{n}\left(\infty\right),\) that is, the embedding matrix for \(A'_{n}\), is expressed as 
\begin{align*}
  A'_{n}\left(\infty\right)
  &=
  \sum_{p=1}^{4}
  F^{p}_{n}\left(\infty\right)
  T_{n-1}\left(\infty\right)
  B^{p}_{n}\left(\infty\right).
\end{align*}

\noindent
Here, \(F^{p}_{n}\) and \(B^{p}_{n}\) were replaced by \(F^{p}_{n}\left(\infty\right)\) and \(B^{p}_{n}\left(\infty\right)\), where

\begin{align}
  F^{p}_{n}\left(\infty\right)
  :=
  \Theta_{1}^{p}
  \mathrm{Id}\left(\infty\right)
  +
  \Theta_{2}^{p}
  F_{d}\left(\infty\right),
  \qquad
  B^{p}_{n}\left(\infty\right)
  :=
  \Theta_{\bar{1}}^{p}
  \mathrm{Id}\left(\infty\right)
  +
  \Theta_{\bar{2}}^{p}
  F_{r}\left(\infty\right).
  \label{emb_Fn_Bn}
\end{align}

\noindent
Therefore, \(T_{n}\left(\infty\right)\) 
%
is obtained as
\begin{align*}
  T_{n}\left(\infty\right)
  =
  \hbar
  \sum_{p=1}^{4}
  F^{p}_{n}\left(\infty\right)
  T_{n-1}\left(\infty\right)
  B^{p}_{n}\left(\infty\right)
  X^{-1}_{n}\left(\infty\right).
\end{align*}
\noindent
By using recursively the above procedure, 
we obtain the following main theorem.

\begin{theorem}[Main result]
  Let \(M\) be a complex two-dimensional locally symmetric K\"{a}hler manifold, \(U\) be an open set of \(M\), and \(\ast\) be a star product with separation of variables on \(M\). 
  For 
  \(
    f,g
    \in
    C^{\infty}
    \left(U\right),
    \
    f
    \ast
    g
  \) 
  is given by
  \begin{align*}
    f
    \ast
    g
    =
    \sum_{n=0}^{\infty}
    \sum_{i,j=1}^{n+1}
    T_{ij}^{n}
    \left\{
    \left(
      D^{1}
    \right)^{n-i+1}
    \left(
      D^{2}
    \right)^{i-1}
    f
    \right\}
    \left\{
    \left(
      D^{\bar{1}}
    \right)^{n-j+1}
    \left(
      D^{\bar{2}}
    \right)^{j-1}
    g
    \right\}.
  \end{align*}
  
  \noindent
  Here 
  each of the coefficient 
   \(
     T_{ij}^{n}
   \)
   is
  \begin{align*}
     {T_{n}\left(\infty\right)}_{ij}
     =
     \begin{cases}
     T_{ij}^{n}
    &
    \
    \left(
      1
      \leq
      i,j
      \leq
      n+1
    \right)
    \\
    0
    &
    \
    \left(
      \mathrm{otherwise}
    \right)
  \end{cases}
  ,
\end{align*}

  \noindent
  and \(T_{n}\left(\infty\right)\) 
  is determined by
  \begin{align}
  T_{n}\left(\infty\right)
  &=
  \hbar^{n}
  \sum_{p_{1},\cdots,p_{n}=1}^{4}
  F^{p_{n}}_{n}\left(\infty\right)
  \cdots
  F^{p_{1}}_{1}\left(\infty\right)
  T_{0}\left(\infty\right)
  \left(
    B^{p_{1}}_{1}\left(\infty\right)
    X^{-1}_{1}\left(\infty\right)
  \right)
  \cdots
  \left(
    B^{p_{n}}_{n}\left(\infty\right)
    X^{-1}_{n}\left(\infty\right)
  \right),
  \label{main_them_N2_Tinf}
\end{align}
   \noindent
   where 
   each \(F^{p_{k}}_{k}\left(\infty\right)\), \(B^{p_{k}}_{k}\left(\infty\right)\) and \(X^{-1}_{k}\left(\infty\right)\) are given as above.
  \label{main_thm_N2}
\end{theorem}

\noindent
 From Theorem \ref{main_thm_N2}, we obtain 
 a deformation quantization with separation of variables for complex two-dimensional locally symmetric K\"{a}hler manifold is realized by this star product.

\subsection{Another formula}
\label{NR_minus}
\ The formula \eqref{main_them_N2_Tinf} was obtained by summation with respect to the index \(i\) of the complex coordinate in the recurrence relation \eqref{rec_rel}. In Subsection \ref{NR_plus}, we made one recurrence relation by adding two recurrence relations and determined the solution using only that recurrence relation. On the other hand, we have not yet considered another recurrence relation. 
In this subsection, we consider another formula obtained by ``subtracting" for \(i=1,2\) rather than ``adding" two recurrence relations for \(i=1,2\) as we did in Subsection \ref{NR_plus}.

\begin{theorem}
  Let \(M\) be a locally symmetric K\"{a}hler manifold, \(U\) be an open set of \(M\),
  and 
  \(
    f,g
    \in
    C^{\infty}
    \left(U\right),
    \
    \ast
  \)
  be a star product with separation of variables on \(U\) such that
  \begin{align*}
    f
    \ast
    g
    =&
    \sum_{n=0}^{\infty}
    \sum_{i,j=1}^{n+1}
    T_{ij}^{n}
    \left\{
    \left(
      D^{1}
    \right)^{n-i+1}
    \left(
      D^{2}
    \right)^{i-1}
    f
    \right\}
    \left\{
    \left(
      D^{\bar{1}}
    \right)^{n-j+1}
    \left(
      D^{\bar{2}}
    \right)^{j-1}
    g
    \right\}.
  \end{align*}
  Then
  \begin{align}
    T_{n}
    Y_{n}
    =
    \hbar
    C'_{n}
    \label{minus_thm_N2_eq}
  \end{align}
  or equivalently
  \begin{align}
    Y_{n}^{\dagger}
    T_{n}
    =
    \hbar
    C_{n}^{'\dagger},
    \label{minus_thm_N2_eq_another}
  \end{align}
  where 
  \begin{align*}
    C'_{n}
    &
    =
    \left(C_{ij}^{'n}\right)
    =
    \left(
      \sum_{k,l=1}^{2}
      \left(
        -1
      \right)^{\delta_{k2}}
      g_{\overline{k}l}
      T^{n-1}_{i-\delta_{2l},j-\delta_{2k}}
    \right)
    \in
    M_{n+1}
    \left(
    \mathbb{C}
    \left[\!
      \left[
        \hbar
      \right]\!
    \right]
  \right)
\end{align*}
and 
\(
  Y_{n}
  \in
  M_{n+1}
  \left(
    \mathbb{C}
    \left[\!
      \left[
        \hbar
      \right]\!
    \right]
  \right)
\)
is a pentadiagonal matrix such that its 
components are given as follows :
\begin{align*}
    Y_{j-2,j}^{n}
    &=
    -
    \tbinom{n-j+3}{2}
    \hbar
    R_{\bar{2} \ \ \, \bar{2}}^{\ \, \bar{1}\bar{1}},
    \\
    Y_{j-1,j}^{n}
    &=
    -
    \left(
      n-j+2
    \right)
    \left(
      j-2
    \right)
    \hbar
    R_{\bar{2} \ \ \, \bar{2}}^{\ \, \bar{2}\bar{1}},
    \\
    Y_{j,j}^{n}
    &=
    n
    -
    2j
    +
    2
    +
    \tbinom{n-j+1}{2}
    \hbar
    R_{\bar{1} \ \ \, \bar{1}}^{\ \, \bar{1}\bar{1}}
    -
    \tbinom{j-1}{2}
    \hbar
    R_{\bar{2} \ \ \, \bar{2}}^{\ \, \bar{2}\bar{2}},
    \\
    Y_{j+1,j}^{n}
    &=
    j
    \left(
      n-j
    \right)
    \hbar
    R_{\bar{1} \ \ \, \bar{1}}^{\ \, \bar{2}\bar{1}},
    \\
    Y_{j+2,j}^{n}
    &=
    \tbinom{j+1}{2}
    \hbar
    R_{\bar{1} \ \ \, \bar{1}}^{\ \, \bar{2}\bar{2}},
    \\
    Y_{j,k}^{n}
    &=
    0
    \
    \left(
      \left|
        j-k
      \right|
      >
      2
    \right).
  \end{align*}
  \label{minus_thm_N2}
\end{theorem}

\begin{prf*}
This is shown in 
the same way as in Theorem \ref{main_thm_N2}.\(_{\blacksquare}\)
\end{prf*}
\noindent
Note that equation \eqref{minus_thm_N2_eq} is also expressed as an embedding version
\begin{align}
  T_{n}\left(\infty\right)
  Y_{n}\left(\infty\right)
  =
  \hbar
  C'_{n}\left(\infty\right)
  \label{minus_thm_N2_eq_rsz}
\end{align}
\noindent
or equivalently
\begin{align}
  Y_{n}^{\dagger}\left(\infty\right)
  T_{n}\left(\infty\right)
  =
  \hbar
  C_{n}^{'\dagger}\left(\infty\right),
  \label{minus_thm_N2_eq_another_rsz}
\end{align}
\noindent
where the embedding \(Y_{n}\left(\infty\right)\) for \(Y_{n}\) is as in the way of \eqref{emb_gene} in Subsection \ref{NR_plus}. 
%
Here we use 
\(
  T_{n}
  =
  T_{n}^{\dagger}.
\) 
Its derivation is in Appendix \ref{prf_Tn_dagger}.
\bigskip
\\
\
It is known from Karabegov's result that there is always a star product with separation of variables on a K\"{a}hler manifold \cite{Kara1,Kara2}. A star product with separation of variables on a locally symmetric K\"{a}hler manifold is 
determined by \eqref{main_them_N2_Tinf} in Theorem \ref{main_thm_N2}. Therefore, \(T_{n}\left(\infty\right)\) given by \eqref{main_them_N2_Tinf} should satisfy \eqref{minus_thm_N2_eq_rsz}. To ensure that the result obtained by Theorem \ref{main_thm_N2} does not contradict Theorem \ref{minus_thm_N2}, we consider the following equation: 
\begin{align}
  Y_{n}^{\dagger}
  T_{n}
  X_{n}
  =
  \hbar
  C_{n}^{'\dagger}
  X_{n}.
  \label{minus_thm_eq_ens}
\end{align}
\noindent
Here, \eqref{minus_thm_eq_ens} is obtained by multiplying \eqref{minus_thm_N2_eq_another} by \(X_{n}\) from the right. This is merely a change to a form that allows direct substitution of the result of Theorem \ref{main_thm_N2}. The reason for using \eqref{minus_thm_N2_eq_another} rather than \eqref{minus_thm_N2_eq} is that \(X^{-1}_{n}\) appears if we try to check \eqref{minus_thm_N2_eq} directly, and it is difficult to calculate because it is an infinite power series matrix. In this discussion, we shall simply show only the cases \(n=1\) and \(2\).
\vspace{3mm}
\\
\noindent
\underline{
  Case:
  \(
    n
    =
    1
  \)
}
\vspace{3mm}
\\
Since Proposition \ref{prop:covar_const_0and1} and
\begin{align*}
  X_{1}
  &=
  \left(
    \begin{array}{cc}
      X_{11}^{1}  &  X_{12}^{1}  \\
      X_{21}^{1}  &  X_{22}^{1}
    \end{array}
  \right)
  =
  \mathrm{Id}_{2},
  \\
  Y_{1}
  &=
  \left(
    \begin{array}{cc}
      Y_{11}^{1}  &  Y_{12}^{1}  \\
      Y_{21}^{1}  &  Y_{22}^{1}
    \end{array}
  \right)
  =
  \left(
    \begin{array}{cc}
      1
      &
      0
      \\
      0
      &
      -
      1
    \end{array}
  \right),
\end{align*}
\noindent
it follows that
\begin{align*}
  Y_{1}^{\dagger}
  T_{1}
  X_{1}
  &=
  \hbar
  \left(
    \begin{array}{cc}
      g_{1\overline{1}}  &  g_{2\overline{1}}  \\
      -g_{1\overline{2}}  &  -g_{2\overline{2}}
    \end{array}
  \right),
  \\
  C_{1}^{\dagger}
  X_{1}
  &=
  \left(
    \begin{array}{cc}
      g_{\overline{1}1}  &  -g_{\overline{2}1}
      \\
      g_{\overline{1}2}  &  -g_{\overline{2}2}
    \end{array}
  \right)^{\dagger}=
  \left(
    \begin{array}{cc}
      g_{1\overline{1}}  &  g_{2\overline{1}}  \\
      -g_{1\overline{2}}  &  -g_{2\overline{2}}
    \end{array}
  \right).
\end{align*}
\noindent
Therefore 
\(
  Y_{1}^{\dagger}
  T_{1}
  X_{1}
  =
  \hbar
  C'_{1}
  X_{1}
\)
, equivalently 
\(
  Y_{1}^{\dagger}\left(\infty\right)
  T_{1}\left(\infty\right)
  X_{1}\left(\infty\right)
  =
  \hbar
  C'_{1}\left(\infty\right)
  X_{1}\left(\infty\right).
\)
\vspace{3mm}
\\
\noindent
\underline{
  Case:
  \(
    n
    =
    2
  \)
}
\vspace{3mm}
\\
Calculations for both sides of \eqref{minus_thm_eq_ens} yields
\begin{align}
  Y_{2}^{\dagger}
  T_{2}
  X_{2}
  =
  \hbar^{2}
  \left(
    \begin{array}{ccc}
      2
      +
      \hbar
      R_{1 \ \ \, 1}^{\ \, 11}
      &
      \hbar
      R_{1 \ \ \, 1}^{\ \, 21}
      &
      \hbar
      R_{1 \ \ \, 1}^{\ \, 22}
      \\
      0
      &
      0
      &
      0
      \\
      -
      \hbar
      R_{2 \ \ \, 2}^{\ \, 11}
      &
      -
      \hbar
      R_{2 \ \ \, 2}^{\ \, 21}
      &
      -2
      -
      \hbar
      R_{2 \ \ \, 2}^{\ \, 22}
    \end{array}
  \right)
  \left(
    \begin{array}{ccc}
      \left(
        g_{\overline{1}1}
      \right)^{2}
      &
      2
      g_{\overline{1}1}
      g_{\overline{2}1}
      &
      \left(
        g_{\overline{2}1}
      \right)^{2}
      \\
      2
      g_{\overline{1}1}
      g_{\overline{1}2}
      &
      2
      \left(
        g_{\overline{1}1}
        g_{\overline{2}2}
        +
        g_{\overline{1}2}
        g_{\overline{2}1}
      \right)
      &
      2
      g_{\overline{2}1}
      g_{\overline{2}2}
      \\
      \left(
        g_{\overline{1}2}
      \right)^{2}
      &
      2
      g_{\overline{1}2}
      g_{\overline{2}2}
      &
      \left(
        g_{\overline{2}2}
      \right)^{2}
    \end{array}
  \right)
\end{align}
\noindent
and
\begin{align}
  \hbar
  C'_{2}
  X_{2}
  =
  \hbar^{2}
  \left(
    \begin{array}{ccc}
      \left(
        g_{\overline{1}1}
      \right)^{2}
      &
      2
      g_{\overline{1}1}
      g_{\overline{2}1}
      &
      \left(
        g_{\overline{2}1}
      \right)^{2}
      \\
      0
      &
      0
      &
      0
      \\
      -
      \left(
        g_{\overline{1}2}
      \right)^{2}
      &
      -2
      g_{\overline{1}2}
      g_{\overline{2}2}
      &
      -
      \left(
        g_{\overline{2}2}
      \right)^{2}
    \end{array}
  \right)
  \left(
    \begin{array}{ccc}
      2
      +
      \hbar
      R_{\bar{1} \ \ \, \bar{1}}^{\ \, \bar{1}\bar{1}}
      &
      2
      \hbar
      R_{\bar{2} \ \ \, \bar{1}}^{\ \, \bar{1}\bar{1}}
      &
      \hbar
      R_{\bar{2} \ \ \, \bar{2}}^{\ \, \bar{1}\bar{1}}
      \\
      \hbar
      R_{\bar{1} \ \ \, \bar{1}}^{\ \, \bar{2}\bar{1}}
      &
      2
      +
      2
      \hbar
      R_{\bar{2} \ \ \, \bar{1}}^{\ \, \bar{2}\bar{1}}
      &
      \hbar
      R_{\bar{2} \ \ \, \bar{2}}^{\ \, \bar{2}\bar{1}}
      \\
      \hbar
      R_{\bar{1} \ \ \, \bar{1}}^{\ \, \bar{2}\bar{2}}
      &
      2
      \hbar
      R_{\bar{2} \ \ \, \bar{1}}^{\ \, \bar{2}\bar{2}}
      &
      2
      +
      \hbar
      R_{\bar{2} \ \ \, \bar{2}}^{\ \, \bar{2}\bar{2}}
    \end{array}
  \right)
\end{align}
\noindent
respectively. Here we use 
\begin{align}
  T_{2}
  X_{2}
  =
  \hbar
  A'_{2}
  =
  \left(
    \begin{array}{ccc}
      \left(
        g_{\overline{1}1}
      \right)^{2}
      &
      2
      g_{\overline{1}1}
      g_{\overline{2}1}
      &
      \left(
        g_{\overline{2}1}
      \right)^{2}
      \\
      2
      g_{\overline{1}1}
      g_{\overline{1}2}
      &
      2
      \left(
        g_{\overline{1}1}
        g_{\overline{2}2}
        +
        g_{\overline{1}2}
        g_{\overline{2}1}
      \right)
      &
      2
      g_{\overline{2}1}
      g_{\overline{2}2}
      \\
      \left(
        g_{\overline{1}2}
      \right)^{2}
      &
      2
      g_{\overline{1}2}
      g_{\overline{2}2}
      &
      \left(
        g_{\overline{2}2}
      \right)^{2}
    \end{array}
  \right).
\end{align}
\noindent
Furthermore, each component 
of both sides of \eqref{minus_thm_eq_ens} coincides with each other from the calculations \eqref{lhs_first}--\eqref{rhs_last} in Appendix \ref{app_minus_n2}. Hence, it is shown that \eqref{minus_thm_eq_ens} for \(n=2\) holds.

\section{Examples}
\label{conc_ex}
\ In this section, we describe the deformation quantization with separation of variables for \(\mathbb{C}^{2}\) and \(\mathbb{C}P^{2}\) as concrete examples realized from Theorem \ref{main_thm_N2}. The deformation quantizations were already given by the other methods as we will see in the next subsection. To compare our results with them, we can check if Theorem \ref{main_thm_N2} works well. Moreover, we show that 
the concrete examples satisfy the identities in Theorem \ref{minus_thm_N2}.

\subsection{Previous results for \(\mathbb{C}^{2}\) and \(\mathbb{C}P^{2}\)}
\label{subsec_prev_C_CP}

\ It is known that noncommutative \(\mathbb{R}^{2N}\) is the most trivial noncommutative manifold. We can say that the quantization map from \(\mathbb{R}^{2N}\) to noncommutative \(\mathbb{R}^{2N}\) is introduced by Dirac \cite{Dirac} because the canonical quantization is the quantization of phase space \(\mathbb{R}^{2N}\). 
In the strict deformation quantization case, the noncommutative plane has also been proposed by Rieffel from the perspective of \(C^{\ast}\)-algebra \cite{Rief2}. On the other hand, the deformation quantization of \(\mathbb{R}^{2N}\) is also the most trivial example in the formal deformation quantization case. In particular, the Moyal product \cite{Moya} and the Voros product \cite{Voro} are well-known examples of star products on \(\mathbb{R}^{2N}\). 
Here, a noncommutative \(\mathbb{R}^{2N}\) can be regarded as a noncommutative \(\mathbb{C}^{N}\). 
For noncommutative \(\mathbb{C}^{N}\), the concrete construction is provided by using the deformation quantization with separation of variables by Karabegov as an example \cite{Kara1}. 
As another concrete example, we shall consider here a deformation quantization with separation of variables for \(\mathbb{C}P^{2}\). Previous works related to noncommutative \(\mathbb{C}P^{N}\) are known from deformation quantizations and fuzzy geometry. Noncommutative \(\mathbb{C}P^{N}\) via a deformation quantizations had constructed by Omori-Maeda-Yoshioka \cite{OMY2}, Bordemann et al. \cite{BBEW}, Sako-Suzuki-Umetsu \cite{SSU1,SSU2}, and Hara-Sako \cite{HS1,HS2}. 
In addition, the (twisted) Fock representation of noncommutative \(\mathbb{C}P^{N}\) was given by Sako-Suzuki-Umetsu \cite{SSU1,SSU2} and Sako-Umetsu \cite{SU1,SU2,SU3}. It is known that this representation is 
essentially equivalent to ordinary Fock representation. 
On the other hand, in fuzzy physics, noncommutative \(\mathbb{C}P^{N}\), called fuzzy \(\mathbb{C}P^{N}\), had constructed via a matrix algebra. 
The construction methods of fuzzy \(S^{2}\), i.e. fuzzy \(\mathbb{C}P^{1}\), had proposed by Hoppe \cite{Hop} and Madore \cite{Mad1}. Hayasaka-Nakayama-Takaya constructed a star product on \(S^{2}\) from their fuzzy \(S^{2}\) \cite{HNT}. 
Furthermore, more general works were proposed by Grosse-Strohmaier \cite{GS}, Alexanian et al. \cite{ABIY} for fuzzy \(\mathbb{C}P^{2}\), and by Balachandran et al. \cite{BDMLO}, Carow-Watamura-Steinacker-Watamura \cite{CSW} for \(\mathbb{C}P^{N}\). See \cite{Mad2} for other references. In particular, the Fock representation on fuzzy \(\mathbb{C}P^{N}\) was also discussed by Alexanian-Pinzul-Stern \cite{APS} and Carow-Watamura-Steinacker-Watamura \cite{CSW}. 
In this subsection, we introduce the results of noncommutative \(\mathbb{C}^{2}\) by Karabegov and Voros, and noncommutative \(\mathbb{C}P^{2}\) by Sako-Suzuki-Umetsu to compare them with results given by Theorem \ref{main_thm_N2} in the following sections.
\begin{theorem}[Karabegov\cite{Kara2}]
For any 
\(
  f,g
  \in
  C^{\infty}
  \left(\mathbb{C}^{2}\right)
\), 
the star product with separation of variables on \(\mathbb{C}^{2}\) is locally given by
\begin{align}
  f
  \ast
  g
  =
  f
  \mathrm{exp}
  \nu
  \left(
    \overleftarrow{
      \partial_{\overline{z}^{1}}
    }
    \overrightarrow{
      \partial_{z^{1}}
    }
    +
    \overleftarrow{
      \partial_{\overline{z}^{2}}
    }
    \overrightarrow{
      \partial_{z^{2}}
    }
  \right)
  g.
  \label{star_prd_C2_Kara_Voro_sol}
\end{align}
\label{star_prd_C2_Kara_Voro}
\end{theorem}
\noindent
It is known as the star product which is called ``Voros product" \cite{Voro} and 
isomorphic to Moyal product \cite{Moya}. 
A simple proof for the isomorphism of them is given by, for example, Section 2.2 in \cite{APS}. 
The star product \eqref{star_prd_C2_Kara_Voro_sol} is sometimes called the ``Wick product" \cite{BBEW,BW1,BW2}. For \(\mathbb{C}P^{N}\), the noncommutative one had been constructed by several previous works. In Section \ref{subsec_NR_CP2}, we compare the star product on \(\mathbb{C}P^{2}\) obtained by Theorem \ref{main_thm_N2} with the previous work by Sako-Suzuki-Umetsu. So we denote the result for \(\mathbb{C}P^{2}\) by them.
\begin{theorem}[Sako-Suzuki-Umetsu\cite{SSU1,SSU2}]
For any 
\(
  f,g
  \in
  C^{\infty}
  \left(\mathbb{C}P^{2}\right)
\), 
the star product with separation of variables on \(\mathbb{C}P^{2}\) is locally given by
\begin{align}
  f
  \ast
  g
  =
  \sum_{n=0}^{\infty}
  \frac{
    \Gamma
    \left(
      1
      -
      n
      +
      \frac{1}{\hbar}
    \right)
    g_{\overline{\mu_{1}}\nu_{1}}
    \cdots
    g_{\overline{\mu_{n}}\nu_{n}}
  }{
    n!
    \Gamma
    \left(
      1
      +
      \frac{1}{\hbar}
    \right)
  }
  \left(
    D^{\nu_{1}}
    \cdots
    D^{\nu_{n}}
    f
  \right)
  \left(
    D^{\overline{\mu_{1}}}
    \cdots
    D^{\overline{\mu_{n}}}
    g
  \right).
  \label{star_prd_C2_SSU_sol}
\end{align}
\label{star_prd_C2_SSU}
\end{theorem}
\noindent
Here \(\mu_{k}, \nu_{k}=1,2\). 
Note that we use the Einstein summation convention in Theorem \ref{star_prd_C2_SSU}. In \cite{SSU1,SSU2}, it is confirmed that the star product on \(\mathbb{C}P^{N}\) coincides with one obtained by Bordemann et al. \cite{BBEW}. See \cite{SSU1,SSU2} for more detail.

\subsection{Example 
of \(\mathbb{C}^{2}\)}
\label{subsec_NR_C2}
\ We shall construct a deformation quantization with separation of variables for \(\mathbb{C}^{2},\) the simplest concrete example via Theorem \ref{main_thm_N2}, and show that it also reproduces the star product with separation of variables by Karabegov \cite{Kara2}. Since \(\mathbb{C}^{2}\) is flat, i.e. \(R_{ABC}^{\quad \ \ \, D}=0,\) each \(X_{k}\left(\infty\right)\) is a constant multiple of the identity matrix, and \(T_{n}\left(\infty\right)\) can be obtained easily.

\begin{prop}
The formula for determining a star product with separation of variables on \(\mathbb{C}^{2}\) is given by
\begin{align}
  T_{n}\left(\infty\right)
  =
  \left(
  \frac{\hbar}{2}
  \right)^{n}
  \frac{1}{n!}
  \left(
    \begin{array}{ccccc}
      \tbinom{n}{0}  &  &  &  \text{\huge{0}}  \\
        &  \ddots  &  &  \\
        &    &  \tbinom{n}{n}  &  \\
      \text{\huge{0}}  &  &  &  0  &  \\
        &  &  &  &  \ddots
    \end{array}
  \right).
  \label{Tn_C}
\end{align}
\label{thm_Tn_C2}
\end{prop}

\begin{prf*}
We now take the canonical coordinates 
\(
  z^{1}=x^{1}+iy^{1},
  z^{2}=x^{2}+iy^{2}
\) 
as the local coordinates. Since 
\(
  \left(
    g_{AB}
  \right)
\), 
the component matrix of the K\"{a}hler metric on \(\mathbb{C}^{2}\), is
\begin{align*}
  \left(
    g_{AB}
  \right)
  =
  \left(
    \begin{array}{cccc}
      g_{11}  &  g_{12}
      &  g_{1\overline{1}}  &  g_{1\overline{2}}  \\
      g_{21}  &  g_{22}
      &  g_{2\overline{1}}  &  g_{2\overline{2}}  \\
      g_{\overline{1}1}  &  g_{\overline{1}2}
      &  g_{\overline{1}\,\overline{1}}  &  g_{\overline{1}\,\overline{2}}  \\
      g_{\overline{2}1}  &  g_{\overline{2}2}
      &  g_{\overline{2}\,\overline{1}}  &  g_{\overline{2}\,\overline{2}}
    \end{array}
  \right)
  :=
  \left(
    \begin{array}{cccc}
      g_{z^{1}z^{1}}  &  g_{z^{1}z^{2}}
      &  g_{z^{1}\overline{z}^{1}}  &  g_{z^{1}\overline{z}^{2}}  \\
      g_{z^{2}z^{1}}  &  g_{z^{2}z^{2}}
      &  g_{z^{2}\overline{z}^{1}}  &  g_{z^{2}\overline{z}^{2}}  \\
      g_{\overline{z}^{1}z^{1}}  &  g_{\overline{z}^{1}z^{2}}
      &  g_{\overline{z}^{1}\overline{z}^{1}}  &  g_{\overline{z}^{1}\overline{z}^{2}}  \\
      g_{\overline{z}^{2}z^{1}}  &  g_{\overline{z}^{2}z^{2}}
      &  g_{\overline{z}^{2}\overline{z}^{1}}  &  g_{\overline{z}^{2}\overline{z}^{2}}  \\
    \end{array}
  \right)
  =
  \left(
    \begin{array}{cccc}
      0  &  0  &  \frac{1}{2}  &  0  \\
      0  &  0  &  0  &  \frac{1}{2}  \\
      \frac{1}{2}  &  0  &  0  &  0  \\
      0  &  \frac{1}{2}  &  0  &  0
    \end{array}
  \right),
\end{align*}

\noindent
so the inverse matrix 
\(
  \left(
    g^{AB}
  \right)
\) 
is 
\begin{align*}
  \left(
    g^{AB}
  \right)
  =
  \left(
    \begin{array}{cccc}
      g^{11}  &  g^{12}
      &  g^{1\overline{1}}  &  g^{1\overline{2}}  \\
      g^{21}  &  g^{22}
      &  g^{2\overline{1}}  &  g^{2\overline{2}}  \\
      g^{\overline{1}1}  &  g^{\overline{1}2}
      &  g^{\overline{1}\,\overline{1}}  &  g^{\overline{1}\,\overline{2}}  \\
      g^{\overline{2}1}  &  g^{\overline{2}2}
      &  g^{\overline{2}\,\overline{1}}  &  g^{\overline{2}\,\overline{2}}
    \end{array}
  \right)
  :=
  \left(
    \begin{array}{cccc}
      g^{z^{1}z^{1}}  &  g^{z^{1}z^{2}}
      &  g^{z^{1}\overline{z}^{1}}  &  g^{z^{1}\overline{z}^{2}}  \\
      g^{z^{2}z^{1}}  &  g^{z^{2}z^{2}}
      &  g^{z^{2}\overline{z}^{1}}  &  g^{z^{2}\overline{z}^{2}}  \\
      g^{\overline{z}^{1}z^{1}}  &  g^{\overline{z}^{1}z^{2}}
      &  g^{\overline{z}^{1}\overline{z}^{1}}  &  g^{\overline{z}^{1}\overline{z}^{2}}  \\
      g^{\overline{z}^{2}z^{1}}  &  g^{\overline{z}^{2}z^{2}}
      &  g^{\overline{z}^{2}\overline{z}^{1}}  &  g^{\overline{z}^{2}\overline{z}^{2}}  \\
    \end{array}
  \right)
  =
  \left(
    \begin{array}{cccc}
      0  &  0  &  2  &  0  \\
      0  &  0  &  0  &  2  \\
      2  &  0  &  0  &  0  \\
      0  &  2  &  0  &  0
    \end{array}
  \right).
\end{align*}
Since 
\(
  g_{AB},
  \
  g^{AB}
\) 
are constant matrices, 
\(
  \Gamma_{BC}^{A}=0
  \
  \left(
    A,B,C
    =
    1,\overline{1},2,\overline{2}
  \right),
\) 
therefore \(R_{ABC}^{\quad \ \ \, D}=0.\) Thus each 
\(
  X_{k}
  \in
  M_{k+1}
  \left(
    \mathbb{C}
    \left[\!
      \left[
        \hbar
      \right]\!
    \right]
  \right)
  \
  \left(
    k=1,\cdots,n
  \right)
\) 
is obtained as 
\(
  X_{k}
  =
  k\mathrm{Id}_{k+1},
\) 
and especially
\(
  X_{k}^{-1}
  =
  \frac{1}{k}
  \mathrm{Id}_{k+1}.
\) 
Next, we must determine 
\(
  F_{k}^{i_{k}}
\) 
and 
\(
  B_{k}^{i_{k}}.
\) 
Since the choice of the functions 
\(
  \Theta_{\overline{\mu}}^{i_{k}}
\) 
and 
\(
  \Theta_{\nu}^{i_{k}}
\) 
is not unique, there are many possible 
\(
  F_{k}^{i_{k}}
\) 
and 
\(
  B_{k}^{i_{k}}.
\) 
 So we shall now simply choose 
\(
  \Theta_{\mu}^{i_{k}}
  =
  \Theta_{\overline{\mu}}^{i_{k}}
  =
  \frac{1}{\sqrt{2}}
  \delta^{i_{k}}_{\ \,\mu}
\) 
here. 
We can confirm that these 
\(
  \Theta_{\overline{\mu}}^{i_{k}}
\) 
and 
\(
  \Theta_{\nu}^{i_{k}}
\) 
recover
\(
  g_{\overline{\mu}\nu}
  =
  \sum_{i_{k}=1}^{4}
  \Theta_{\overline{\mu}}^{i_{k}}
  \Theta_{\nu}^{i_{k}}.
\) 
So in this case
\begin{align*}
  F_{k}^{i_{k}}
  =
  \frac{1}{\sqrt{2}}
  \left(
    \delta^{i_{k}}_{\ \,1}
    \mathrm{Id}_{k+1}
    +
    \delta^{i_{k}}_{\ \,2}
    F_{d,k+1}
  \right),
  \
  B_{k}^{i_{k}}
  =
  \frac{1}{\sqrt{2}}
  \left(
    \delta^{i_{k}}_{\ \,1}
    \mathrm{Id}_{k+1}
    +
    \delta^{i_{k}}_{\ \,2}
    F_{d,k+1}
  \right).
\end{align*}

\noindent
By 
\(
  X^{-1}_{k}
  =
  \frac{1}{k}
  \mathrm{Id}_{k}
\), 
%
we obtain
\begin{align*}
  T_{n}\left(\infty\right)
  =&
  \hbar^{n}
  \sum_{i_{1},\cdots,i_{n}=1}^{4}
  F^{i_{n}}_{n}\left(\infty\right)
  \cdots
  F^{i_{1}}_{1}\left(\infty\right)
  T_{0}\left(\infty\right)
  \left(
    B^{i_{1}}_{1}\left(\infty\right)
    X^{-1}_{1}\left(\infty\right)
  \right)
  \cdots
  \left(
    B^{i_{n}}_{n}\left(\infty\right)
    X^{-1}_{n}\left(\infty\right)
  \right).
  \\
  =&
  \frac{\hbar^{n}}{n!}
  \left(
    \frac{1}{\sqrt{2}}
  \right)^{2n}
  \sum_{i_{1},\cdots,i_{n}=1}^{4}
  \left(
    \delta^{i_{n}}_{\ \,1}
    \mathrm{Id}\left(\infty\right)
    +
    \delta^{i_{n}}_{\ \,2}
    F_{d}\left(\infty\right)
  \right)
  \cdots
  \left(
    \delta^{i_{1}}_{\ \,1}
    \mathrm{Id}\left(\infty\right)
    +
    \delta^{i_{1}}_{\ \,2}
    F_{d}\left(\infty\right)
  \right)
  T_{0}\left(\infty\right)
  \\
  &
  \cdot
  \left(
    \delta^{i_{1}}_{\ \,1}
    \mathrm{Id}\left(\infty\right)
    +
    \delta^{i_{1}}_{\ \,2}
    F_{r}\left(\infty\right)
  \right)
  \cdots
  \left(
    \delta^{i_{n}}_{\ \,1}
    \mathrm{Id}\left(\infty\right)
    +
    \delta^{i_{n}}_{\ \,2}
    F_{r}\left(\infty\right)
  \right).
\end{align*}

\noindent
By using
\begin{align*}
  \left(
    \begin{array}{c|c}
      A_{k+1} &  0  \\  \hline
      0  &  0
    \end{array}
  \right)
  \left(
    \begin{array}{c|c}
      X^{-1}_{k} &  0  \\  \hline
      0  &  0
    \end{array}
  \right)
  =
  \left(
    \begin{array}{c|c}
      A_{k+1} &  0  \\  \hline
      0  &  0
    \end{array}
  \right)
  \left(
    \begin{array}{c|c}
      X^{-1}_{k} &  0  \\  \hline
      0  &  B
    \end{array}
  \right)
\end{align*}

\noindent
for any 
\(
  A_{k+1}
  \in
  M_{k+1}
  \left(
    \mathbb{C}
    \left[\!
      \left[
        \hbar
      \right]\!
    \right]
  \right)
\)
and 
\(
  B
  \in
  M_{\infty}
  \left(
    \mathbb{C}
    \left[\!
      \left[
        \hbar
      \right]\!
    \right]
  \right)
\). 
The above equation can be written as
\begin{flalign*}
  T_{n}\left(\infty\right)
  =&
  \left(
  \frac{\hbar}{2}
  \right)^{n}
  \frac{1}{n!}
  \sum_{\mu_{1},\nu_{1},\cdots,\mu_{n},\nu_{n}=1}^{2}
  \delta_{\mu_{1}\nu_{1}}
  \cdots
  \delta_{\mu_{n}\nu_{n}}
  \left(
    F_{d}\left(\infty\right)
  \right)^{
    \gamma
    \left(
      \nu^{\left(n\right)}
    \right)
  }
  T_{0}\left(\infty\right)
  \left(
    F_{r}\left(\infty\right)
  \right)^{
    \gamma
    \left(
      \mu^{\left(n\right)}
    \right)
  }
  \\
  =&
  \left(
  \frac{\hbar}{2}
  \right)^{n}
  \frac{1}{n!}
  \sum_{\mu_{1},\cdots,\mu_{n}=1}^{2}
  \left(
    F_{d}\left(\infty\right)
  \right)^{
    \gamma
    \left(
      \mu^{\left(n\right)}
    \right)
  }
  T_{0}\left(\infty\right)
  \left(
    F_{r}\left(\infty\right)
  \right)^{
    \gamma
    \left(
      \mu^{\left(n\right)}
    \right)
  },
\end{flalign*}

\noindent
where 
\(
  \gamma\left(\cdot\right)
  :
  \left\{
    1,2
  \right\}^{n}
  \to
  \mathbb{Z}_{\geq0}
\) 
is defined by 
\begin{align*}
  \gamma
  \left(
    \mu^{\left(n\right)}
  \right)
  :=
  \sum_{k=1}^{n}
  \left[
    \frac{\mu_{k}}{2}
  \right]
  =
  \sum_{k=1}^{n}
  \delta_{\mu_{k}2}
\end{align*}

\noindent
for multi-indices 
\(
  \mu^{\left(n\right)}
  =
  \left(
    \mu_{1},
    \cdots,
    \mu_{n}
  \right),
\) 
and 
\(
  \left[
    \cdot
  \right]
  :
  \mathbb{R}
  \to
  \mathbb{Z}
\) 
is Gauss symbol. We define 
\(
  \gamma
  \left(
    \nu^{\left(n\right)}
  \right)
\) 
in the same way. In other words, 
\(
  \gamma\left(\cdot\right)
  :
  \left\{
    1,2
  \right\}^{n}
  \to
  \mathbb{Z}_{\geq0}
\) 
gives the number of each 
\(
  \mu_{k}
\) 
(or 
\(
  \nu_{k}
\)
) 
that is 2. We can rewrite
\begin{align*}
  \left(
  \frac{\hbar}{2}
  \right)^{n}
  \frac{1}{n!}
  \sum_{\mu_{1},\cdots,\mu_{n}=1}^{2}
  \left(
    F_{d}\left(\infty\right)
  \right)^{
    \gamma
    \left(
      \mu^{\left(n\right)}
    \right)
  }
  T_{0}\left(\infty\right)
  \left(
    F_{r}\left(\infty\right)
  \right)^{
    \gamma
    \left(
      \mu^{\left(n\right)}
    \right)
  }
  =&
  \left(
  \frac{\hbar}{2}
  \right)^{n}
  \frac{1}{n!}
  \sum_{j=0}^{n}
  \binom{n}{j}
  \left(
    F_{d}\left(\infty\right)
  \right)^{j}
  T_{0}\left(\infty\right)
  \left(
    F_{r}\left(\infty\right)
  \right)^{j}
\end{align*}

\noindent
since there are \(\tbinom{n}{j}\) combinations of each \(\mu_{k}\) such that 
\(
  \gamma
  \left(
    \mu^{\left(n\right)}
  \right)
  =
  j.
\) 
After summing over 
\(
  j
  \in
  \left\{1,\cdots,n\right\},
\) 
this is concretely expressed as 
\begin{align*}
  \left(
  \frac{\hbar}{2}
  \right)^{n}
  \frac{1}{n!}
  \sum_{j=0}^{n}
  \binom{n}{j}
  \left(
    F_{d}\left(\infty\right)
  \right)^{j}
  T_{0}\left(\infty\right)
  \left(
    F_{r}\left(\infty\right)
  \right)^{j}
  =
  \left(
  \frac{\hbar}{2}
  \right)^{n}
  \frac{1}{n!}
  \left(
    \begin{array}{ccccc}
      \binom{n}{0}  &  &  &  \text{\huge{0}}  \\
        &  \ddots  &  &  \\
        &    &  \binom{n}{n}  &  \\
      \text{\huge{0}}  &  &  &  0  &  \\
        &  &  &  &  \ddots
    \end{array}
  \right)
\end{align*}

\noindent
by matrix representations. This proof was completed.\(_\blacksquare\)
\end{prf*}

\noindent
By Theorem \ref{thm_Tn_C2}, the coefficients
\(
  T_{\overrightarrow{\alpha_{n}},\overrightarrow{\beta_{n}^{\ast}}}^{n}
\) 
(or 
\(
  T_{ij}^{n}
\)
) 
is determined as
\begin{align}
  T_{\overrightarrow{\alpha_{n}},\overrightarrow{\beta_{n}^{\ast}}}^{n}
  =
  T_{ij}^{n}
  =
  \begin{cases}
    \displaystyle
    \frac{
      \binom{n}{i-1}
      \hbar^{n}
    }{
      2^{n}
      \cdot
      n!
    }
    &
    \
    \left(
      i
      =
      j
    \right)
    \\
    0
    &
    \
    \left(
      i
      \neq
      j
    \right)
  \end{cases},
  \label{coeff_Tn_C}
\end{align}

\noindent
where 
\(
  \alpha_{2}^{n}
  =
  i-1
\) 
and 
\(
  \beta_{2}^{n}
  =
  j-1.
\) 
Hence, the star product with separation of variables on \(\mathbb{C}^{2}\) is obtained as follows :
\begin{align}
  f
  \ast
  g
  =&
  \sum_{n=0}^{\infty}
  \sum_{\overrightarrow{\alpha_{n}},\overrightarrow{\beta_{n}^{\ast}}}
  T_{\overrightarrow{\alpha_{n}},\overrightarrow{\beta_{n}^{\ast}}}^{n}
  \left(
    D^{\overrightarrow{\alpha_{n}}}f
  \right)
  \left(
    D^{\overrightarrow{\beta_{n}^{\ast}}}g
  \right)
  \nonumber
  \\
  =&
  \sum_{n=0}^{\infty}
  \sum_{i=1}^{n+1}
  T_{ii}^{n}
  \left\{
    \left(
      D^{1}
    \right)^{
      n
      -
      \left(
        i-1
      \right)
    }
    \left(
      D^{2}
    \right)^{i-1}
    f
  \right\}
  \left\{
    \left(
      D^{\overline{1}}
    \right)^{
      n
      -
      \left(
        i-1
      \right)
    }
    \left(
      D^{\overline{2}}
    \right)^{i-1}
    g
  \right\}.
  \label{star_C_1}
\end{align}

\noindent
Substituting \eqref{coeff_Tn_C} into the right-hand side of \eqref{star_C_1},

\begin{flalign}
  &
  \sum_{n=0}^{\infty}
  \sum_{i=1}^{n+1}
  T_{ii}^{n}
  \left\{
    \left(
      D^{1}
    \right)^{
      n
      -
      \left(
        i-1
      \right)
    }
    \left(
      D^{2}
    \right)^{i-1}
    f
  \right\}
  \left\{
    \left(
      D^{\overline{1}}
    \right)^{
      n
      -
      \left(
        i-1
      \right)
    }
    \left(
      D^{\overline{2}}
    \right)^{i-1}
    g
  \right\}
  &
  \nonumber
  \\
  &=
  \sum_{n=0}^{\infty}
  \sum_{i=1}^{n+1}
  \left(
    \frac{\hbar}{2}
  \right)^{n}
  \frac{1}{n!}
  \binom{n}{i-1}
  \left\{
    \left(
      g^{1\overline{1}}
      \partial_{\overline{1}}
    \right)^{
      n
      -
      \left(
        i-1
      \right)
    }
    \left(
      g^{2\overline{2}}
      \partial_{\overline{2}}
    \right)^{
      i-1
    }
    f
  \right\}
  \left\{
    \left(
      g^{\overline{1}1}
      \partial_{1}
    \right)^{
      n
      -
      \left(
        i-1
      \right)
    }
    \left(
      g^{\overline{2}2}
      \partial_{2}
    \right)^{
      i-1
    }
    g
  \right\}
  &
  \nonumber
  \\
  &=
  \sum_{n=0}^{\infty}
  \frac{
    \left(
      2\hbar
    \right)^{n}
  }{n!}
  \sum_{j=0}^{n}
  \binom{n}{j}
  \left(
    \partial_{\overline{z}^{1}}^{
      n
      -
      j
    }
    \partial_{\overline{z}^{2}}^{
      j
    }
    f
  \right)
  \left(
    \partial_{z^{1}}^{
      n
      -
      j
    }
    \partial_{z^{2}}^{
      j
    }
    g
  \right),
  &
  \label{star_C_2}
\end{flalign}

\noindent
where we rewrite \(i-1\) as \(j\) in the last line in \eqref{star_C_2}. Here, using the left (right) differential operator

\begin{align*}
  f
  \overleftarrow{
    \partial_{z^{A}}
  }
  :=
  \partial_{z^{A}}
  f,
  \
  \overrightarrow{
    \partial_{z^{A}}
  }
  g
  :=
  \partial_{z^{A}}
  g
  \
  \left(
    A
    \in
    \left\{
      1,2,\bar{1},\bar{2}
    \right\}
  \right),
\end{align*}

\noindent
we obtain the star product on \(\mathbb{C}^{2}\) as
\begin{align}
  f
  \ast
  g
  &
  =
  \sum_{n=0}^{\infty}
  \frac{
    \left(
      2\hbar
    \right)^{n}
  }{n!}
  \sum_{j=0}^{n}
  \binom{n}{j}
  \left(
    \partial_{\overline{z}^{1}}^{
      n
      -
      j
    }
    \partial_{\overline{z}^{2}}^{
      j
    }
    f
  \right)
  \left(
    \partial_{z^{1}}^{
      n
      -
      j
    }
    \partial_{z^{2}}^{
      j
    }
    g
  \right)
  \nonumber
  \\
  &=
  \sum_{n=0}^{\infty}
  \frac{
    \left(
      2\hbar
    \right)^{n}
  }{n!}
  f
  \left\{
    \sum_{j=0}^{n}
    \binom{n}{j}
    \left(
      \overleftarrow{
        \partial_{\overline{z}^{1}}
      }
      \overrightarrow{
        \partial_{z^{1}}
      }
    \right)^{
      n-j
    }
    \left(
      \overleftarrow{
        \partial_{\overline{z}^{2}}
      }
      \overrightarrow{
        \partial_{z^{2}}
      }
    \right)^{j}
  \right\}
  g.
  \label{star_C_3}
\end{align}
\noindent
Let us compare the previous result \eqref{star_prd_C2_Kara_Voro_sol} to verify that Theorem \ref{main_thm_N2} works well. By using the binomial theorem for the left and right differential operators in \eqref{star_C_3}, we have
\begin{align}
  \sum_{n=0}^{\infty}
  \frac{
    \left(
      2\hbar
    \right)^{n}
  }{n!}
  f
  \left\{
    \sum_{j=0}^{n}
    \binom{n}{j}
    \left(
      \overleftarrow{
        \partial_{\overline{z}^{1}}
      }
      \overrightarrow{
        \partial_{z^{1}}
      }
    \right)^{
      n-j
    }
    \left(
      \overleftarrow{
        \partial_{\overline{z}^{2}}
      }
      \overrightarrow{
        \partial_{z^{2}}
      }
    \right)^{j}
  \right\}
  g
  =&
  \sum_{n=0}^{\infty}
  \frac{
    \left(
      2\hbar
    \right)^{n}
  }{n!}
  f
  \left(
    \overleftarrow{
      \partial_{\overline{z}^{1}}
    }
    \overrightarrow{
      \partial_{z^{1}}
    }
    +
    \overleftarrow{
      \partial_{\overline{z}^{2}}
    }
    \overrightarrow{
      \partial_{z^{2}}
    }
  \right)^{n}
  g
  \nonumber
  \\
  =&
  f
  \mathrm{exp}
  2\hbar
  \left(
    \overleftarrow{
      \partial_{\overline{z}^{1}}
    }
    \overrightarrow{
      \partial_{z^{1}}
    }
    +
    \overleftarrow{
      \partial_{\overline{z}^{2}}
    }
    \overrightarrow{
      \partial_{z^{2}}
    }
  \right)
  g.
  \label{star_C_4}
\end{align}
\noindent
Here, if we put the formal parameter \(\nu\) as \(\nu:=2\hbar\), we obtain
\(
  f
  \ast
  g
  =
  f
  \mathrm{exp}
  \nu
  \left(
    \overleftarrow{
      \partial_{\overline{z}^{1}}
    }
    \overrightarrow{
      \partial_{z^{1}}
    }
    +
    \overleftarrow{
      \partial_{\overline{z}^{2}}
    }
    \overrightarrow{
      \partial_{z^{2}}
    }
  \right)
  g,
\) 
which is the star product with separation of variables by Karabegov 
\eqref{star_prd_C2_Kara_Voro_sol}.
\bigskip
\\
\ We stated in Subsection \ref{NR_plus} that \(T_{n}\left(\infty\right)\) given in \eqref{main_them_N2_Tinf} solved via Theorem \ref{main_thm_N2} satisfies the formula \eqref{minus_thm_N2_eq_rsz} in Theorem \ref{minus_thm_N2}. So we shall directly verify that \(T_{n}\left(\infty\right)\) obtained via Theorem \ref{main_thm_N2} actually satisfies the formula \eqref{minus_thm_N2_eq_rsz} in Theorem \ref{minus_thm_N2} in 
the case of \(\mathbb{C}^{2}\). 
By Theorem \ref{minus_thm_N2}, \(Y_{n}\) is obtained as
\begin{align*}
  Y_{n}
  =
  \left(
    \begin{array}{ccccc}
      n  &  &  &  &  \text{\huge{0}} \\
        &  n-2  &  &  &  \\
        &  &  \ddots  &  &  \\
        &    &  &  -(n-2)  &  \\
      \text{\huge{0}}  &  &  &  &  -n
    \end{array}
  \right).
\end{align*}

\noindent
The embedding \(Y_{n}\left(\infty\right)\) for \(Y_{n}\) is as in the way of \eqref{emb_gene}. 
%
Thus,
\begin{align*}
  T_{n}\left(\infty\right)
  Y_{n}\left(\infty\right)
  =&
  \left(
  \frac{\hbar}{2}
  \right)^{n}
  \frac{1}{n!}
  \left(
    \begin{array}{ccccccc}
      n \binom{n}{0}
      &  &  &  &  \text{\huge{0}}  &  &  \\
        &  (n-2) \binom{n}{1}  &  &  &  &  &  \\
        &  &  \ddots  &  &  &  &  \\
        &    &  &  -(n-2) \binom{n}{n-1}  &  &  &  \\
      \text{\huge{0}}  &  &  &  &  -n \binom{n}{n}  &  &  \\
        &  &  &  &  &  0  &  \\
        &  &  &  &  &  &  \ddots  
    \end{array}
  \right).
\end{align*}

\noindent
On the other hand, 
\(
  \hbar
  C_{n}^{'}
  \left(\infty\right)
\) 
is calculated as

\begin{align*}
  \hbar
  C_{n}^{'}\left(\infty\right)
  =&
  \hbar
  \left(
    \frac{1}{2}
    T_{n-1}\left(\infty\right)
    -
    \frac{1}{2}
    F_{d}\left(\infty\right)
    T_{n-1}\left(\infty\right)
    F_{r}\left(\infty\right)
  \right).
\end{align*}

\noindent
We substitute \eqref{Tn_C} for \(n-1\), \eqref{emb_Fd}, and \eqref{emb_Fr} into the the right-hand side of the above. Then,

\begin{align*}
  &
  \left(
  \frac{\hbar}{2}
  \right)^{n}
  \frac{1}{\left(n-1\right)!}
  \left\{
    \left(
      \begin{array}{cccccc}
        \binom{n-1}{0}  &  &  &  \text{\huge{0}}  &  &  \\
          &  \ddots  &  &  &  &  \\
          &    &  \binom{n-1}{n-1}  &  &  &  \\
        \text{\huge{0}}  &  &  &  0  &  &  \\
          &  &  &  &  0  &  \\
          &  &  &  &  &  \ddots  
      \end{array}
    \right)
    -
    \left(
      \begin{array}{cccccc}
        0  &  &  &  \text{\huge{0}}  &  &  \\
          &  \binom{n-1}{0}  &  &  &  &  \\
          &    &  \ddots  &  &  &  \\
        \text{\huge{0}}  &  &  &  \binom{n-1}{n-1}  &  &  \\
          &  &  &  &  0  &  \\
          &  &  &  &  &  \ddots  
      \end{array}
    \right)
  \right\}
  \\
  &=
  \left(
  \frac{\hbar}{2}
  \right)^{n}
  \frac{1}{n!}
  \left(
    \begin{array}{ccccccc}
      n \binom{n}{0}
      &  &  &  &  \text{\huge{0}}  &  &  \\
        &  (n-2) \binom{n}{1}  &  &  &  &  &  \\
        &  &  \ddots  &  &  &  &  \\
        &    &  &  -(n-2) \binom{n}{n-1}  &  &  &  \\
      \text{\huge{0}}  &  &  &  &  -n \binom{n}{n}  &  &  \\
        &  &  &  &  &  0  &  \\
        &  &  &  &  &  &  \ddots  
    \end{array}
  \right).
\end{align*}

\noindent
Hence, for any 
\(
  n
  \in
  \mathbb{Z}_{\geq0},
\) 
the matrix 
\(
  T_{n}\left(\infty\right)
\) for \(\mathbb{C}^{2}\) satisfies \eqref{minus_thm_N2_eq_another_rsz}.

\subsection{
Example of \(\mathbb{C}P^{2}\)
}
\label{subsec_NR_CP2}
\ 
In this subsection, we confirm that the deformation quantization with separation of variables for \(\mathbb{C}P^{2}\) obtained via our main result (Theorem \ref{main_thm_N2}) satisfies the identities in Theorem \ref{minus_thm_N2} or equivalently \eqref{minus_thm_N2_eq_rsz} in Subsection \ref{NR_minus}.
\begin{prop}
The formula for determining a star product with separation of variables on \(\mathbb{C}P^{2}\) is given by
\begin{align}
  T_{n}\left(\infty\right)
  =
  \frac{
    \Gamma
    \left(
      1
      -
      n
      +
      \frac{1}{\hbar}
    \right)
  }{
    n!
    \Gamma
    \left(
      1
      +
      \frac{1}{\hbar}
    \right)
  }
  \sum_{\mu_{1},\nu_{1},\cdots,\mu_{n},\nu_{n}=1}^{2}
  g_{\overline{\mu_{1}}\nu_{1}}
  \cdots
  g_{\overline{\mu_{n}}\nu_{n}}
  \left(
    F_{d}\left(\infty\right)
  \right)^{
    \gamma
    \left(
      \nu^{\left(n\right)}
    \right)
  }
  T_{0}\left(\infty\right)
  \left(
    F_{r}\left(\infty\right)
  \right)^{
    \gamma
    \left(
      \mu^{\left(n\right)}
    \right)
  },
  \label{formula_CP2}
\end{align}
where each \(g_{\overline{\mu_{k}}\nu_{k}}\) is Fubini-Study metric and 
\(
  \gamma\left(\cdot\right)
  :
  \left\{
    1,2
  \right\}^{n}
  \to
  \mathbb{Z}_{\geq0}
\) 
is defined by
\begin{align}
  \gamma
  \left(
    \mu^{\left(n\right)}
  \right)
  :=
  \sum_{k=1}^{2}
  \left[
    \frac{\mu_{k}}{2}
  \right]
  =
  \sum_{k=1}^{2}
  \delta_{\mu_{k}2},
  \qquad
  \gamma
  \left(
    \nu^{\left(n\right)}
  \right)
  :=
  \sum_{k=1}^{2}
  \left[
    \frac{\nu_{k}}{2}
  \right]
  =
  \sum_{k=1}^{2}
  \delta_{\nu_{k}2}
  \label{gamma_gauss}
\end{align}

\noindent
for multi-indices
\(
  \mu^{\left(n\right)}
  =
  \left(
    \mu_{1},
    \cdots,
    \mu_{n}
  \right),
  \nu^{\left(n\right)}
  =
  \left(
    \nu_{1},
    \cdots,
    \nu_{n}
  \right).
\)
\label{thm_Tn_CP2}
\end{prop}

\noindent
We can prove this proposition by direct calculations of \eqref{main_them_N2_Tinf} in Theorem \ref{main_thm_N2}. See \cite{OS_Varna} for 
detailed proof of it. By Proposition \ref{thm_Tn_CP2}, the star product with separation of variables on \(\mathbb{C}P^{2}\) is obtained by
\begin{align}
  f
  \ast
  g
  =
  \sum_{n=0}^{\infty}
  \frac{
    \Gamma
    \left(
      1
      -
      n
      +
      \frac{1}{\hbar}
    \right)
    g_{\overline{\mu_{1}}\nu_{1}}
    \cdots
    g_{\overline{\mu_{n}}\nu_{n}}
  }{
    n!
    \Gamma
    \left(
      1
      +
      \frac{1}{\hbar}
    \right)
  }
  \left(
    D^{\nu_{1}}
    \cdots
    D^{\nu_{n}}
    f
  \right)
  \left(
    D^{\overline{\mu_{1}}}
    \cdots
    D^{\overline{\mu_{n}}}
    g
  \right).
  \label{star_CP2}
\end{align}

\noindent
Detailed calculations for deriving equation \eqref{star_CP2} from \eqref{formula_CP2} are also given in \cite{OS_Varna}. This star product coincides with 
\eqref{star_prd_C2_SSU_sol} by Sako-Suzuki-Umetsu 
. This result implies that Theorem \ref{main_thm_N2} works well.
\bigskip
\\
\ As in the case of \(\mathbb{C}^{2},\) we shall see that \eqref{formula_CP2} satisfies \eqref{minus_thm_N2_eq_rsz} in Theorem \ref{minus_thm_N2}. 
By Theorem \ref{minus_thm_N2}, the components of \(Y_{n}\) for \(\mathbb{C}P^{2}\) are calculated as
\begin{align*}
  Y_{i,j}^{n}
  =
  \begin{cases}
    \hbar
    \left\{
      n
      -
      2
      \left(
        i-1
      \right)
    \right\}
    \left(
      1
      -
      n
      +
      \frac{1}{\hbar}
    \right)
    &
    \left(
      i=j
    \right)
    \\
    0
    &
    \left(
      i
      \neq
      j
    \right)
  \end{cases}.
\end{align*}

%
%
%

\noindent
\(Y_{n}\left(\infty\right)\) is written as
\begin{align}
  Y_{n}\left(\infty\right)
  =
  \hbar
  \left(
    1
    -
    n
    +
    \frac{1}{\hbar}
  \right)
  \left(
    \begin{array}{ccccccc}
      n  &  &  &  &  \text{\huge{0}}  &  &  \\
        &  n-2  &  &  &  &  &  \\
        &  &  \ddots  &  &  &  &  \\
        &    &  &  -(n-2)  &  &  &  \\
      \text{\huge{0}}  &  &  &  &  -n  &  &  \\
        &  &  &  &  &  0  &  \\
        &  &  &  &  &  &  \ddots  
    \end{array}
  \right).
  \label{Yn_CP2}
\end{align}

\noindent

\noindent
To calculate 
\(
  T_{n}\left(\infty\right)
  Y_{n}\left(\infty\right)
\) 
we use the following proposition.

\begin{prop}
\begin{align*}
  &
  \left(
    F_{d}\left(\infty\right)
  \right)^{
  k
  }
  T_{0}\left(\infty\right)
  \left(
    F_{r}\left(\infty\right)
  \right)^{
    l
  }
  \left(
    \begin{array}{ccccccc}
      n  &  &  &  &  \text{\huge{0}}  &  &  \\
        &  n-2  &  &  &  &  &  \\
        &  &  \ddots  &  &  &  &  \\
        &    &  &  -(n-2)  &  &  &  \\
      \text{\huge{0}}  &  &  &  &  -n  &  &  \\
        &  &  &  &  &  0  &  \\
        &  &  &  &  &  &  \ddots  
    \end{array}
  \right)
  \\
  &=
  \left(
    n
    -
    2
    l
  \right)
  \left(
    F_{d}\left(\infty\right)
  \right)^{
    k
  }
  T_{0}\left(\infty\right)
  \left(
    F_{r}\left(\infty\right)
  \right)^{
    l
  },
  \qquad
  \left(
    k,
    l
    =
    0,
    \cdots,
    n
  \right).
\end{align*}
\label{Prop_Fd_Fr}
\end{prop}

\begin{prf*}
Since the properties of \(F_{d}\left(\infty\right)\) and \(F_{r}\left(\infty\right)\) in Remark \ref{rem_Fd_Fr},
\begin{align}
  \left(
    F_{d}\left(\infty\right)
  \right)^{k}
  T_{0}\left(\infty\right)
  \left(
    F_{r}\left(\infty\right)
  \right)^{l}
  =
  \left(
    \begin{array}{c|c}
      \left(
        \delta_{i,k}
        \delta_{j,l}
      \right)
      &
      0
      \\ \hline
      0
      &
      0
    \end{array}
  \right),
  \label{calc_Fd_T_Fr}
\end{align}

\noindent
where 
\(
  \left(
    \delta_{i,k}
    \delta_{j,l}
  \right)
  \in
  M_{n+1}
  \left(
  \mathbb{C}
  \left[\!
    \left[
      \hbar
    \right]\!
  \right]
  \right)
\) 
is the square matrix of order \(n+1\) such that the 
\(\left(k,l\right)\) component
is 1 and the others are 0. And note that the indices \(k\) and \(l\) are fixed. Then we get the equation. \(_{\blacksquare}\)
\end{prf*}

\noindent
By Proposition \ref{Prop_Fd_Fr}, 
\(
  T_{n}\left(\infty\right)
  Y_{n}\left(\infty\right)
\) 
is calculated as
\begin{flalign}
  &
  T_{n}\left(\infty\right)
  Y_{n}\left(\infty\right)
  &
  \nonumber
  \\
  &=
  \frac{
    \hbar
    \left(
      1
      -
      n
      +
      \frac{1}{\hbar}
    \right)
    \Gamma
    \left(
      1
      -
      n
      +
      \frac{1}{\hbar}
    \right)
  }{
    n!
    \Gamma
    \left(
      1
      +
      \frac{1}{\hbar}
    \right)
  }
  \sum_{\mu_{1},\nu_{1},\cdots,\mu_{n},\nu_{n}=1}^{2}
  g_{\overline{\mu_{1}}\nu_{1}}
  \cdots
  g_{\overline{\mu_{n}}\nu_{n}}
  &
  \nonumber
  \\
  &
  \mbox{}
  \quad
  \cdot
  \left(
    F_{d}\left(\infty\right)
  \right)^{
    \gamma
    \left(
      \nu^{\left(n\right)}
    \right)
  }
  T_{0}\left(\infty\right)
  \left(
    F_{r}\left(\infty\right)
  \right)^{
    \gamma
    \left(
      \mu^{\left(n\right)}
    \right)
  }
  \left(
    \begin{array}{ccccccc}
      n  &  &  &  &  \text{\huge{0}}  &  &  \\
        &  n-2  &  &  &  &  &  \\
        &  &  \ddots  &  &  &  &  \\
        &    &  &  -(n-2)  &  &  &  \\
      \text{\huge{0}}  &  &  &  &  -n  &  &  \\
        &  &  &  &  &  0  &  \\
        &  &  &  &  &  &  \ddots  
    \end{array}
  \right)
  &
  \nonumber
  \\
  &=
  \frac{
    \hbar
    \Gamma
    \left(
      1
      -
      \left(
        n-1
      \right)
      +
      \frac{1}{\hbar}
    \right)
  }{
    \left(
      n-1
    \right)!
    \Gamma
    \left(
      1
      +
      \frac{1}{\hbar}
    \right)
  }
  \!\!
  \sum_{
    \substack{
      \mu_{1},\cdots,\mu_{n}=1
      \\
      \nu_{1},\cdots,\nu_{n}=1
    }
  }^{2}
  \!\!\!\!\!
  g_{\overline{\mu_{1}}\nu_{1}}
  \cdots
  g_{\overline{\mu_{n}}\nu_{n}}
  \left(
    1
    -
    \frac{2}{n}
    \gamma
    \left(
      \mu^{\left(n\right)}
    \right)
  \right)
  \left(
    F_{d}\left(\infty\right)
  \right)^{
    \gamma
    \left(
      \nu^{\left(n\right)}
    \right)
  }
  T_{0}\left(\infty\right)
  \left(
    F_{r}\left(\infty\right)
  \right)^{
    \gamma
    \left(
      \mu^{\left(n\right)}
    \right)
  }.
  \label{another_lhs_CP2_calc}
\end{flalign}

\noindent
%

Next, we calculate the right-hand side of \eqref{minus_thm_N2_eq_rsz}. By Theorem \ref{minus_thm_N2},

\begin{align}
  \hbar
  C_{n}^{'}\left(\infty\right)
  =&
  \hbar
  \left(
    g_{\bar{1}1}
    T_{n-1}\left(\infty\right)
    +
    g_{\bar{1}2}
    F_{d}\left(\infty\right)
    T_{n-1}\left(\infty\right)
    -
    g_{\bar{2}1}
    T_{n-1}\left(\infty\right)
    F_{r}\left(\infty\right)
    -
    g_{\bar{2}2}
    F_{d}\left(\infty\right)
    T_{n-1}\left(\infty\right)
    F_{r}\left(\infty\right)
  \right),
  \nonumber
  \\
  =&
  \frac{
    \hbar
    \Gamma
    \left(
      1
      -
      \left(
        n-1
      \right)
      +
      \frac{1}{\hbar}
    \right)
  }{
    \left(
      n-1
    \right)!
    \Gamma
    \left(
      1
      +
      \frac{1}{\hbar}
    \right)
  }
  \sum_{p_{n}=1}^{4}
  \sum_{\mu_{1},\nu_{1},\cdots,\mu_{n-1},\nu_{n-1}=1}^{2}
  g_{\overline{\mu_{1}}\nu_{1}}
  \cdots
  g_{\overline{\mu_{n-1}}\nu_{n-1}}
  \left(
    \Theta_{1}^{p_{n}}
    \mathrm{Id}\left(\infty\right)
    +
    \Theta_{2}^{p_{n}}
    F_{d}\left(\infty\right)
  \right)
  \nonumber
  \\
  &
  \cdot
  \left(
    F_{d}\left(\infty\right)
  \right)^{
    \gamma
    \left(
      \nu^{\left(n-1\right)}
    \right)
  }
  T_{0}\left(\infty\right)
  \left(
    F_{r}\left(\infty\right)
  \right)^{
    \gamma
    \left(
      \mu^{\left(n-1\right)}
    \right)
  }
  \left(
    \Theta_{\bar{1}}^{p_{n}}
    \mathrm{Id}\left(\infty\right)
    -
    \Theta_{\bar{2}}^{p_{n}}
    F_{r}\left(\infty\right)
  \right)
  .
  \label{Cn_CP2_mid}
\end{align}

\noindent
Here 
\(
  \Theta_{\bar{\mu}}^{i_{k}},
  \Theta_{\nu}^{i_{k}}
  \in
  C^{\infty}
  \left(
    U
  \right)
\) 
are the \(C^\infty\) functions such that 
\(
  g_{\overline{\mu}\nu}
  =
  \sum_{i_{k}=1}^{4}
  \Theta_{\bar{\mu}}^{i_{k}}
  \Theta_{\nu}^{i_{k}}
\). 
We can choose 
\(
  \Theta_{\bar{\mu}}^{i_{k}}
\) 
and 
\(
  \Theta_{\nu}^{i_{k}}
\) 
concretely 
\cite{OS_Varna}, but we do not have to. This is because we can carry out the following calculations without using concrete expressions. The right-hand side of \eqref{Cn_CP2_mid} can be rewritten as
\begin{align}
  &
  \sum_{p_{n}=1}^{4}
  \left(
    \Theta_{1}^{p_{n}}
    \mathrm{Id}\left(\infty\right)
    +
    \Theta_{2}^{p_{n}}
    F_{d}\left(\infty\right)
  \right)
  \left(
    F_{d}\left(\infty\right)
  \right)^{
    \gamma
    \left(
      \nu^{\left(n-1\right)}
    \right)
  }
  T_{0}\left(\infty\right)
  \left(
    F_{r}\left(\infty\right)
  \right)^{
    \gamma
    \left(
      \mu^{\left(n-1\right)}
    \right)
  }
  \left(
    \Theta_{\bar{1}}^{p_{n}}
    \mathrm{Id}\left(\infty\right)
    -
    \Theta_{\bar{2}}^{p_{n}}
    F_{r}\left(\infty\right)
  \right)
  &
  \nonumber
  \\
  &=
  \sum_{\mu_{n},\nu_{n}=1}^{2}
  \left(
    1
    -
    2
    \left[
      \frac{\mu_{n}}{2}
    \right]
  \right)
  g_{\overline{\mu_{n}}\nu_{n}}
  \left(
    F_{d}\left(\infty\right)
  \right)^{
    \left(
      \sum_{k=1}^{n-1}
      \left[
        \frac{\nu_{k}}{2}
      \right]
    \right)
    +
    \left[
      \frac{\nu_{n}}{2}
    \right]
  }
  T_{0}\left(\infty\right)
  \left(
    F_{r}\left(\infty\right)
  \right)^{
    \left(
      \sum_{k=1}^{n-1}
      \left[
        \frac{\mu_{k}}{2}
      \right]
    \right)
    +
    \left[
      \frac{\mu_{n}}{2}
    \right]
  },
  \label{Cn_rhs_mtx}
\end{align}
\noindent
where we use 
\(
  \gamma
  \left(
    \mu^{\left(n-1\right)}
  \right)
  =
  \sum_{k=1}^{n-1}
  \left[
    \frac{\mu_{k}}{2}
  \right]
\), 
\(
  \gamma
  \left(
    \nu^{\left(n-1\right)}
  \right)
  =
  \sum_{k=1}^{n-1}
  \left[
    \frac{\nu_{k}}{2}
  \right]
\) 
and also introduce a new indices \(\mu_{n}\) and \(\nu_{n}\). 
\(
  \left(
    1
    -
    2
    \left[
      \frac{\mu_{n}}{2}
    \right]
  \right)
\) 
in the right-hand side of \eqref{Cn_rhs_mtx} is sign \(\pm1\) using the Gauss symbol and \(-1\) when \(\mu_{n}=2\). Then we obtain
\begin{flalign*}
  \hbar
  C_{n}^{'}\left(\infty\right)
  =&
  \frac{
    \hbar
    \Gamma
    \left(
      1
      -
      \left(
        n-1
      \right)
      +
      \frac{1}{\hbar}
    \right)
  }{
    \left(
      n-1
    \right)!
    \Gamma
    \left(
      1
      +
      \frac{1}{\hbar}
    \right)
  }
  \left(
    1
    -
    2
    \left[
      \frac{\mu_{n}}{2}
    \right]
  \right)
  \sum_{\mu_{1},\nu_{1},\cdots,\mu_{n},\nu_{n}=1}^{2}
  g_{\overline{\mu_{1}}\nu_{1}}
  \cdots
  g_{\overline{\mu_{n}}\nu_{n}}
  \\
  &
  \cdot
  \left(
    F_{d}\left(\infty\right)
  \right)^{
    \sum_{k=1}^{n-1}
    \left(
    \left[
      \frac{\nu_{k}}{2}
    \right]
    \right)
    +
    \left[
      \frac{\nu_{n}}{2}
    \right]
  }
  T_{0}\left(\infty\right)
  \left(
    F_{r}\left(\infty\right)
  \right)^{
    \sum_{k=1}^{n-1}
    \left(
    \left[
      \frac{\mu_{k}}{2}
    \right]
    \right)
    +
    \left[
      \frac{\mu_{n}}{2}
    \right]
  }.
\end{flalign*}
\noindent
Using
\begin{align*}
  \gamma
  \left(
    \mu^{\left(n\right)}
  \right)
  =
  \left(
    \sum_{k=1}^{n-1}
    \left[
      \frac{\mu_{k}}{2}
    \right]
  \right)
  +
  \left[
    \frac{\mu_{n}}{2}
  \right],
  \qquad
  \gamma
  \left(
    \nu^{\left(n\right)}
  \right)
  =
  \left(
    \sum_{k=1}^{n-1}
    \left[
      \frac{\nu_{k}}{2}
    \right]
  \right)
  +
  \left[
    \frac{\nu_{n}}{2}
  \right],
\end{align*}

\noindent
we have
\begin{align}
  &
  \frac{
    \hbar
    \Gamma
    \left(
      1
      -
      \left(
        n-1
      \right)
      +
      \frac{1}{\hbar}
    \right)
  }{
    \left(
      n-1
    \right)!
    \Gamma
    \left(
      1
      +
      \frac{1}{\hbar}
    \right)
  }
  \sum_{\mu_{1},\nu_{1},\cdots,\mu_{n},\nu_{n}=1}^{2}
  \left(
    1
    -
    2
    \left[
      \frac{\mu_{n}}{2}
    \right]
  \right)
  g_{\overline{\mu_{1}}\nu_{1}}
  \cdots
  g_{\overline{\mu_{n-1}}\nu_{n-1}}
  g_{\overline{\mu_{n}}\nu_{n}}
  \nonumber
  \\
  &
  \mbox{}
  \quad
  \cdot
  \left(
    F_{d}\left(\infty\right)
  \right)^{
    \left(
      \sum_{k=1}^{n-1}
      \left[
        \frac{\nu_{k}}{2}
      \right]
    \right)
    +
    \left[
      \frac{\nu_{n}}{2}
    \right]
  }
  T_{0}\left(\infty\right)
  \left(
    F_{r}\left(\infty\right)
  \right)^{
    \left(
      \sum_{k=1}^{n-1}
      \left[
        \frac{\mu_{k}}{2}
      \right]
    \right)
    +
    \left[
      \frac{\mu_{n}}{2}
    \right]
  }
  \nonumber
  \\
  &=
  \frac{
    \hbar
    \Gamma
    \left(
      1
      -
      \left(
        n-1
      \right)
      +
      \frac{1}{\hbar}
    \right)
  }{
    \left(
      n-1
    \right)!
    \Gamma
    \left(
      1
      +
      \frac{1}{\hbar}
    \right)
  }
  \sum_{\mu_{1},\nu_{1},\cdots,\mu_{n},\nu_{n}=1}^{2}
  \left(
    1
    -
    2
    \left[
      \frac{\mu_{n}}{2}
    \right]
  \right)
  g_{\overline{\mu_{1}}\nu_{1}}
  \cdots
  g_{\overline{\mu_{n}}\nu_{n}}
  \left(
    F_{d}\left(\infty\right)
  \right)^{
    \gamma
    \left(
      \nu^{\left(n\right)}
    \right)
  }
  T_{0}\left(\infty\right)
  \left(
    F_{r}\left(\infty\right)
  \right)^{
    \gamma
    \left(
      \mu^{\left(n\right)}
    \right)
  }
  \nonumber
  \\
  &=
  \frac{
    \hbar
    \Gamma
    \left(
      1
      -
      \left(
        n-1
      \right)
      +
      \frac{1}{\hbar}
    \right)
  }{
    \left(
      n-1
    \right)!
    \Gamma
    \left(
      1
      +
      \frac{1}{\hbar}
    \right)
  }
  \sum_{\mu_{1},\nu_{1},\cdots,\mu_{n},\nu_{n}=1}^{2}
  \left(
    1
    -
    \frac{2}{n}
    \cdot
    n
    \left[
      \frac{\mu_{n}}{2}
    \right]
  \right)
  g_{\overline{\mu_{1}}\nu_{1}}
  \cdots
  g_{\overline{\mu_{n}}\nu_{n}}
  \nonumber
  \\
  &
  \mbox{}
  \quad
  \cdot
  \left(
    F_{d}\left(\infty\right)
  \right)^{
    \gamma
    \left(
      \nu^{\left(n\right)}
    \right)
  }
  T_{0}\left(\infty\right)
  \left(
    F_{r}\left(\infty\right)
  \right)^{
    \gamma
    \left(
      \mu^{\left(n\right)}
    \right)
  }.
  \label{another_CP2_before}
\end{align}

\noindent
By replacing the indices, 
the above is expressed as
\begin{align}
  &
  \sum_{\mu_{1},\nu_{1},\cdots,\mu_{n},\nu_{n}=1}^{2}
  \left(
    1
    -
    \frac{2}{n}
    \cdot
    n
    \left[
      \frac{\mu_{n}}{2}
    \right]
  \right)
  g_{\overline{\mu_{1}}\nu_{1}}
  \cdots
  g_{\overline{\mu_{n}}\nu_{n}}
  \left(
    F_{d}\left(\infty\right)
  \right)^{
    \gamma
    \left(
      \nu^{\left(n\right)}
    \right)
  }
  T_{0}\left(\infty\right)
  \left(
    F_{r}\left(\infty\right)
  \right)^{
    \gamma
    \left(
      \mu^{\left(n\right)}
    \right)
  }
  \nonumber
  \\
  &=
  \sum_{\mu_{1},\nu_{1},\cdots,\mu_{n},\nu_{n}=1}^{2}
  \left(
    1
    -
    \underbrace{
      \frac{2}{n}
      \left[
        \frac{\mu_{n}}{2}
      \right]
      -
      \cdots
      -
      \frac{2}{n}
      \left[
        \frac{\mu_{n}}{2}
      \right]
    }_{n}
  \right)
  g_{\overline{\mu_{1}}\nu_{1}}
  \cdots
  g_{\overline{\mu_{n}}\nu_{n}}
  \left(
    F_{d}\left(\infty\right)
  \right)^{
    \gamma
    \left(
      \nu^{\left(n\right)}
    \right)
  }
  T_{0}\left(\infty\right)
  \left(
    F_{r}\left(\infty\right)
  \right)^{
    \gamma
    \left(
      \mu^{\left(n\right)}
    \right)
  }
  \nonumber
  \\
  &=
  \sum_{\mu_{1},\nu_{1},\cdots,\mu_{n},\nu_{n}=1}^{2}
  \left(
    1
    -
    \frac{2}{n}
    \left[
      \frac{\mu_{1}}{2}
    \right]
    -
    \cdots
    -
    \frac{2}{n}
    \left[
      \frac{\mu_{n}}{2}
    \right]
  \right)
  g_{\overline{\mu_{1}}\nu_{1}}
  \cdots
  g_{\overline{\mu_{n}}\nu_{n}}
  \left(
    F_{d}\left(\infty\right)
  \right)^{
    \gamma
    \left(
      \nu^{\left(n\right)}
    \right)
  }
  T_{0}\left(\infty\right)
  \left(
    F_{r}\left(\infty\right)
  \right)^{
    \gamma
    \left(
      \mu^{\left(n\right)}
    \right)
  }.
  \label{gamma_trans_CP2}
\end{align}

\noindent
Substituting \eqref{gamma_trans_CP2} into \eqref{another_CP2_before}, we have
\begin{flalign*}
  &
  \frac{
    \hbar
    \Gamma
    \left(
      1
      -
      \left(
        n-1
      \right)
      +
      \frac{1}{\hbar}
    \right)
  }{
    \left(
      n-1
    \right)!
    \Gamma
    \left(
      1
      +
      \frac{1}{\hbar}
    \right)
  }
  \sum_{\mu_{1},\nu_{1},\cdots,\mu_{n},\nu_{n}=1}^{2}
  \left(
    1
    -
    \frac{2}{n}
    \cdot
    n
    \left[
      \frac{\mu_{n}}{2}
    \right]
  \right)
  g_{\overline{\mu_{1}}\nu_{1}}
  \cdots
  g_{\overline{\mu_{n}}\nu_{n}}
  &
  \\
  &
  \mbox{}
  \quad
  \cdot
  \left(
    F_{d}\left(\infty\right)
  \right)^{
    \gamma
    \left(
      \nu^{\left(n\right)}
    \right)
  }
  T_{0}\left(\infty\right)
  \left(
    F_{r}\left(\infty\right)
  \right)^{
    \gamma
    \left(
      \mu^{\left(n\right)}
    \right)
  }
  &
  \\
  &=
  \frac{
    \hbar
    \Gamma
    \left(
      1
      -
      \left(
        n-1
      \right)
      +
      \frac{1}{\hbar}
    \right)
  }{
    \left(
      n-1
    \right)!
    \Gamma
    \left(
      1
      +
      \frac{1}{\hbar}
    \right)
  }
  \sum_{\mu_{1},\nu_{1},\cdots,\mu_{n},\nu_{n}=1}^{2}
  \left(
    1
    -
    \frac{2}{n}
    \left[
      \frac{\mu_{1}}{2}
    \right]
    -
    \cdots
    -
    \frac{2}{n}
    \left[
      \frac{\mu_{n}}{2}
    \right]
  \right)
  g_{\overline{\mu_{1}}\nu_{1}}
  \cdots
  g_{\overline{\mu_{n}}\nu_{n}}
  &
  \\
  &
  \mbox{}
  \quad
  \cdot
  \left(
    F_{d}\left(\infty\right)
  \right)^{
    \gamma
    \left(
      \nu^{\left(n\right)}
    \right)
  }
  T_{0}\left(\infty\right)
  \left(
    F_{r}\left(\infty\right)
  \right)^{
    \gamma
    \left(
      \mu^{\left(n\right)}
    \right)
  }
  \\
  &=
  \frac{
    \hbar
    \Gamma
    \left(
      1
      -
      \left(
        n-1
      \right)
      +
      \frac{1}{\hbar}
    \right)
  }{
    \left(
      n-1
    \right)!
    \Gamma
    \left(
      1
      +
      \frac{1}{\hbar}
    \right)
  }
  \sum_{\mu_{1},\nu_{1},\cdots,\mu_{n},\nu_{n}=1}^{2}
  \left(
    1
    -
    \frac{2}{n}
    \gamma
    \left(
      \mu^{\left(n\right)}
    \right)
  \right)
  g_{\overline{\mu_{1}}\nu_{1}}
  \cdots
  g_{\overline{\mu_{n}}\nu_{n}}
  &
  \\
  &
  \mbox{}
  \quad
  \cdot
  \left(
    F_{d}\left(\infty\right)
  \right)^{
    \gamma
    \left(
      \nu^{\left(n\right)}
    \right)
  }
  T_{0}\left(\infty\right)
  \left(
    F_{r}\left(\infty\right)
  \right)^{
    \gamma
    \left(
      \mu^{\left(n\right)}
    \right)
  }.
  &
\end{flalign*}

\noindent
Hence,
\begin{align}
  \hbar
  C_{n}^{'}\left(\infty\right)
  =&
  \frac{
    \hbar
    \Gamma
    \left(
      1
      -
      \left(
        n-1
      \right)
      +
      \frac{1}{\hbar}
    \right)
  }{
    \left(
      n-1
    \right)!
    \Gamma
    \left(
      1
      +
      \frac{1}{\hbar}
    \right)
  }
  \sum_{\mu_{1},\nu_{1},\cdots,\mu_{n},\nu_{n}=1}^{2}
  \left(
    1
    -
    \frac{2}{n}
    \gamma
    \left(
      \mu^{\left(n\right)}
    \right)
  \right)
  g_{\overline{\mu_{1}}\nu_{1}}
  \cdots
  g_{\overline{\mu_{n}}\nu_{n}}
  \nonumber
  \\
  &
  \cdot
  \left(
    F_{d}\left(\infty\right)
  \right)^{
    \gamma
    \left(
      \nu^{\left(n\right)}
    \right)
  }
  T_{0}\left(\infty\right)
  \left(
    F_{r}\left(\infty\right)
  \right)^{
    \gamma
    \left(
      \mu^{\left(n\right)}
    \right)
  }.
  \label{another_rhs_CP2}
\end{align}

\noindent
From \eqref{another_lhs_CP2_calc} and \eqref{another_rhs_CP2} it is shown 
that \eqref{formula_CP2} satisfies \eqref{minus_thm_N2_eq_rsz} in Subsection \ref{NR_minus} for the \(\mathbb{C}P^{2}\) case.

\section{Discussions}
\label{DQ_Discussion}
\ In this paper, we obtained the explicit star products with separation of variables for two-dimensional locally symmetric K\"{a}hler manifolds by solving the recurrence relations in \cite{HS1,HS2}. By using Theorem \ref{main_thm_N2}, we can construct a deformation quantization with separation of variables for any complex two-dimensional locally symmetric K\"{a}hler manifold. Furthermore, we gave the noncommutative \(\mathbb{C}^{2}\) and \(\mathbb{C}P^{2}\) as concrete examples. We verified these noncommutative \(\mathbb{C}^{2}\) and \(\mathbb{C}P^{2}\) coincide with the 
previous results by Karabegov \cite{Kara1} and Sako-Suzuki-Umetsu \cite{SSU1,SSU2}, respectively.
\bigskip
\\
\ In this section, we discuss what can be done in the future using the theorems obtained in this paper. We give an example of a two-dimensional 
symmetric K\"{a}hler manifold to which Theorem \ref{main_thm_N2} is applicable and for which deformation quantization has not yet been obtained. The concrete example of two dimensional locally symmetric K\"{a}hler manifold is a quadric surface 
\(
  Q^{2}
  \left(\mathbb{C}\right)
  =
  SO(4)/
  \left(
    SO(2)
    \times
    U(1)
  \right)
\). 
However, the deformation quantization for \(Q^{2}\left(\mathbb{C}\right)\) have not yet been constructed. 
It is known that the K\"{a}hler potential of \(Q^{2}\left(\mathbb{C}\right)\) is given by
\begin{align}
  \Phi_{Q^{2}\left(\mathbb{C}\right)}
  =
  \log
  \left(
    1
    +
    \sum_{k=1}^{2}
    |z^{k}|^{2}
    +
    \frac{1}{4}
    \sum_{k,l=1}^{2}
    (
      \overline{z}^{k}
    )^{2}
    (
      z^{l}
    )^{2}
  \right).
\label{kahler_p_Q2}
\end{align}
\noindent
See \cite{Ell,HN1,Nitta} for more detail. 
Since K\"{a}hler metric and curvature are expressed by using a K\"{a}hler potential, it is 
possible to explicitly construct the noncommutative 
\(Q^{2}\left(\mathbb{C}\right)\) by using Theorem \ref{main_thm_N2}, explicitly. Supersymmetric nonlinear sigma models whose target spaces are K\"{a}hler manifolds are studied by Higashijima-Nitta \cite{HN1,HN2,HN3}, Higashijima-Kimura-Nitta \cite{HKN} and Kondo-Takahashi \cite{KT}. Some models of them are defined on \(\mathbb{C}P^{N}\) and \(Q^{2}\left(\mathbb{C}\right)\) as the target spaces. 
It should be possible to extend the nonlinear sigma models to ones on noncommutative \(\mathbb{C}P^{N}\) and \(Q^{2}\left(\mathbb{C}\right)\), naively.
\bigskip
\\
\ We can expect to develop the physics on noncommutative locally symmetric K\"{a}hler manifolds from our result. For example, the constructions of some field theories and gauge theories on noncommutative ones are expected. 
Someone might think that in contrast to 
strict deformation quantization, 
formal deformation quantization is difficult to 
interpret some physics. However, it is possible to construct Fock representations from the formal deformation quantization in our cases. In fact, (twisted) Fock representations can be constructed by using a deformation quantization with separation of variables as discussed in Section \ref{DQ_Review}. See \cite{SSU1,SSU2,SU1,SU2,SU3} for more detailed discussions. Recalling that field theories can be 
described by using Fock representations, we can expect to propose the field 
theories on noncommutative K\"{a}hler manifold. 
In addition, it is expected to clarify the relationship between fuzzy manifolds and deformation quantization, since Fock representations can be interpreted by using a matrix representation. 
Furthermore, we can concretely construct gauge theories by using (twisted) Fock representations. 
They have been already studied by Maeda-Sako-Suzuki-Umetsu \cite{MSSU} and Sako-Suzuki-Umetsu \cite{SSU3} on noncommutative homogeneous K\"{a}hler manifolds. See \cite{Sako} for the review 
of these facts. For example, if the noncommutative \(Q^{2}\left(\mathbb{C}\right)\) are given by Theorem \ref{main_thm_N2}, we can propose gauge theories on noncommutative \(Q^{2}\left(\mathbb{C}\right)\).
\bigskip
\\
\ As described above, it is expected that various physical theories on complex two-dimensional noncommutative locally symmetric K\"{a}hler manifold can be obtained by using the Theorem \ref{main_thm_N2}. They are left for future 
work.

\section*{Acknowledgement}
\ A.S.\ was supported by JSPS KAKENHI Grant Number 21K03258. 
The authors are grateful to Yasufumi Nitta, Shunsuke Saito, and Yohei Ito for useful advice. The authors appreciate the referees of IJGMMP for their thoughtful feedback.

\appendix
\section{Some properties from K\"{a}hler geometry}
\label{app_kaehler}
\ We review some properties for K\"{a}hler manifolds that we use in this paper. See Kobayashi-Nomizu for more details \cite{KN}. Let \(M\) be a complex \(N\)-dimensional K\"{a}hler manifold, \(U\) be a holomorphic coordinate neighborhood of \(M\), and \(\nabla\) be the Levi-Civita connection on \(M\). For 
\(
  \partial_{B},
  \partial_{C}
  \in
  \Gamma\left(
    TM
    \left.
    \right|_{U}
  \right),
\) 
the Christoffel symbol \(\Gamma_{BC}^{A}\) is defined by 
\(
  \nabla_{\partial_{B}}
  \partial_{C}
  =
  \Gamma_{BC}^{A}
  \partial_{A},
\)
where 
\(
  A,B,C
  \in
  \left\{
    1,
    \cdots,
    N,
    \overline{1},
    \cdots,
    \overline{N}
  \right\},
  \
  \partial_{A}
  :=
  \frac{\partial}{\partial z^{A}},
\) 
and 
\(
  z^{\overline{i}}
  :=
  \overline{z}^{i}
\) 
for 
\(
  i
  \in
  \left\{
    1,
    \cdots,
    N
  \right\}.
\) 
In particular, \(\Gamma_{BC}^{A}\) is given by
\begin{align*}
  \Gamma_{BC}^{A}
  =
  \frac{1}{2}
  g^{AD}
  \left(
    \partial_{B}
    g_{DC}
    +
    \partial_{C}
    g_{DB}
    -
    \partial_{D}
    g_{BC}
  \right)
\end{align*}

\noindent
by using the components of the K\"{a}hler metric \(g\). Since \(M\) is K\"{a}hler, the non-trivial Christoffel symbols are \(\Gamma^{k}_{ij}\) and \(\Gamma^{\overline{k}}_{\overline{i}\,\overline{j}}\), where 
\(
  i,
  j,
  k
  \in
  \left\{
    1,
    \cdots,
    N
  \right\}.
\) 
Furthermore, the covariant derivative for \((k,l)\)-tensor field 
\(
  Y
  ^{
    A_{1}
    \cdots
    A_{k}
  }
  _{
    B_{1}
    \cdots
    B_{l}
  }
  \partial_{A_{1}}
  \otimes
  \cdots
  \otimes
  \partial_{A_{k}}
  \otimes
  dz^{B_{1}}
  \otimes
  \cdots
  \otimes
  dz^{B_{l}}
  \in
  \Gamma
  (
    (
      TM
    )^{
      \otimes
      k
    }
    \otimes
    (
      T^{^{\ast}}M
    )^{
      \otimes
      l
    }
  )
\) 
is given by using the Christoffel symbol as follows : 
\begin{align*}
  &
  \nabla_{\partial_{C}}
  \left(
    Y
    ^{
      A_{1}
      \cdots
      A_{k}
    }
    _{
      B_{1}
      \cdots
      B_{l}
    }
    \partial_{A_{1}}
    \otimes
    \cdots
    \otimes
    \partial_{A_{k}}
    \otimes
    dz^{B_{1}}
    \otimes
    \cdots
    \otimes
    dz^{B_{l}}
  \right)
  \\
  =&
  \left(
    \partial_{C}
    Y
    ^{
      \mu_{1}
      \cdots
      \mu_{k}
      \overline{\mu_{1}}
      \cdots
      \overline{\mu_{m}}
    }
    _{
      \nu_{1}
      \cdots
      \nu_{l}
      \overline{\nu_{1}}
      \cdots
      \overline{\nu_{n}}
    }
    +
    \sum_{q=1}^{k}
    \Gamma^{A_{q}}_{CD}
    Y
    ^{
      A_{1}
      \cdots
      A_{q-1}
      D
      A_{q+1}
      \cdots
      A_{k}
    }
    _{
      B_{1}
      \cdots
      B_{l}
    }
    -
    \sum_{q=1}^{l}
    \Gamma^{D}_{CB_{q}}
    Y
    ^{
      A_{1}
      \cdots
      A_{k}
    }
    _{
      B_{1}
      \cdots
      B_{q-1}
      D
      B_{q+1}
      \cdots
      B_{l}
    }
  \right)
  \\
  &
  \times
  \partial_{A_{1}}
  \otimes
  \cdots
  \otimes
  \partial_{A_{k}}
  \otimes
  dz^{B_{1}}
  \otimes
  \cdots
  \otimes
  dz^{B_{l}},
\end{align*}

\noindent
where 
\(
  A_{1},
  \cdots,
  A_{k},
  B_{1},
  \cdots,
  B_{l},
  C,
  D
  \in
  \left\{
    1,
    \cdots,
    N,
    \overline{1},
    \cdots,
    \overline{N}
  \right\}.
\)
\\
\ Next, we define the Riemann curvature tensor 
\(
  R^{\nabla}
  :
  \Gamma
  \left(
    TM
  \right)
  \times
  \Gamma
  \left(
    TM
  \right)
  \times
  \Gamma
  \left(
    TM
  \right)
  \to
  \Gamma
  \left(
    TM
  \right)
\) 
on \(M\) for the vector fields 
\(
  X,Y
  \in
  \Gamma
  \left(
    TM
  \right)
\) 
and its component 
\(
  R^{\nabla}
  \left(
    \partial_{A},
    \partial_{B}
  \right)
  \partial_{C}
\) 
by 
\begin{align*}
  R^{\nabla}
  \left(
  X,Y
  \right)
  &:=
  \nabla_{X}\nabla_{Y}
  -
  \nabla_{Y}\nabla_{X}
  -
  \nabla_{\left[X,Y\right]},
  \\
  R^{\nabla}
  \left(
    \partial_{A},
    \partial_{B}
  \right)
  \partial_{C}
  &:=
  R_{ABC}^{\quad \ \ \, D}
  \partial_{D},
\end{align*}

\noindent
respectively. 
Note that the notation \(R_{ABC}^{\quad \ \ \, D}\) used in this paper can be expressed by the relation
\begin{align}
  R_{ABC}^{\quad \ \ \, D}
  =
  \mathfrak{R}_{\ \,CAB}^{D}
  \label{HS_KN}
\end{align}

\noindent
using the notation
\begin{align*}
  R^{\nabla}
  \left(
    \partial_{A},
    \partial_{B}
  \right)
  \partial_{C}
  :=
  \mathfrak{R}_{\ \,CAB}^{D}
  \partial_{D}
\end{align*}

\noindent
by Kobayashi-Nomizu \cite{KN}. In this paper, we fix the position of the indices of the components of the Riemann curvature tensor in the above equation. The components \(R_{ABC}^{\quad \ \ \, D}\) are also given by
\begin{align*}
  R_{ABC}^{\quad \ \ \, D}
  =
  \partial_{A}
  \Gamma_{BC}^{D}
  -
  \partial_{B}
  \Gamma_{AC}^{D}
  +
  \Gamma_{BC}^{E}
  \Gamma_{AE}^{D}
  -
  \Gamma_{AC}^{E}
  \Gamma_{BE}^{D}
\end{align*}

\noindent
by using the Christoffel symbols, and their non-trivial ones are
\begin{align*}
  R_{i\overline{j}k}^{\quad \, l}
  &=
  -
  \partial_{\overline{j}}
  \Gamma_{ik}^{l},
  \\
  R_{i\bar{j}\bar{k}}^{\quad \, \bar{l}}
  &=
  \partial_{i}
  \Gamma_{\bar{j}\bar{k}}^{\bar{l}}.
\end{align*}

\noindent
For 
\(
  R_{ABCD}
  =
  g_{DE}
  R_{ABC}^{\quad \ \ \, E},
\) 
it is also confirmed that the non-trivial components are 
\(
  R_{i\overline{j}k\bar{l}}
\) 
and 
\(
  R_{i\bar{j}\bar{k}l},
\) 
which are given by
\begin{align*}
  R_{i\overline{j}k\overline{l}}
  &=
  -
  g_{\overline{l}p}
  \partial_{\overline{j}}
  \Gamma_{ik}^{p}
  =
  -
  \partial_{i}\partial_{\overline{j}}\partial_{k}\partial_{\overline{l}}\Phi
  +
  g^{p\overline{q}}
  \left(
    \partial_{i}\partial_{\overline{q}}\partial_{k}\Phi
  \right)
  \left(
    \partial_{p}\partial_{\overline{j}}\partial_{\overline{l}}\Phi
  \right),
  \\
  R_{i\overline{j}\overline{k}l}
  &=
  g_{l\overline{p}}
  \partial_{i}
  \Gamma_{\overline{j}\overline{k}}^{\overline{p}}
  =
  \partial_{\overline{j}}\partial_{i}\partial_{\overline{k}}\partial_{l}\Phi
  -
  g^{\overline{p}q}
  \left(
    \partial_{\overline{j}}\partial_{q}\partial_{\overline{k}}\Phi
  \right)
  \left(
    \partial_{\overline{p}}\partial_{i}\partial_{l}\Phi
  \right)
\end{align*}

\noindent
respectively, where \(\Phi\) is the K\"{a}hler potential, i.e. \(\Phi\) is a function on \(M\) such that 
\(
  g_{i\overline{j}}
  =
  \partial_{i}
  \partial_{\overline{j}}
  \Phi.
\)
\\
\ Here we consider 
\(
  R_{\bar{i} \ \ \, \bar{c}}^{\ \, \bar{k}\bar{l}}
  =
  g^{\overline{k}p}
  g^{\overline{l}q}
  R_{\overline{i}pq\overline{c}}
\) 
which often appears in this paper. We refer to 
\(
  R_{\bar{i} \ \ \, \bar{c}}^{\ \, \bar{k}\bar{l}}
\) 
simply as ``curvature". This curvature 
\(
  R_{\bar{i} \ \ \, \bar{c}}^{\ \, \bar{k}\bar{l}}
\) 
has the following symmetries 
concerning the indices 
\(
  \overline{i},
  \overline{c},
  \overline{k},
\)
and
\(
  \overline{l}
\)
 :
\begin{align}
  R_{\bar{i} \ \ \, \bar{c}}^{\ \, \bar{k}\bar{l}}
  =
  R_{\bar{c} \ \ \, \bar{i}}^{\ \, \bar{k}\bar{l}}
  =
  R_{\bar{i} \ \ \, \bar{c}}^{\ \, \bar{l}\bar{k}}
  =
  R_{\bar{c} \ \ \, \bar{i}}^{\ \, \bar{l}\bar{k}}.
  \label{curv_kaehler_prop}
\end{align}

\noindent
This property plays an important role in this paper.
\bigskip
\\
\ We now turn our attention to a locally symmetric K\"{a}hler manifold : the fact that \(M\) is locally symmetric is equivalent to the fact that 
\(
  \nabla_{\partial_{E}}
    R_{ABC}^{\quad \ \ \, D}
    =
    0,
\) 
for 
\(
  A,B,C,D,E
  \in
  \left\{
    1,
    \cdots,
    N,
    \overline{1}.
    \cdots,
    \overline{N}
  \right\}.
\)

\section{Calculations for both sides of \eqref{minus_thm_eq_ens} in Subsection \ref{NR_minus}}
\label{app_minus_n2}

\ We show that \eqref{minus_thm_eq_ens} for \(n=2\) holds in Subsection \ref{NR_minus}. In this appendix, we denote the detailed calculations for each component of both sides of \eqref{minus_thm_eq_ens} in Subsection \ref{NR_minus}. They can be enumerated as follows.
\vspace{3mm}
\\
\noindent
\underline{
  The left-hand side
}
\begin{align}
  \left(1,1\right)
  &=
  \hbar^{2}
  \left\{
    \left(
      2
      +
      \hbar
      R_{1 \ \ \, 1}^{\ \, 11}
    \right)
    \left(
      g_{\overline{1}1}
    \right)^{2}
    +
    2
    \hbar
    R_{1 \ \ \, 1}^{\ \, 21}
    g_{\overline{1}1}
    g_{\overline{1}2}
    +
    \hbar
    R_{1 \ \ \, 1}^{\ \, 22}
    \left(
      g_{\overline{1}2}
    \right)^{2}
  \right\}
  \nonumber
  \\
  &=
  \hbar^{2}
  \left\{
    2
    \left(
      g_{\overline{1}1}
    \right)^{2}
    +
    \hbar
    R_{1\bar{1}\bar{1}1}
  \right\},
  \label{lhs_first}
  \\
  \left(1,2\right)
  &=
  2\hbar^{2}
  \left\{
    \left(
	  2
      +
      \hbar
      R_{1 \ \ \, 1}^{\ \, 11}
    \right)
    g_{\overline{1}1}
    g_{\overline{2}1}
    +
    \hbar
    R_{1 \ \ \, 1}^{\ \, 21}
    \left(
      g_{\overline{1}1}
      g_{\overline{2}2}
      +
      g_{\overline{1}2}
      g_{\overline{2}1}
    \right)
    +
    \hbar
    R_{1 \ \ \, 1}^{\ \, 22}
    g_{\overline{1}2}
    g_{\overline{2}2}
  \right\}
  \nonumber
  \\
  &=
  2\hbar^{2}
  \left\{
    g_{\overline{1}1}
    g_{\overline{2}1}
    +
    \hbar
    R_{1\bar{2}\bar{1}1}
  \right\},
  \\
  \left(1,3\right)
  &=
  \hbar^{2}
  \left\{
    \left(
	  2
      +
      \hbar
      R_{1 \ \ \, 1}^{\ \, 11}
    \right)
    \left(
      g_{\overline{2}1}
    \right)^{2}
    +
    2
    \hbar
    R_{1 \ \ \, 1}^{\ \, 21}
    g_{\overline{2}1}
    g_{\overline{2}2}
    +
    \hbar
    R_{1 \ \ \, 1}^{\ \, 22}
    \left(
      g_{\overline{2}2}
    \right)^{2}
  \right\}
  \nonumber
  \\
  &=
  \hbar^{2}
  \left\{
    2
    \left(
      g_{\overline{2}1}
    \right)^{2}
    +
    \hbar
    R_{1\bar{2}\bar{2}1}
  \right\},
  \\
  \left(2,1\right)
  &=
  0,
  \\
  \left(2,2\right)
  &=
  0,
  \\
  \left(2,3\right)
  &=
  0,
  \\
  \left(3,1\right)
  &=
  -
  \hbar^{2}
  \left\{
    \hbar
    R_{2 \ \ \, 2}^{\ \, 11}
    \left(
      g_{\overline{1}1}
    \right)^{2}
    +
    2
    \hbar
    R_{2 \ \ \, 2}^{\ \, 21}
    g_{\overline{1}1}
    g_{\overline{1}2}
    +
    \left(
      2
      +
      \hbar
      R_{2 \ \ \, 2}^{\ \, 22}
    \right)
    \left(
      g_{\overline{1}2}
    \right)^{2}
  \right\}
  \nonumber
  \\
  &=
  -\hbar^{2}
  \left\{
    2
    \left(
      g_{\overline{1}2}
    \right)^{2}
    +
    \hbar
    R_{2\bar{1}\bar{1}2}
  \right\},
  \\
  \left(3,2\right)
  &=
  -
  2
  \hbar^{2}
  \left\{
    \hbar
    R_{2 \ \ \, 2}^{\ \, 11}
    g_{\overline{1}1}
    g_{\overline{2}1}
    +
    \hbar
    R_{2 \ \ \, 2}^{\ \, 21}
    \left(
      g_{\overline{1}1}
      g_{\overline{2}2}
      +
      g_{\overline{1}2}
      g_{\overline{2}1}
    \right)
    +
    \left(
      2
      +
      \hbar
      R_{2 \ \ \, 2}^{\ \, 22}
    \right)
    g_{\overline{1}2}
    g_{\overline{2}2}
  \right\}
  \nonumber
  \\
  &=
  -2\hbar^{2}
  \left\{
    2
    g_{\overline{1}2}
    g_{\overline{2}2}
    +
    \hbar
    R_{2\bar{2}\bar{1}2}
  \right\},
  \\
  \left(3,3\right)
  &=
  -
  \hbar^{2}
  \left\{
    \hbar
    R_{2 \ \ \, 2}^{\ \, 11}
    \left(
      g_{\overline{2}1}
    \right)^{2}
    +
    2
    \hbar
    R_{2 \ \ \, 2}^{\ \, 21}
    g_{\overline{2}1}
    g_{\overline{2}2}
    +
    \left(
      2
      +
      \hbar
      R_{2 \ \ \, 2}^{\ \, 22}
    \right)
    \left(
      g_{\overline{2}2}
    \right)^{2}
  \right\}
  \nonumber
  \\
  &=
  -\hbar^{2}
  \left\{
    2
    \left(
      g_{\overline{2}2}
    \right)^{2}
    +
    \hbar
    R_{2\bar{2}\bar{2}2}
  \right\}.
\end{align}
\noindent
\underline{
  The right-hand side
}
\begin{align}
  \left(1,1\right)
  &=
  \hbar^{2}
  \left\{
    \left(
      2
      +
      \hbar
      R_{\bar{1} \ \ \, \bar{1}}^{\ \, \bar{1}\bar{1}}
    \right)
    \left(
      g_{\overline{1}1}
    \right)^{2}
    +
    2
    \hbar
    R_{\bar{1} \ \ \, \bar{1}}^{\ \, \bar{2}\bar{1}}
    g_{\overline{1}1}
    g_{\overline{2}1}
    +
    \hbar
    R_{\bar{1} \ \ \, \bar{1}}^{\ \, \bar{2}\bar{2}}
    \left(
      g_{\overline{2}1}
    \right)^{2}
  \right\}
  \nonumber
  \\
  &=
  \hbar^{2}
  \left\{
    2
    \left(
      g_{\overline{1}1}
    \right)^{2}
    +
    \hbar
    R_{1\bar{1}\bar{1}1}
  \right\},
  \\
  \left(1,2\right)
  &=
  2\hbar^{2}
  \left\{
    \hbar
    R_{\bar{2} \ \ \, \bar{1}}^{\ \, \bar{1}\bar{1}}
    \left(
      g_{\overline{1}1}
    \right)^{2}
    +
    \left(
      2
      +
      2
      \hbar
      R_{\bar{2} \ \ \, \bar{1}}^{\ \, \bar{2}\bar{1}}
    \right)
    g_{\overline{1}1}
    g_{\overline{2}1}
    +
    \hbar
    R_{\bar{2} \ \ \, \bar{1}}^{\ \, \bar{2}\bar{2}}
    \left(
      g_{\overline{2}1}
    \right)^{2}
  \right\}
  \nonumber
  \\
  &=
  2\hbar^{2}
  \left\{
    g_{\overline{1}1}
    g_{\overline{2}1}
    +
    \hbar
    R_{1\bar{2}\bar{1}1}
  \right\},
  \\
  \left(1,3\right)
  &=
  \hbar^{2}
  \left\{
    \hbar
    R_{\bar{2} \ \ \, \bar{2}}^{\ \, \bar{1}\bar{1}}
    \left(
      g_{\overline{1}1}
    \right)^{2}
    +
    2
    \hbar
    R_{\bar{2} \ \ \, \bar{2}}^{\ \, \bar{2}\bar{1}}
    g_{\overline{1}1}
    g_{\overline{2}1}
    +
    \left(
      2
      +
      \hbar
      R_{\bar{2} \ \ \, \bar{2}}^{\ \, \bar{2}\bar{2}}
    \right)
    \left(
      g_{\overline{2}1}
    \right)^{2}
  \right\}
  \nonumber
  \\
  &=
  \hbar^{2}
  \left\{
    2
    \left(
      g_{\overline{2}1}
    \right)
    +
    \hbar
    R_{1\bar{2}\bar{2}1}
  \right\},
  \\
  \left(2,1\right)
  &=
  0,
  \\
  \left(2,2\right)
  &=
  0,
  \\\left(2,3\right)
  &=
  0,
  \\
  \left(3,1\right)
  &=
  -
  \hbar^{2}
  \left\{
    \left(
      2
      +
      \hbar
      R_{\bar{1} \ \ \, \bar{1}}^{\ \, \bar{1}\bar{1}}
    \right)
    \left(
      g_{\overline{1}2}
    \right)^{2}
    +
    2
    \hbar
    R_{\bar{1} \ \ \, \bar{1}}^{\ \, \bar{2}\bar{1}}
    g_{\overline{1}2}
    g_{\overline{2}2}
    +
    \hbar
    R_{\bar{1} \ \ \, \bar{1}}^{\ \, \bar{2}\bar{2}}
    \left(
      g_{\overline{2}2}
    \right)^{2}
  \right\}
  \nonumber
  \\
  &=
  -\hbar^{2}
  \left\{
    2
    \left(
      g_{\overline{1}2}
    \right)^{2}
    +
    \hbar
    R_{2\bar{1}\bar{1}2}
  \right\},
  \\
  \left(3,2\right)
  &=
  -2
  \hbar^{2}
  \left\{
    \hbar
    R_{\bar{2} \ \ \, \bar{1}}^{\ \, \bar{1}\bar{1}}
    \left(
      g_{\overline{1}2}
    \right)^{2}
    +
    \left(
      2
      +
      2
      \hbar
      R_{\bar{2} \ \ \, \bar{1}}^{\ \, \bar{2}\bar{1}}
    \right)
    g_{\overline{1}2}
    g_{\overline{2}2}
    +
    \hbar
    R_{\bar{2} \ \ \, \bar{1}}^{\ \, \bar{2}\bar{2}}
    \left(
      g_{\overline{2}2}
    \right)^{2}
  \right\}
  \nonumber
  \\
  &=
  -2\hbar^{2}
  \left\{
    2
    g_{\overline{1}2}
    g_{\overline{2}2}
    +
    \hbar
    R_{2\bar{2}\bar{1}2}
  \right\},
  \\
  \left(3,3\right)
  &=
  -
  \hbar^{2}
  \left\{
    \hbar
    R_{\bar{2} \ \ \, \bar{2}}^{\ \, \bar{1}\bar{1}}
    \left(
      g_{\overline{1}2}
    \right)^{2}
    +
    2
    \hbar
    R_{\bar{2} \ \ \, \bar{2}}^{\ \, \bar{2}\bar{1}}
    g_{\overline{1}2}
    g_{\overline{2}2}
    +
    \left(
      2
      +
      \hbar
      R_{\bar{2} \ \ \, \bar{2}}^{\ \, \bar{2}\bar{2}}
    \right)
    \left(
      g_{\overline{2}2}
    \right)^{2}
  \right\}
  \nonumber
  \\
  &=
  -\hbar^{2}
  \left\{
    2
    \left(
      g_{\overline{2}2}
    \right)^{2}
    +
    \hbar
    R_{2\bar{2}\bar{2}2}
  \right\}.
  \label{rhs_last}
\end{align}

\noindent
Note that we denote each \(\left(i,j\right)\) component of both sides simply as \(\left(i,j\right)\). Hence, the calculations \eqref{lhs_first}--\eqref{rhs_last} show that \eqref{minus_thm_eq_ens} in Subsection \ref{NR_minus} for \(n=2\) holds.

\section{
Hermiteness of 
\(
  T_{n}
\)
}
\label{prf_Tn_dagger}

In Subsection \ref{NR_minus}, we use the fact that \(T_{n}\) is Hermitian conjugate, i.e. 
\(
  T_{n}
  =
  T_{n}^{\dagger}.
\) 
So we shall derive this property of \(T_{n}\) in this appendix.

\begin{lemma}[Sako-Umetsu\cite{SU1}]
  The coefficients \(T_{\overrightarrow{\alpha_{n}},\,\overrightarrow{\beta_{n}^{\ast}}}^{n}\) satisfy
  \begin{align*}
    T_{\overrightarrow{\alpha_{n}},\,\overrightarrow{\beta_{n}^{\ast}}}^{n}
    =
    \overline{
      T_{\overrightarrow{\beta_{n}},\,\overrightarrow{\alpha_{n}^{\ast}}}^{n}
    }
  \end{align*}
  
  \noindent
  or equivalently
  \begin{align}
    T_{n}
    =
    T_{n}^{\dagger},
    \label{star_conj_mtx}
  \end{align}
\end{lemma}

\begin{prf*}
From Proposition 3.1. in \cite{SU1},
\begin{align}
  \overline{
    f
    \ast
    g
  }
  =
  \overline{g}
  \ast
  \overline{f}.
  \label{star_conj}
\end{align}

\noindent
Substituting \eqref{star_expand} into the left-hand side of \eqref{star_conj}, we have
\begin{align*}
  \overline{
    f
    \ast
    g
  }
  =&
  \sum_{n=0}^{\infty}
  \sum_{
    \overrightarrow{\alpha_{n}},
    \overrightarrow{\beta_{n}^{\ast}}
  }
  \overline{
    T_{\overrightarrow{\alpha_{n}},\overrightarrow{\beta_{n}^{\ast}}}^{n}
  }
  \left\{
    \left(
      D^{\overline{1}}
    \right)^{\alpha_{1}^{n}}
    \cdots
    \left(
      D^{\overline{N}}
    \right)^{\alpha_{N}^{n}}
    \overline{f}
  \right\}
  \left\{
    \left(
      D^{1}
    \right)^{\beta_{1}^{n}}
    \cdots
    \left(
      D^{N}
    \right)^{\beta_{N}^{n}}
    \overline{g}
  \right\}
  \\
  =&
  \sum_{n=0}^{\infty}
  \sum_{
    \overrightarrow{\alpha_{n}},
    \overrightarrow{\beta_{n}^{\ast}}
  }
  \overline{
    T_{\overrightarrow{\beta_{n}},\overrightarrow{\alpha_{n}^{\ast}}}^{n}
  }
  \left(
    D^{\overrightarrow{\alpha_{n}}}\overline{g}
  \right)
  \left(
    D^{\overrightarrow{\beta_{n}^{\ast}}}\overline{f}
  \right).
\end{align*}

\noindent
Here we rewrote 
the dummy indices in the last equation. On the other hand, since 
\begin{align*}
  \overline{g}
  \ast
  \overline{f}
  =
  \sum_{n=0}^{\infty}
  \sum_{
    \overrightarrow{\alpha_{n}},
    \overrightarrow{\beta_{n}^{\ast}}
  }
  T_{\overrightarrow{\alpha_{n}},\overrightarrow{\beta_{n}^{\ast}}}^{n}
  \left(
    D^{\overrightarrow{\alpha_{n}}}\overline{g}
  \right)
  \left(
    D^{\overrightarrow{\beta_{n}^{\ast}}}\overline{f}
  \right)
\end{align*}

\noindent
by the assumption for the star product with separation of variables, we obtain 
\begin{align*}
  T_{\overrightarrow{\alpha_{n}},\,\overrightarrow{\beta_{n}^{\ast}}}^{n}
  =
  \overline{
    T_{\overrightarrow{\beta_{n}},\,\overrightarrow{\alpha_{n}^{\ast}}}^{n}
  }
\end{align*}

\noindent
by comparing the coefficients on both sides of \eqref{star_conj}. The above equation can be expressed as an equivalent equation 
\(
  T_{n}
  =
  T_{n}^{\dagger}.
  _{\blacksquare}
\)
\end{prf*}


\begin{thebibliography}{99}
    \bibitem{ABIY} G.~Alexanian, A.~P.~Balachandran, G.~Immirzi and B.~Ydri, ``Fuzzy \(\mathbb{C}P^{2}\)", J.~Geom.~Phys.~\textbf{42}(1-2), (2002), 28--53.
    \bibitem{APS} G.~Alexanian, A.~Pinzul and A.~Stern, ``Generalized coherent state approach to star products and applications to the fuzzy sphere", Nucl.~Phys.~B~\textbf{600}(3), (2001), 531--547.
    \bibitem{BDMLO} A.~P.~Balachandran, B.~P.~Dolan, J.~H.~Lee, X.~Martin and D.~O'Connor, ``Fuzzy complex projective spaces and their star products", 
J.~Geom.~Phys.~\textbf{43}, (2002), 184--204.
    \bibitem{BFFLS} F.~Bayen, M.~Flato, C.~Fronsdal, A.~Lichnerowicz and D.~Sternheimer, ``Deformation theory and quantization. I. Deformations of symplectic structures", Ann.~Phys.~\textbf{111}, (1978), 61--111.
    \bibitem{Bere1} F.~A.~Berezin, ``Quantization", MATH~USSR--Izv, \textbf{8}(5), (1974), 1109-–1165.
    \bibitem{Bere2} F.~A.~Berezin, ``General concept of quantization", Comm.~Math.~Phys.~\textbf{40}(2), (1975), 153--174.
	\bibitem{BBEW} M.~Bordemann, M.~Brischle, C.~Emmrich and S.~Waldmann, ``Phase Space Reduction for Star-Products: An Explicit Construction for \(\mathbb{C}P^{n}\)", Lett.~Math.~Phys.~\textbf{36}, (1995), 357--371.
    \bibitem{BMS} M.~Bordemann, E.~Meinrenken and M.~Schlichenmaier, ``Toeplitz quantization of K\"{a}hler manifolds and \(gl(N), N\to\infty\) limits", Comm. Math. Phys.~\textbf{165}(2), (1994), 281--296.
    \bibitem{BW1} M.~Bordemann and S.~Waldmann, ``A Fedosov star product of the Wick type for K\"{a}hler manifolds", Lett.~Math.~Phys.~\textbf{41}(3), (1997), 243--253.
    \bibitem{BW2} M.~Bordemann and S.~Waldmann, ``Formal GNS construction and states in deformation quantization", Comm.~Math.~Phys.~\textbf{195}(3), (1998), 549--583.
    \bibitem{CGR1} M.~Cahen, S.~Gutt, and J.~Rawnsley, ``Quantization of K\"{a}hler manifolds. I", J.~Geom.~Phys.~\textbf{7}(1), (1991), 45--62.
    \bibitem{CGR2} M.~Cahen, S.~Gutt, and J.~Rawnsley, ``Quantization of K\"{a}hler manifolds. II", Trans.~Amer.~Math.~Soc.~\textbf{337}, (1993), 73--98.
    \bibitem{CGR3} M.~Cahen, S.~Gutt, and J.~Rawnsley, ``Quantization of K\"{a}hler manifolds. III", Lett.~Math.~Phys.~\textbf{30}, (1994), 291--305.
    \bibitem{CGR4} M.~Cahen, S.~Gutt, and J.~Rawnsley, ``Quantization of K\"{a}hler manifolds. IV", Lett.~Math.~Phys.~\textbf{34}, (1995), 159--168.
    \bibitem{CSW} U.~Carow-Watamura, H.~Steinacker and S.~Watamura, ``Monopole bundles over fuzzy complex projective spaces", J.~Geom.~Phys.~\textbf{54}(4), (2005), 373--399.
    \bibitem{DL} M.~de Wilde and P.~B.~A.~Lecomte, ``Existence of star-products and of formal deformations of the Poisson Lie algebra of arbitrary symplectic manifolds", Lett.~Math.~Phys.~\textbf{7}, (1983), 487--496.
    \bibitem{Dirac} P.~A.~M.~Dirac, ``The fundamental equations of quantum mechanics", Proc.~R.~Soc.~Lond.~A~\textbf{109}(752), (1925), 642--653.
    \bibitem{Ell} U.~Ellwanger, ``Supersymmetric \(\sigma\)-models in four dimensions as quantum theories", Nucl.~Phys.~B, \textbf{281}(3-4), (1987), 489--508.
    \bibitem{Fedo} B.~V.~Fedosov, ``A simple geometrical construction of deformation quantization", J.~Differ.~Geom.~\textbf{40(2)}, (1994), 213--238.
    \bibitem{GS} H.~Grosse and A.~Strohmaier, ``Towards a nonperturbative covariant regularization in 4D quantum field theory", Lett.~Math.~Phys.~\textbf{48}, (1999), 163--179.
    \bibitem{GUT} S.~Gutt, ``Deformation quantization of Poisson manifolds", Geometry and Topology Monographs~\textbf{17}, (2011), 171--220.
    \bibitem{HS1} K.~Hara and A.~Sako, ``Noncommutative Deformations of Locally Symmetric K\"{a}hler manifolds", J.~Geom.~Phys.~\textbf{114}, (2017), 554--569.
    \bibitem{HS2} K.~Hara and A.~Sako, ``Quantization of Locally Symmetric K\"ahler manifolds", Geometry Integrability and Quantization~\textbf{19}, (2018), 122--131.
    \bibitem{HNT} K.~Hayasaka, R.~Nakayama and Y.~Takaya, ``A New Noncommutative Product on the Fuzzy Two-Sphere Corresponding to the Unitary Representation of \(SU(2)\) and the Seiberg-Witten Map", Phys.~Lett.~B~\textbf{553(1--2)}, (2003), 109--118.
    \bibitem{HKN} K.~Higashijima, T.~Kimura and M.~Nitta, ``Gauge theoretical construction of non-compact Calabi–Yau manifolds", Ann.~Phys.~(N.Y.)~\textbf{296}(2), (2002), 347--370.
    \bibitem{HN1} K.~Higashijima and M.~Nitta, ``Supersymmetric nonlinear sigma models as gauge theories", Prog.~Theor.~Phys.~\textbf{103}(3), (2000), 635--663.
    \bibitem{HN2} K.~Higashijima and M.~Nitta, ``Quantum equivalence of auxiliary field methods in supersymmetric theories", Prog.~Theor.~Phys.~\textbf{103}(4), (2000), 833-846.
    \bibitem{HN3} K.~Higashijima and M.~Nitta, ``Supersymmetric nonlinear sigma models", in Proc.~Int.~Symp.~Quantum~Chromodynamics~and~Color~Confinement, eds. H.~Suganuma, M.~Fukushima and H.~Toki, (World Scientific Publishing Co Ptd Ltd, Osaka, 2000), 279--286pp.
    \bibitem{Hop} J.~Hoppe, ``Quantum theory of a massless relativistic surface and a two-dimensional bound state problem", Soryushiron~Kenkyu~Electronics~\textbf{80}(3), (1989), 145--202.
    \bibitem{Kara1} A.~V.~Karabegov, ``On the Deformation Quantization, on a K\"{a}hler Manifold, Associated with a Berezin Quantization", Funct.~Anal.~its~Appl.~\textbf{30(2)}, (1996), 142--144.
    \bibitem{Kara2} A.~V.~Karabegov, ``Deformation quantization with separation of variables on a K\"{a}hler manifold", Comm.~Math.~Phys.~\textbf{180}, (1996), 745--755.
    \bibitem{KS1} A.~V.~Karabegov and M.~Schlichenmaier, ``Identification of Berezin-Toeplitz deformation quantization", Crelle's Journal~\textbf{2001}, (2000), 49--76.
    \bibitem{KN} S.~Kobayashi and K.~Nomizu, ``Foundations of Differential Geometry II", New York:~John~Wiley\&Sons, (1969), 470 pp.
    \bibitem{KT} A.~Kondo and T.~Takahashi, ``Supersymmetric nonlinear sigma models as anomalous gauge theories", Phys.~Rev.~D~\textbf{102}(2), 025014, (2020), 1--13.
    \bibitem{Kont1} M.~Kontsevich, ``Deformation quantization of Poisson manifolds", Lett.~Math.~Phys.~\textbf{66}, (2003), 157--216.
    \bibitem{Kos} B.~Kostant, ``Quantization and unitary representations", Berlin: Springer, Lectures in Modern Analysis and Applications III. Lecture Notes in Mathematics \textbf{170}, (1970), 87--208pp.
    \bibitem{Mad1} J.~Madore, ``The Fuzzy sphere", Class~Quantum~Gravity~\textbf{9}, (1992), 69--87.
    \bibitem{Mad2} J.~Madore, ``An Introduction to Noncommutative Differential Geometry and its Physical Applications", (2nd ed., London Mathematical Society Lecture Note Series). Cambridge:~Cambridge~University~Press, (1999), 371pp.
    \bibitem{MSSU} Y.~Maeda, A.~Sako, T.~Suzuki and H.~Umetsu, ``Gauge theories in noncommutative homogeneous K\"{a}hler manifolds", J.~Math.~Phys.~\textbf{55}(9), 092301, (2014), 1--27.
    \bibitem{Mor1} C.~Moreno, ``\(\ast\)-products on some K\"{a}hler manifolds", Lett.~Math.~Phys.~\textbf{11}, (1986),361--372.
    \bibitem{Mor2} C.~Moreno, ``Invariant star products and representations of compact semisimple Lie groups", Lett.~Math.~Phys.~\textbf{12}, (1986), 217--229.
    \bibitem{Moya} J.~E.~Moyal, ``Quantum Mechanics as a Statistical Theory", Proceedings of the Cambridge Philosophical Society~\textbf{45(1)}, (1949), 99--124.
    \bibitem{Nitta} M.~Nitta, ``K\"{a}hler Potential for Global Symmetry Breaking in Supersymmetric Theories", hep-th/9903174, (1999).
    \bibitem{OMY1} H.~Omori, Y.~Maeda and A.~Yoshioka, ``Weyl Manifolds and Deformation Quantization", Adv.~Math.~\textbf{85}, (1991), 224--255.
    \bibitem{OMY2} H.~Omori, Y.~Maeda and A.~Yoshioka, ``Non-Commutative Complex Projective Space", Adv.~Math.~\textbf{22}, (1993), 133--152.
    \bibitem{OMMY} H.~Omori, Y.~Maeda, N.~Miyazaki and A.~Yoshioka, ``Poincar\'{e}--Cartan Class and Deformation Quantization of K\"{a}hler Manifolds", Comm.~Math.~Phys.~textbf{194}, (1998), 207--230.
    \bibitem{OS_Varna} T.~Okuda and A.~Sako, ``Deformation quantization with separation of variables for complex two-dimensional locally symmetric K\"{a}hler manifold", J.~Geom.~Symmetry~Phys.~\textbf{64}, (2022), 39--49.
    \bibitem{Pere} A.~M.~Perelomov, ``Generalized coherent states and their applications", Berlin: Springer, (1986), 320pp.
    \bibitem{Rawn} J.~H.~Rawnsley, ``Coherent states and K\"{a}hler manifolds", Q.~J.~Math.~\textbf{28}(4), (1977), 403--415.
    \bibitem{Rief1} M.~Rieffel, ``Deformation quantization of Heisenberg manifolds", Comm.~Math.~Phys.~\textbf{122}(4), (1989), 531--562.
    \bibitem{Rief2} M.~Rieffel, ``Deformation quantization for actions of \(\mathbb{R}^{d}\)", Mem.~Amer.~Math.~Soc.~\textbf{106}(506), 1993, 93pp.
    \bibitem{Rief3} M.~Rieffel, ``Quantization and C*-algebras", Contemporary Mathematics~\textbf{167}, (1994), 67--97.
    \bibitem{Sako} A.~Sako, ``A Recipe To Construct A Gauge Theory On A Noncommutative K\"{a}hler Manifold", Noncommutative Geometry and Physics \textbf{4}, (2017), 361--404.
    \bibitem{SSU1} A.~Sako, T.~Suzuki and H.~Umetsu, ``Explicit Formulae for Noncommutative Deformations of \(\mathbb{C}P^{N}\) and \(\mathbb{C}H^{N}\)", J.~Math.~Phys.~\textbf{53}(7), 073502, (2012), 1--16.
    \bibitem{SSU2} A.~Sako, T.~Suzuki and H.~Umetsu, ``Noncommutative \(\mathbb{C}P^{N}\) and \(\mathbb{C}H^{N}\) and their physics", J.~Phys.~Conf.~Ser.~\textbf{442}, 012052, (2013), 1--10.
    \bibitem{SSU3} A.~Sako, T.~Suzuki and H.~Umetsu, ``Gauge theories on noncommutative \(\mathbb{C}P^{N}\) and Bogomol'nyi-Prasad-Sommerfield-like equations", J.~Math.~Phys.~\textbf{56}, 113506, (2015), 1--22.
    \bibitem{SU1} A.~Sako and H.~Umetsu, ``Twisted Fock Representations of Noncommutative K\"{a}hler Manifolds", J.~Math.~Phys.~\textbf{57}(9), 093501, (2016), 1--30.
    \bibitem{SU2} A.~Sako and H.~Umetsu, ``Deformation Quantization of K\"{a}hler Manifolds and Their Twisted Fock Representation", Geometry, Integrability and Quantization~\textbf{18}, (2016), 225--240.
    \bibitem{SU3} A.~Sako and H.~Umetsu, ``Fock Representations and Deformation Quantization of K\"ahler Manifolds", Advances in Applied Clifford Algebras~\textbf{27}, (2017), 2769--2794.
    \bibitem{Schli1} M.~Schlichenmaier, ``Deformation quantization of compact K\"{a}hler manifolds via Berezin-Toeplitz operators", in Group 21: physical applications and mathematical aspects of geometry, groups and algebras; proceeding of the XXI International Colloquium on Group Theoretical Methods in Physics, July 15--20, 1996, Goslar, (1996), 396--400.
    \bibitem{Schli2} M.~Schlichenmaier, ``Berezin-Toeplitz quantization of compact K\"{a}hler manifolds", in Quantization, coherent states and Poisson structures: proceedings of the XIV Workshop on Geometric Methods in Physics, Bia\l
    owie\.{z}a, July 9--15, 1995, (1998), 101--115.
    \bibitem{Schli3} M.~Schlichenmaier, ``Berezin-Toeplitz quantization and Berezin symbols for arbitrary compact K\"{a}hler manifolds", in Coherent states, quantization and gravity : proceedings of the XVII Workshop on Geometric Methods in Physics, Bia\l
    owie\.{z}a, July 3--9, 1998, (2001), 45--46.
    \bibitem{Schli4} M.~Schlichenmaier, ``Deformation quantization of compact K\"{a}hler manifolds by Berezin-Toeplitz quantization", Conf.~Mosh\'{e}~Flato~1999, (Springer, Dordrecht), (2000), 289--306.
    \bibitem{Schli5} M.~Schlichenmaier, ``Berezin-Toeplitz quantization and Berezin transform", Long Time Behaviour of Classical and Quantum Systems, (World Scientific), (2001), 271--287.
    \bibitem{Schli6} M.~Schlichenmaier, ``Berezin-Toeplitz quantization for compact K\"{a}hler manifolds. A review of results", Adv.~Math.~Phys.~\textbf{2010}, (2010), 643--680.
    \bibitem{Schli7} M.~Schlichenmaier, ``Berezin–Toeplitz quantization and star products for compact K\"{a}hler manifolds", Contemp.~Math.~\textbf{583}, (2012), 257--294.
    \bibitem{Sou} J.~Souriau, ``Structure des syst\'{e}mes dynamiques", Paris: Dunod, (1970), 414pp.
    \bibitem{Voro} A.Voros, ``Wentzel-Kramers-Brillouin method in the Bargmann representation", Phys.~Rev.~A~\textbf{40} (1989), 6814--6825.
\end{thebibliography}
\end{document}